\documentclass[final,3p,11pt, times]{article}
\usepackage{authblk}
\usepackage{amsmath}
\usepackage{amssymb,latexsym}
\usepackage[german,english]{babel}
\usepackage{url}
\usepackage{graphicx}
\usepackage{gastex}
\usepackage{longtable}
\usepackage{lscape}
\usepackage{tabularx}
\usepackage{multicol}
\usepackage{verbatim}
\usepackage{multirow}
\usepackage{graphicx}
\usepackage{epsfig,graphics,graphicx}
\usepackage{geometry}
\usepackage[usenames]{xcolor}
\usepackage{authblk}
\usepackage{verbatim}
\usepackage{multirow}
\usepackage{ragged2e}
\usepackage{xcolor}
\usepackage[T1]{fontenc}
\usepackage{authblk}
\usepackage{lipsum}
\usepackage{hyperref}
\hypersetup{colorlinks=pink,linkcolor=blue,filecolor=red,urlcolor=blue,}
\usepackage{tikz}

\usepackage{lineno}

\hbadness=\maxdimen
\newtheorem{eg}{{\bf Example}}[subsection]
\newtheorem{rem}{{\bf Remark}}
\newtheorem{deff}{Definition}[subsection]
\newtheorem{lem}{{\bf Lemma}}[subsection]
\newtheorem{thm}{{\bf Theorem}}[section]
\newtheorem{prop}{{\bf Proposition}}[section]
\newtheorem{cor}{{\bf Corollary}}[section]
\newtheorem{claim}{{\bf Claim}}[section]

\begin{document}
	
	\title{{Some doubly semi-equivelar maps on the plane and the torus}}
	
	\author{Yogendra Singh, Anand Kumar Tiwari}

	\affil{\small Department of Applied Science, Indian Institute of Information Technology, Allahabad 211015, India. E-mail: $\{$rss2018503,anand$\}$@iiita.ac.in}
	
	\maketitle
	
	\hrule
	
	\begin{abstract} 
		A vertex $v$ in a map $M$ has the face-sequence $(p_1 ^{n_1}. \ldots. p_k^{n_k})$, if there are consecutive $n_i$ numbers of $p_i$-gons incident at $v$ in the given cyclic order, for $1 \leq i \leq k$.  A map $M$ is called a semi-equivelar map if each of its vertex has same face-sequence. Doubly semi-equivelar maps are a generalization of semi-equivelar maps which have precisely 2 distinct face-sequences. In this article, we enumerate the types of doubly semi-equivelar maps on the plane and torus which have combinatorial curvature 0. Further, we present classification of doubly semi-equivelar maps on the torus and illustrate this classification for those doubly semi-equivelar maps which comprise of face-sequence pairs $\{(3^6), (3^3.4^2)\}$ and $\{(3^3.4^2), (4^4)\}$.

\end{abstract}

\smallskip



\hrule

\section{Introduction}

An embedding of a simple graph $G$ into a surface $S$ is called a 2-cell embedding if the closure of each component of $F \setminus G$ is a $p$-gonal 2-disk, denoted as $D_p$ $(p \geq 3)$. These disks are called the faces of the embedding. A map $M$ is a 2-cell embedding on a surface $S$ such that the intersection of any two faces is either empty, a vertex, or an edge. Let $M_1$ and $M_2$ be maps with the vertex set $V(M_1)$ and $V(M_2)$ respectively. Then $M_1$ is isomorphic to $M_2$ ($M_1 \cong M_2$) if there is a bijective map $f : V(M_1) \to V(M_2)$ such that $v_1$-$v_2$ is an edge in $M_1$ if and only if $f(v_1)$-$f(v_2)$ is an edge in $M_2$, $[v_1 v_2 \ldots v_n]$ is a face in $M_1$ if and only $[f(v_1) f(v_2) \ldots f(v_n)]$ is a face in $M_2$, and $f$ preserves the incidence of edges and incidence of faces.


\smallskip

Let $v$ be a vertex in a map $M$. The face-sequence of $v$, denoted as $f_{seq}(v)$, is $f_{seq}(v)=(p_1 ^{n_1}. \ldots. p_k^{n_k})$, if there are $n_i$ numbers of $p_i$-gonal faces incident at $v$ in the given cyclic order, for $1 \leq i \leq k$. The combinatorial curvature of the vertex $v$ is defined as $\phi(v) = 1 - (\sum_{i=1}^{k} n_i)/2 + \sum_{i=1}^{k} (n_i/p_i)$. We say that a map $M$ has the combinatorial curvature $k$ if $\phi(v) = k$ for each vertex $v$ in $M$. Some maps of positive curvature are discussed in \cite{rbk,zhang}.

\smallskip

The union of all the faces incidents at a vertex $v$ in a map $M$ form a 2-disk $D_p$ with the boundary cycle $C_p$ (cycle of length $p$). The cycle $C_p$ is called the link of $v$, and is denoted as ${\rm lk}(v)$. Let $v$ be a vertex with ${\rm lk}(v) = C_p(v_1, \ldots, v_p)$. The face-sequence of ${\rm lk}(v)$ is  $(f_1^{n_1}. \ldots. f_{l}^ {n_l})$,  if the consecutive $n_i$ numbers of vertices in ${\rm lk}(v)$ have the face-sequence $f_i$ in the given cyclic order, for $1 \leq i \leq l$.

\smallskip

A map is called a semi-equivelar map of type $[p_1 ^{n_1}. \ldots. p_k^{n_k}]$ if the face-sequence of each vertex is $(p_1 ^{n_1}. \ldots. p_k^{n_k})$. Maps on the surfaces of Platonic solids and Archimedean solids are examples of semi-equivelar maps on the 2-sphere. The well known eleven types Archimedean tilings \cite{gru(1987)} on the plane provide semi-equivelar maps of types $[3^6]$, $[3^4.6]$, $[3^3.4^2]$, $[3^2.4.3.4]$, $[3.4.6.4]$, $[3.6.3.6]$, $[3.12^2]$, $[4^4]$, $[4.6.12]$, $[4.8^2]$, and $[6^3]$. Further, by \cite{CoMo(1980)}, we know that there is a semi-equivelar map on the plane $\mathbb{R}^2$ of type $[p^q]$ for every $p, q >0$ such that $1/p + 1/q <1/2$.  Thus we have: 

\begin{prop}\label{p1}
	There are infinitely many semi-equivelar maps on the plane.
\end{prop}

Since the plane $\mathbb{R}^2$ is the universal cover of the torus, the semi-equivelar maps on the plane suggest to explore the same types maps on the torus. Altshuler \cite{alt(1973)} constructed semi-equivelar maps of types $[3^6]$ and $[6^3]$ on the torus. Kurth \cite{kur(1986)} enumerated these types semi-equivelar maps including $[4^4]$ on the torus. Brehm and Kuhnel \cite{bk(2008)} classified these three types semi-equivelar maps on the torus using a different approach. Tiwari and Upadhyay \cite{tu(2017)} classified semi-equivelar maps on the torus of types $[3^4.6]$, $[3^3.4^2]$, $[3^2.4.3.4]$, $[3.4.6.4]$, $[3.6.3.6]$, $[3.12^2]$, $[4.6.12]$, $[4.8^2]$, and $[6^3]$. Recently, Datta and Maity \cite{dm(2017)} showed that only those types semi-equivelar maps exist on the torus which have combinatorial curvature 0: In other words, they proved: 

\begin{prop}\label{p2}
	If $X$ denote the type of a semi-equivelar map on the torus then $X$ is $[3^6]$, $[3^4.6]$, $[3^3.4^2]$, $[3^2.4.3.4]$, $[3.4.6.4]$, $[3.6.3.6]$, $[3.12^2]$, $[4^4]$, $[4.6.12]$, $[4.8^2]$ or $[6^3]$.
\end{prop}  

\smallskip 

Doubly semi-equivelar maps are a generalization of semi-equivelar maps which have two distinct face-sequences. More precisely, a map $M$ with face-sequences $f_1$ and $f_2$ is called a doubly semi-equivelar map (in short DSEM) of type $[f_{1}^{F} :f_2^{G}]$, if every vertex with the face sequence $f_1$ has the face-sequence of its link as $F =(f_{11}^{l_1}. \ldots.f_{1k_1}^{l_{k_1}})$ and every vertex with face sequence $f_2$ has the face-sequence of its link as $G=(f_{21}^{m_1}.\ldots.f_{2k_2}^{m_{k_2}})$. Here $f_{ij}$ is either $f_1$ or $f_2$ for each $i,j$. For $f_1 = f_2 = f$, the map $M$ turns into a semi-equivelar map of type $f$. Note that, if we stack (a process of splitting a $p$-gonal face into $p$ numbers of triangles by introducing a new vertex inside the face and joining it to all the vertices of the face by an edge) each face of a semi-equivelar map of type $[p^q]$, we get a DSEM of type $[{3^p}^{F}: {3^{2q}}^{G}]$ for some $F$ and $G$. Other than these, the twenty 2-uniform tilings \cite{gru(1987)} of the plane give doubly semi-equivelar maps on the plane. Thus we have: 

\begin{prop} \label{p3}
	There are infinitely many doubly semi-equivelar maps on the plane.
\end{prop}

Here we are interested in the case of doubly semi-equivelar maps on the plane which have combinatorial curvature 0. We show:


\begin{thm} \label{t1}
	There are exactly 22 types of DSEMs, say  $T_1, T_2, \ldots, T_{22}$ on the plane $\mathbb{R}^2$ with combinatorial curvature 0. The types $T_i$'s, $1 \leq i \leq 22$, are given explicitly in the proof. For each these types, except $T_{22}$, there is a unique DSEM $N_i$ on the plane, as shown in example Sec. \ref{s3}. For the type $T_{22}$, we get infinitely many non-isomorphic DSEMs on the plane.

\end{thm}

Now, if we take the quotient of these maps (given in example Sec. \ref{s3}) by suitable translation groups, we get possible types of doubly semi-equivelar map of curvature 0 on the torus (see for example the shaded region in Figure 3, given in Sec. \ref{s3}).  Thus, we have:

\begin{cor} \label{co1}
	There are exactly 22 types of DSEMs with combinatorial curvature 0 on the torus.
\end{cor} \hfill $\Box$

We describe classification of DSEMs on the torus which consist the pairs of face-sequences $\{(3^6), (3^3.4^2)\}$ and $\{(3^3.4^2), (4^4)\}$. Further, we enumerate such DSEMs for a small number of vertices. Thus we show:

\begin{thm} \label{t2}
	Let $M_i$ be a DSEM of type $T_i$, $i \in \{1,2,3,4\}$, with $n$ vertices. Then $M_i$ can be classified on the torus up to isomorphism.
\end{thm}


\begin{cor} \label{co2}
	Let $M_i$ be a DSEM of type $T_i$, $1 \leq i \leq 4$ on the torus with $|V(M_i)|$ vertices. Then (a) there are 14 DSEMs for $i=1$ and $|V(M_i)| \leq 18$ (b) there are 12 DSEMs for $i=2$ and $|V(M_i)| \leq 14$ (c) there are 12 DSEMs for $i=3$ and $|V(M_i)| \leq 18$, and (d) there are 12 DSEMs for $i=4$ and $|V(M_i)| \leq 24$. See these list in Tables \ref{table:1} - \ref{table:4}. 
\end{cor} 


\section{Definitions and Notations}\label{s2}
For the graph theory related terms, we refer to Bondy and Murty \cite{bondy(2008)}. 
Let $G$ be a graph, with the vertex set $V(G)$ and the edge set $E(G)$, embedded on a surface $S$. The notation $u$-$v$ denotes an edge joining $u, v \in V(G)$. The notation $P(u_1, \ldots, u_n)$ denotes a path $u_1$-$u_2$- $\cdots$ -$u_n$, where the vertices $u_1$ and $u_2$ are called the boundary vertices and the vertices $u_i$, $2 \leq i \leq n-1$, are called the inner vertices of the path. A path $P_1$ is called an extension of another path $P_2$ if $P_2$ is a proper subgraph of $P_1$. The notation $C = C_n(u_1, \ldots, u_n)$ denotes a cycle of length $n$. A cycle $C$ is contractible if it is the boundary of a 2-disk otherwise called non-contractible. If a cycle $C_n(v_1, v_2, \ldots, v_n)$ is the boundary of an $n$-gonal face $F$, then we denote $F$ as $[v_1 v_2 \ldots v_n]$. The notations $G_1 \cup G_2$ and $G_1 \cap G_2$ denote the usual union and intersection of two graphs $G_1$ and $G_2$.

Let $v$ be a vertex in a map such that the link of $v$ is a $k$-length cycle. Then we use bold appearance for those vertices in the expression of link of $v$ which are connected with $v$ by an edge. For example if $f_{seq}(v) = (3^3.4^2)$ then the link of $v$ is a 7-length cycle and therefore we denote ${\rm lk}(v) = C_7(\boldsymbol{v_1}, v_2, v_3, v_4, v_5, \boldsymbol{v_6}, v_7)$. 

\section{Examples} \label{s3}

The following maps  are doubly semi-equivelar maps on the plane $\mathbb{R}^2$. The shaded region in Figure 3 shows a DSEM of type $T_3$ on the torus. 


\vspace{-.5cm}




\vspace{8.5cm}

\section{Proof of Theorem \ref{t1}}
To prove Theorem \ref{t1}, we first prove Lemma \ref{l1} and Lemma \ref{l2}. 

\begin{lem} \label{l1}
	Let $v$ be a vertex with the face-sequence $f$ such that $\phi (v)=0$. Then $f \in X = \{(3^3.4^2)$, $(3^6)$, $(3.4^2.6)$, $(3^2.6^2)$, $(3^4.6)$, $(3^2.4.3.4)$, $(3.6.3.6)$, $(4^4)$, $(3.4.6.4)$, $(3^2.4.12)$, $(4.8^2)$, $(3.12^2)$, $(6^3)$, $(5^2.10)$, $(3.8.24)$, $(3.9.18)$, $(3.10.15)$, $(4.5.20)$, $(3.7.42)$, $(4.6.12)$, $(3.4.3.12) \}$.
\end{lem}

{\bf Proof:} Let $v$ be a vertex with face-sequence $f= (p_1 ^{n_1}. p_2^{n_2}. \ldots. p_k^{n_k})$ and
$\phi (v) = 0$. Then 
$$ \phi(v) = 1 - (\sum_{i=1}^{k} n_i)/2 + \sum_{i=1}^{k} (n_i/p_i) = 0 \Longrightarrow (\frac{1}{2} - \frac{1}{p_1})n_1 + (\frac{1}{2} - \frac{1}{p_2})n_2  + \ldots + (\frac{1}{2} - \frac{1}{p_k})n_k = 1.$$

By Equation (1) [\cite{dm(2017)}, Lemma 2.1], we see that $f \in X$. 	\hfill $\Box$

\vspace{.3cm} 

\begin{lem}\label{l2}
	Let $N_i$, $1 \leq i \leq 22$, be the doubly semi-equivelar maps on the plane $\mathbb{R}^2$ as shown in Section \ref{s3}. Let $f_1 = (3^6)$, $f_2 = (3^3.4^2)$, $f_3 = (4^4)$, $f_4 = (3.4.6.4)$, $f_5 = (3^2.4.3.4)$, $f_6 = (3^2.6^2)$, $f_7 = (3^2.4.12)$, $f_8 = (3^4.6)$, $f_9 = (3.4^2.6)$, $f_{10} = (4.6.12)$, $f_{11}=(3.4.3.12)$, $f_{12}=(3.12^2)$, and $f_{13} =(3.6.3.6)$.
	
	$1.$ If $M_1$, $M_2$, $M_3$, $M_4$, and $M_5$ are DSEMs of types $T_1= [f_1^{({f_1}^2.{f_2}. {f_1}^2.{f_2})}:f_2^{({f_1}^5.{f_2}^2)}]$, $ T_2=[f_1^{({f_1}^2.{f_2}^4)}:f_2^{({f_1}^5.{f_2}^2)}]$, $T_3=[f_2^{({f_2}^4.{f_3}^3)}:f_3^{({f_2}^3.{f_3}.{f_2}^3.{f_3} )}]$,  $T_4=[f_2^{({f_2}^4.{f_3}^3)}:f_3^{({f_2}^3.{f_3}^5)}]$, and $T_5=[f_2^{({f_4}^2. {f_2}.{f_4}^2.{f_2}^2)}:f_4^{({f_4}^5.{f_2}^4)}]$ respectively, then $M_1 \cong N_1$, $M_2 \cong N_2$, $M_3 \cong N_3$, $M_4 \cong N_4$,  and $M_5 \cong N_5$.
	
	$2.$ If $M_6$ and $M_7$ are DSEMs of types $T_6= [f_2^{( {f_2}.{f_5}^6)}:f_5^{( {f_2}^2.{f_5}^3.{f_2}.f_5)}]$ and $T_7=[f_2^{(f_2.{f_5}^3.f_2.{f_5}^2)}:f_5^{( {f_2}^3.f_5.{f_2}^2.f_5)}]$ respectively, then $M_6 \cong N_6$ and $M_7 \cong N_7$.
	
	$3.$ If $M_8$, $M_9$, and $M_{10}$ are DSEMs of types   $T_8=[f_1^{({f_6}^6)}:f_6^{(f_1.{f_6}^9)}]$, $T_9=[f_1^{({f_7}^6)}:f_7^{({f_1}.{f_7}^{13})}]$, 
	and $T_{10}=[f_1^{({f_5}^6)}:f_5^{(f_1.{f_5}^{6})}]$ respectively,  then $M_8 \cong N_8$, $M_9 \cong N_9$, and $M_{10} \cong N_{10}$.
	
	$4.$ If $M_{11},M_{12}$, and $M_{13}$ are the DSEMs of type $T_{11}= [{f_1}^{({{f_1}}^2.{f_8}^2.{f_1}.{f_8})}:{f_8}^{({{f_8}}^5.{f_1}^3)}]$, $T_{12}= [f_1^{({f_8}^6)}:{f_8}^{({{f_8}}^5.{f_1}.{{f_8}}.{f_1})}]$,  and $T_{13}= [f_1^{({{f_8}}^6)}:{f_8}^{({{f_8}}^7.{f_1})}]$ respectively, then $M_{11} \cong N_{11}, M_{12} \cong N_{12}$, and $M_{13} \cong N_{13}$.
	
	$5.$ If $M_{14},M_{15}$, and $M_{16}$ are DSEMs of types $T_{14}=[f_4^{({f_4}^5.{f_9}^{4})}:f_9^{({f_4}^2.{f_9}^{7})}]$, $T_{15}=[f_4^{({f_4}^5.{f_5}^4)}:f_5^{({f_4}^2.{f_5}.{f_4}^2.{f_5}^2)}]$, and $T_{16}=[f_4^{({f_{10}}^{3}.{f_4}^{2}.{f_{10}}^{3}.{f_4})}:f_{10}^{({f_{10}}^{11}.{f_4}.{f_{10}}^{2}.{f_4}^2)}]$ respectively, then $M_{14} \cong N_{14}, M_{15} \cong N_{15}$, and $M_{16} \cong N_{16}$.

	$6.$ If $M_{17}$ and $M_{18}$ are DSEMs of types  $T_{17}=[f_{11}^{F_1}:f_{12}^{F_2}]$ and  $T_{18}=[f_{8}^{(f_{6}^2.f_{8}.f_{6}.f_{8}.f_{6}^2.f_{8})}:f_{6}^{(f_{6}^2.f_{8}^3.f_{6}^2.f_{8}.f_{6}.f_{8})}]$ respectively, where 	$F_1 = (f_{12}^2.f_{11}.f_{12}^2.f_{11}.f_{12}^2.f_{11}.f_{12}^2.f_{11}^3)$ and $F_2 = (f_{11}.f_{12}^2.f_{11}.f_{12}^2$. $f_{11}.f_{12}^2.f_{11}.f_{12}.f_{11}. f_{12}^2.f_{11}.f_{12}^2.f_{11}.f_{12}^2.f_{11})$,  then $M_{17} \cong N_{17}$ and $M_{18} \cong N_{18}$.
	
	$7.$ If $M_{19},M_{20}$, and $M_{21}$ are DSEMs of types $T_{19} = [f_8^{(f^4_8.f_{13}.f_8^2.f_{13})}:f_{13}^{(f^{4}_8.f_{13}.f_8^4.f_{13})}]$, $T_{20} = [(3^4.6)^{(f^7_8.f_{13})}:(3.6.3.6)^{(f^{10}_8)}]$, and $T_{21}$ = $[f_{13}^{({f_6}^4.{f_{13}}.{f_6}^4.{f_{13}})}:f_6^{({f_6}^2.{f_{13}}.{f_6}.{f_{13}}.{f_6}^2.{f_{13}}.{f_6}.{f_{13}})}]$ respectively, then $M_{19} \cong N_{19}$, $M_{20} \cong N_{20}$, and $M_{21} \cong N_{21}$.
	 
\end{lem}

{\bf Proof:} Let $M_1$ be a DSEM of type $T_1=[f_1^{({f_1}^2.{f_2}. {f_1}^2.{f_2})}:f_2^{({f_1}^5.{f_2}^2)}]$, where $f_1 = (3^6)$ and $f_2 = (3^3.4^2)$. Let $w^{0}_{0} \in V(M_1)$ with $f_{seq}(w^{0}_{0}) = (3^3.4^2)$. Suppose the quadrangular and triangular faces incident at $w^{0}_{0}$ are $[w^{0}_{0} w^{-1}_{0} w^{-1}_{1} w^{0}_{1}]$,  $[w^{0}_{0} w^{-1}_{0} w^{-1}_{-1} w^{0}_{-1}]$,  $[w^{0}_{0} w^{0}_{1} w^{1}_{0}]$, $[w^{0}_{0} w^{1}_{0} w^{1}_{-1}]$,  and $[w^{0}_{0} w^{1}_{-1} w^{0}_{-1}]$. Since the face-sequence of ${\rm lk}(w^{0}_{0})$ is $({f_1}^5.{f_2}^2)$, we have $ f_{seq}(w^{0}_{-1}) = f_{seq}(w^{-1}_{-1}) = f_{seq}(w^{-1}_{0}) = f_{seq}(w^{-1}_{1}) = f_{seq}(w^{0}_{1}) = (3^3.4^2)$ and $f_{seq}(w^{1}_{0}) = f_{seq}(w^{1}_{-1}) =(3^6)$. Then the quadrangular and triangular faces incident at $w^{0}_{-1}$ are $[w^{0}_{-1} w^{0}_{0}w^{-1}_{0} w^{-1}_{-1}]$,  $[w^{0}_{-1} w^{-1}_{-1} w^{-1}_{-2} w^0_{-2}]$, $[w^{0}_{-1} w^{1}_{-1} w^{0}_{0}]$,  $[w^{0}_{-1} w^{1}_{-1} w^{1}_{-2}]$, and $[w^{0}_{-1} w^0_{-2} w^{1}_{-2}]$, whereas the quadrangular and triangular faces incident at $w^{0}_{1}$ are $[w^{0}_{1} w^{0}_{0}w^{-1}_{0} w^{-1}_{1} ]$, $[w^{0}_{1} w^{-1}_{1} w^{-1}_{2} w^0_{2}]$, $[w^{0}_{1} w^{1}_{0} w^{0}_{0}]$,  $[w^{0}_{1} w^{1}_{0} w^{1}_{1}]$, and $[w^{0}_{1} w^0_{2} w^{1}_{1}]$. Continuing in this way, we get a path $P_0 =P(  \cdots -w^0_{-2}-w^{0}_{-1}-w^{0}_{0}-w^{0}_{1}-w^0_{2}- \cdots)$ such that all triangles incident at each vertex of $P_0$ lie on one side and all quadrangles incident with the same vertex lie on the other side of $P_0$. If the path $P_0$ contains a cycle $C$, then it gives a map on 2-disk $\mathbb{D}^2$, say $M(\mathbb{D}^2)$, whose boundary vertices satisfy any of the following two conditions:

$(i)$ every boundary vertices of $M(\mathbb{D}^2)$ contains two quadrangles only.

$(ii)$ every boundary vertices of $M(\mathbb{D}^2)$ contains three triangles only.

Assume that there are $n=n_{1}+n_{2}$ inner vertices and $m$ boundary vertices in $M(\mathbb{D}^2)$, where $n_1$ and $n_2$ denote the number of inner vertices with face-sequences  $(3^6)$ and $(3^3.4^2)$, respectively.

Let $v, e$ and $f$ denote the number of vertices, edges and faces in $M(\mathbb{D}^2)$. If the boundary vertices of $M(\mathbb{D}^2)$ satisfy $(i)$, then $v=m+n_{1}+n_{2}$, $e= 3m/2 + 6n_{1}/2 + 5n_{2}/2$ and  $f=2m/4 + 6n_{1}/3 + 3n_{2}/3 + 2n_{2}/4$, whereas if $(ii)$ holds, then $v= m + n_{1}+ n_{2}$, $e= 4m/2 + 6n_{1}/2 + 5n_{2}/2$, and  $f= 3m/3 + 6n_{1}/3 + 2n_{2}/4 +3n_{2}/3$. In both the conditions, we get $v-e-f=0$. This is not possible as the Euler characteristic of  $\mathbb{D}^2$ is 1. Therefore, $P_0$ is an infinite path.

Using the same procedure for the vertex $w^{1}_{0}$ with face-sequence $(3^6)$ and the fact that the face-sequence ${\rm lk}(w^{1}_{0})$ is $({f_1}^2.{f_2}. {f_1}^2.{f_2})$, we get a path $P_1 = P (\cdots -w^{1}_{-2}-w^{1}_{-1}-w^{1}_{0}-w^{1}_{1}-w^{1}_{2}- \cdots)$ such that every inner vertex of the path has face-sequence $(3^6)$. Similarly as above, $P_1$ does not contains any cycle. Thus, $P_1$ is an infinite path.

The faces that are incident at the vertices of the path $P_0$ (resp. $P_1$) create an endless strip that is bounded by two infinite paths $P_{-1}$ and $P_1$ (resp. $P_{0}$ and $P_2$). Note that the faces  between $P_0$ and $P_{-1}$ (resp. $P_0$ and $P_1$) are quadrangles (resp. triangles). 

Continuing in this manner, we obtain paths $ \cdots, P_{-2}, P_{-1}, P_{0}, P_{1}, P_{2}, \cdots$ such that $(a)$ the faces between $P_{3k}$ and $P_{3k-1}$ are quadrangles $(b)$ the faces between $P_{3k}$ and $P_{3k+1}$ (resp. $P_{3k+1}$ and $P_{3k+2}$) are triangles $(c)$ the faces at $w^{3k}_{i}$ are $[w^{3k}_{i} w^{3k-1}_{i} w^{3k-1}_{i+1} w^{3k}_{i+1}]$, $[w^{3k}_{i}w^{3k-1}_{i}w^{3k-1}_{i-1}w^{3k}_{i-1}]$, $[w^{3k}_{i}w^{3k}_{i+1}w^{3k+1}_{i}]$, $[w^{3k}_{i}w^{3k+1}_{i}w^{3k+1}_{i-1}]$ and 
$[w^{3k}_{i}w^{3k+1}_{i-1}w^{3k}_{i-1}]$ $(d)$ the faces at $w^{3k+1}_{i}$ are 
$[w^{3k+1}_{i}w^{3k+2}_{i-1}$ $w^{3k+2}_{i}]$, $[w^{3k+1}_{i}w^{3k+2}_{i}w^{3k+1}_{i+1}]$, $[w^{3k+1}_{i}w^{3k+2}_{i-1}w^{3k+1}_{i-1}]$,
$[w^{3k+1}_{i}w^{3k+1}_{i-1}w^{3k}_{i}]$, $[w^{3k+1}_{i}w^{3k}_{i}w^{3k}_{i+1}]$, $[w^{3k+1}_{i}w^{3k}_{i+1} \linebreak w^{3k+1}_{i+1}]$ $(e)$ the faces at $w^{3k+2}_{i}$ are $[w^{3k+2}_{i}w^{3k+2}_{i+1}w^{3k+3}_{i+1}w^{3k+3}_{i}]$, $[w^{3k+2}_{i}w^{3k+3}_{i}w^{3k+3}_{i-1}w^{3k+2}_{i-1}]$, $[w^{3k+2}_{i}w^{3k+2}_{i-1} \linebreak w^{3k+1}_{i}]$, $[w^{3k+2}_{i}w^{3k+1}_{i}w^{3k+1}_{i+1}]$ and $[w^{3k+2}_{i}w^{3k+1}_{i+1}w^{3k+2}_{i+1}]$. Now the mapping $f : V(M_1) \to  V(N_1)$ defined by $f(w^r_{s}) = v^r_{s}$ for $r,s \in \mathbb{Z}$ is an isomorphism. Thus, $M_1 \cong N_1$.

Let $M_2$ be a  DSEM of the type $T_2= [f_1^{({f_1}^2.{f_2}^4)}:f_2^{({f_1}^5.{f_2}^2)}]$. By the same argument used for $M_1$, we get $M_2 \cong N_2$. In order to show that the paths in $M_2$ corresponding to the thick black paths in $N_2$ are infinite, we use the fact that there is no map on the 2-disk  $\mathbb{D}^2$ that fulfills $(i)$ and $(ii)$.

\smallskip
If $M_3$, $M_4$, and $M_5$ are DSEMs of types $T_3=[f_2^{({f_2}^4.{f_3}^3)}:f_3^{({f_2}^3.{f_3}.{f_2}^3.{f_3} )}]$, $T_4=[f_2^{({f_2}^4.{f_3}^3)}:f_3^{({f_2}^3.{f_3}^5)}]$, and  $T_5=[f_2^{({f_4}^2. {f_2}.{f_4}^2.{f_2}^2)}:f_4^{({f_4}^5.{f_2}^4)}]$ respectively, then by using the same reasoning as used for $M_1$, we get $M_3 \cong N_3$,  $M_4 \cong N_4$, and $M_5 \cong N_5$.

\smallskip
Let $M_6$ be a DSEM of type $T_6= [f_2^{( {f_2}.{f_5}^6)}:f_5^{( {f_2}^2.{f_5}^3.{f_2}.f_5)}]$, where $f_2 = (3^3.4^2)$ and $f_5 = (3^2.4.3.4)$. Let $w_{0,0} \in V(M_6)$ and $f_{seq}(w_{0,0}) = (3^2.4.3.4)$. Suppose the quadrangular and triangular faces at $w^0_{0}$ are $[w^0_{0}w^{0}_{1}w^{1}_{1}w^{1}_{0}]$, $[w^0_{0}w^{-1}_{0}w^{-1}_{-1}w^{0}_{-1}]$, $[w^0_{0}w^{0}_{1}w^{0}_{-1}]$, $[w^0_{0}w^{-1}_{0}x^0_{0}]$, and $[w^0_{0}x^0_{0}w^{0}_{1}]$. Since the face-sequence of ${\rm lk}(w^0_{0})$ is ${( {f_2}^2.f^3_5.{f_2}.{f_5})}$, we have $ f_{seq}(w^{0}_{1}) = f_{seq}(w^{1}_{1}) = f_{seq}(w^{1}_{0}) = f_{seq}(w^{-1}_{0}) = (3^2.4.3.4)$ and $f_{seq}(w^{0}_{-1}) = f_{seq}(w^{-1}_{-1}) = f_{seq}(x^{0}_{0}) =(3^3.4^2)$.

Proceeding similarly, as in $M_1$, we get a path $P_0 = P(\cdots -w^{0}_{-2}-w^{0}_{-1}-w^{0}_{0}-w^{0}_{1}-w^{0}_{2}- \cdots)$. If the path $P_0$ contains a cycle $C$, then it gives a map on 2-disk $\mathbb{D}^2$, say $M(\mathbb{D}^2)$, whose boundary vertices satisfy any of the following two conditions: 

$(i)$ the boundary vertices of $M(\mathbb{D}^2)$ are of the type $(3^3,4^0), (3^1,4^1)$ and $(3^2,4^1)$. 

$(ii)$ the boundary vertices of $M(\mathbb{D}^2)$ are of the type  $(3^0,4^2), (3^2,4^1)$ and $(3^1,4^1)$.

Assume that there are $n=n_{1}+n_{2}$ inner vertices and $m$ boundary vertices in $M(\mathbb{D}^2)$, where $n_1$ and $n_2$ denote the number of inner vertices with face-sequences  $(3^3.4^2)$ and $(3^2.4.3.4)$, respectively. 
If the boundary vertices of $M(\mathbb{D}^2)$ satisfy $(i)$ or $(ii)$, then in both the cases, we get $v-e-f \neq 1$. As a result, $P_0$ is an infinite path.

Using the same procedure for the vertex $w^{-1}_{-1}$ with face-sequence $(3^3.4^2)$ and the fact that the face-sequence of ${\rm lk}(w^{-1}_{-1})$ is $({f_1}.{f_2}^6)$, we get a path $P_{-1} = P (\cdots -w^{-1}_{-2}-w^{-1}_{-1}-w^{-1}_{0}-w^{-1}_{1}-w^{-1}_{2}- \cdots)$. Similarly as above, $P_{-1}$ is an infinite path. Continuing in this way, we obtain paths $ \cdots, P_{-2}, P_{-1}, P_{0}, P_{1}, P_{2}, \cdots$ and the vertices $ \ldots, x^{2i}_{-2}, x^{2i}_{-1}, x^{2i}_0, x^{2i}_1, x^{2i}_2, \ldots$ having face-sequence $(3^3.4^2)$ lying between $P_{2i}$ and $P_{2i-1}$ such that $(a)$ the paths $P_{2i}$'s appear as shown in Figure 4.1, $(b)$ the paths $P_{2i-1}$'s appear as shown in Figure 4.2. Now, the function $f : V(M_6) \to  V(N_6)$ defined by $f(w^r_{s}) = v^r_{s}$ and $f(y^{2i}_{j}) = x^{2i}_{j}$ for $r,s,i,j \in \mathbb{Z}$ is an isomorphism. Thus, $M_6 \cong N_6$.

For the remaining types maps, similarly, we show that the paths in $M_i$ (for $7 \leq i \leq 22$) corresponding to the thick paths shown in the respective $N_i$ are infinite. This further gives suitable isomorphism between $M_i$ and $N_i$. Thus the proof.   

\vspace{0cm}

\begin{picture}(0,0)(-34,20)
\setlength{\unitlength}{7mm}

\drawpolygon(.5,0)(.5,1)(-.5,1)(-.5,0)
\drawpolygon(-.5,1)(.5,1)(0,2)


\drawpolygon(.5,0)(1.5,-.5)(1.5,.5)
\drawpolygon(.5,0)(1.5,.5)(.5,1)
\drawpolygon(1.5,-.5)(2.5,-.5)(2.5,.5)(1.5,.5)
\drawpolygon(2.5,-.5)(3.5,-.5)(3.5,.5)(2.5,.5)


\drawpolygon(-.5,0)(-1.5,-.5)(-1.5,.5)
\drawpolygon(-.5,0)(-1.5,.5)(-.5,1)
\drawpolygon(-1.5,-.5)(-2.5,-.5)(-2.5,.5)(-1.5,.5)
\drawpolygon(-2.5,-.5)(-3.5,-.5)(-3.5,.5)(-2.5,.5)

\drawpolygon(0,2)(.5,1)(1,2)
\drawpolygon(.5,1)(1.5,.5)(2,1.5)(1,2)
\drawpolygon(1.5,.5)(2.5,.5)(2,1.5)
\drawpolygon(2.5,.5)(3,1.5)(2,1.5)
\drawpolygon(2.5,.5)(3.5,.5)(3,1.5)

\drawpolygon(0,2)(-.5,1)(-1,2)
\drawpolygon(-.5,1)(-1.5,.5)(-2,1.5)(-1,2)
\drawpolygon(-1.5,.5)(-2.5,.5)(-2,1.5)
\drawpolygon(-2.5,.5)(-3,1.5)(-2,1.5)
\drawpolygon(-2.5,.5)(-3.5,.5)(-3,1.5)

\drawline[AHnb=0](3.5,.5)(3.7,.6)
\drawline[AHnb=0](3.8,.65)(4,.75)
\drawline[AHnb=0](4.1,.8)(4.3,.9)
\drawline[AHnb=0](3.5,.5)(3.7,.6)
\drawline[AHnb=0](3.8,.65)(4,.75)
\drawline[AHnb=0](4.1,.8)(4.3,.9)

\drawline[AHnb=0](-3.5,.5)(-3.7,.6)
\drawline[AHnb=0](-3.8,.65)(-4,.75)
\drawline[AHnb=0](-4.1,.8)(-4.3,.9)
\drawline[AHnb=0](-3.5,.5)(-3.7,.6)
\drawline[AHnb=0](-3.8,.65)(-4,.75)
\drawline[AHnb=0](-4.1,.8)(-4.3,.9)

\drawpolygon[fillcolor=black](-3.5,.5)(-1.5,.5)(-.5,1)(.5,1)(1.5,.5)(3.5,.5)(3.5,.6)(1.5,.6)(.5,1.1)(-.5,1.1)(-1.5,.6)(-3.5,.6)

\put(-4.3,.35){\scriptsize {\tiny$z_{j-1}$}}
\put(-3,.25){\scriptsize {\tiny$w_{i}$}}
\put(-2,.25){\scriptsize {\tiny$z_{j}$}}
\put(-1.5,1.1){\scriptsize {\tiny$z_{j+1}$}}
\put(.7,1.1){\scriptsize {\tiny$z_{j+2}$}}
\put(1.55,.15){\scriptsize {\tiny$z_{j+3}$}}
\put(2.55,.15){\scriptsize {\tiny$w_{i+1}$}}
\put(3.55,.15){\scriptsize {\tiny$z_{j+4}$}}

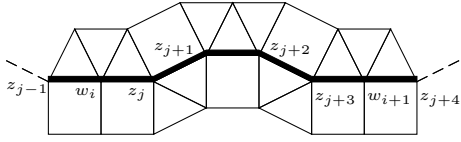
\put(-4.2,-1.2){\scriptsize {\tiny {\bf Figure 4.1:} Path of type $P_{2i}$ (indicated by thick line)}}
\end{picture}

\begin{picture}(0,0)(-114,4)
\setlength{\unitlength}{7mm}

\drawpolygon(.5,0)(.5,-1)(0,-2)(-.5,-1)(-.5,0)
\drawline[AHnb=0](3,-1.5)(2,-1.5)(1,-2)(0,-2)(-1,-2)(-2,-1.5)(-3,-1.5)(-2.5,-.5)(-2,-1.5)(-1.5,-.5)(-.5,-1)(-1,-2)

\drawline[AHnb=0](1,-2)(.5,-1)(1.5,-.5)(2,-1.5)(2.5,-.5)(3,-1.5)(3.5,-.5)

\drawline[AHnb=0](-3,-1.5)(-3.5,-.5)


\drawpolygon(.5,0)(1.5,-.5)(1.5,.5)

\drawpolygon(1.5,-.5)(2.5,-.5)(2.5,.5)(1.5,.5)
\drawpolygon(2.5,-.5)(3.5,-.5)(3.5,.5)(2.5,.5)


\drawpolygon(-.5,0)(-1.5,-.5)(-1.5,.5)

\drawpolygon(-1.5,-.5)(-2.5,-.5)(-2.5,.5)(-1.5,.5)
\drawpolygon(-2.5,-.5)(-3.5,-.5)(-3.5,.5)(-2.5,.5)


\drawline[AHnb=0](3.5,-.5)(3.7,-.6)
\drawline[AHnb=0](3.8,-.65)(4,-.75)
\drawline[AHnb=0](4.1,-.75)(4.3,-.85)


\drawline[AHnb=0](-3.5,-.5)(-3.7,-.6)
\drawline[AHnb=0](-3.8,-.65)(-4,-.75)
\drawline[AHnb=0](-4.1,-.75)(-4.3,-.85)


\drawpolygon[fillcolor=black](-3.5,-.5)(-1.5,-.5)(-.5,-1)(.5,-1)(1.5,-.5)(3.5,-.5)(3.5,-.4)(1.5,-.4)(.5,-.9)(-.5,-.9)(-1.5,-.4)(-3.5,-.4)

\put(-4.3,-.2){\scriptsize {\tiny$z_{j-1}$}}
\put(-3,-.2){\scriptsize {\tiny$w_{i}$}}
\put(-2,-.2){\scriptsize {\tiny$z_{j}$}}
\put(-1.5,-1.15){\scriptsize {\tiny$z_{j+1}$}}
\put(.7,-1.1){\scriptsize {\tiny$z_{j+2}$}}
\put(1.55,-.15){\scriptsize {\tiny$z_{j+3}$}}
\put(2.55,-.2){\scriptsize {\tiny$w_{i+1}$}}
\put(3.55,-.15){\scriptsize {\tiny$z_{j+4}$}}

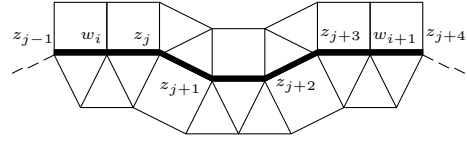
\put(-4.2,-2.75){\scriptsize {\tiny {\bf Figure 4.2:} Path of type $P_{2i-1}$ (indicated by thick line)}}
\end{picture}

\vspace{2.1cm}

\hfill $\Box$


\vspace{.3cm}

{\bf Proof of Theorem \ref{t1}:} Suppose $M$ is a DSEM on the plane consists of face-sequence pair $\{f_1, f_2\} \in X$, as in Lemma \ref{l1}. Then for a fix $f_1 \in X$, we find $f_2 \in X$ that is compatible with $f_1$. For example, let $u$ be a vertex with face-sequence $f_1 = (3^3.4^2)$ and  ${\rm lk}(u) = C_7(\boldsymbol{u_1},u_2,u_3,u_4,u_5,\boldsymbol{u_6},u_7)$. Then, for $f_2$ we have one of the four cases: $(i)$ there must be at least one 4-gon, $(ii)$ two consecutive 3-gons $(iii)$ successive 3-gon and 4-gon, $(iv)$ two consecutive 4-gons. 
This implies that $f_2 \in \{(3^6), (3.4^2.6), (3^2.6^2), (3^4.6), (3^2.4.3.4), (4^4), (3.4.6.4),(3^2.4.12)$, $(4.8^2),(4.5.20),(4.6.12),(3.4.3.12)\}$. Similarly for each $f_1 \in X$, we obtain $f_2$ compatible with $f_1$. This gives $(f_1, f_2) \in A \cup B$, where 

\smallskip

$A =  \{((3^3.4^2)$, $(3^2.6^2))$, $((3^3.4^2)$, $(3.4^2.6))$, $((3^3.4^2)$, $(3^4.6))$, $((3^3.4^2)$, $(3^2.4.12))$, $((3^3.4^2)$, $(4.8^2))$, $((3^3.4^2)$, $(4.5.20))$, $((3^3.4^2)$, $(4.6.12))$, $((3^3.4^2)$, $(3.4.3.12))$, $((3.4^2.6)$, $(3^2.6^2))$, $((3.4^2.6)$, $(3^4.6))$, $((3.4^2.6)$, $(3^2.4.3.4))$, $(3.4^2.6)$, $(4^4)$, $((3.4^2.6)$, $(3^2.4.12))$, $((3.4^2.6)$, $(4.8^2))$, $((3.4^2.6)$, $(6^3))$, $((3.4^2.6)$, $(4.5.20))$, $((3.4^2.6)$, $(4.6.12))$, $((3.4^2.6)$, $(3.4.3.12))$, $((3^2.6^2)$, $(3^2.4.3.4))$, $((3^2.6^2)$, $(3^2.4.12))$, $((3^2.6^2)$, $(6^3))$, $((3^4.6)$, $(3^2.4.3.4))$, $((3^4.6)$, $(3^2.4.12))$, $((3^4.6)$, $(6^3))$, $((3^4.6)$, $(4.6.12))$, $((3^2.4.3.4)$, $(4^4))$, $((3^2.4.3.4)$, $(3^2.4.12))$, $((3^2.4.3.4)$, $(3.4.3.12))$,  $(3^2.4.3.4)$, $(4.8^2))$, $((3.6.3.6)$, $(4.6.12))$, $((3.6.3.6)$, $(6^3))$, $((4^4)$, $(3.4.6.4))$, $((4^4)$, $(3^2.4.12))$, $((4^4)$, $(4.8^2))$, $((4^4)$, $(4.5.20))$, $((4^4)$, $(4.6.12))$, $((4^4)$, $(3.4.3.12))$, $((3.4.6.4)$, $(3^2.4.12))$, $((3.4.6.4)$, $(6^3))$,$((3.4.6.4)$, $(4.8^2))$, $((3.4.6.4)$, $(4.5.20))$, $((3.4.6.4)$, $(3.4.3.12))$, $((3^2.4.12)$, $(3.12^2))$, $((3^2.4.12)$, $(3.4.3.12))$, $((3^2.4.12)$, $(4.5.20))$,   $((3^2.4.12)$, $(4.6.12))$, $((3^2.4.12)$, $(4.8^2))$,  $((4.8^2)$, $(4.5.20))$, $((4.8^2)$, $(4.6.12))$, $((4.8^2)$, $(3.4.3.12))$, $((3.12^2)$, $(3.4.3.12))$, $((6^3)$, $(4.6.12))$, $((5^2.10)$, $(4.5.20)) \}$, 

\smallskip

$B = \{ ((3^3.4^2)$, $(3^6))$, $((3^3.4^2)$, $(4^4))$, $((3^3.4^2)$, $(3^2.4.3.4))$, $((3^3.4^2)$, $(3.4.6.4))$, $((3^6)$, $(3^2.6^2))$, $((3^6)$, $(3^4.6))$, $((3^6)$, $(3^2.4.12))$, $((3^6)$, $(3^2.4.3.4))$, $((3.12^2)$, $(3.4.3.12))$,  $((3^2.6^2)$, $(3.6.3.6))$, $((3.4^2.6)$, $(3.6.3.6))$, $((3.4^2.6)$, $(3.4.6.4))$, $((3^2.4.3.4)$, $(3.4.6.4))$, $((3.4.6.4)$, $(4.6.12))$, $((3^2.6^2),(3^4.6))$, $((3^4.6)$, $(3.6.3.6))\}$.

\begin{claim}\label{c1}
	There exists no DSEM for  $(f_1, f_2) \in A$.
\end{claim}

\noindent{\bf Case 1:} First consider the pairs $((3^3.4^2), f_2) \in A$. Without loss of generality let $u$ be a vertex in $M$ with $f_{seq}(u)=f_1 = (3^3.4^2)$ such that ${\rm lk}(u) = C_7(\boldsymbol{u_1},u_2,u_3,u_4,u_5,\boldsymbol{u_6},u_7)$. Then:

\smallskip

\noindent{\bf Subcase 1.1:} Suppose $f_2 = (3^2.6^2)$. Then at least one vertex in $\{u_3, u_4\}$ has the face-sequence $(3^2.6^2)$. If $f_{seq}(u_3) = (3^2.6^2)$ (or $f_{seq}(u_4) = (3^2.6^2)$), then $u_2$ (resp. $u_5$) has the face-sequence $f_3 \neq f_1, f_2$. Thus $f_2 \neq (3^2.6^2)$, that is, for the pair $((3^3.4^2), (3^2.6^2))$, $M$ does not exist.

\smallskip

\noindent{\bf Subcase 1.2:} Suppose $f_2 = (3.4^2.6)$. Then at least one vertex in the set $\{u_1,u_2,u_5,u_6, u_7\}$ has the face-sequence $(3.4^2.6)$.  If $f_{seq}(u_2) = (3.4^2.6)$ (or $f_{seq}(u_5) = (3.4^2.6)$), then $f_{seq}(u_3) \neq f_1, f_2$ (resp. $f_{seq}(u_4) \neq f_1,f_2)$. Thus $f_{seq}(u_2) \neq (3.4^2.6)$ and $f_{seq}(u_5) \neq (3.4^2.6)$.

If $f_{seq}(u_7) = (3.4^2.6)$, then ${\rm lk}(u_7) = C_9(u_1,\boldsymbol{u_2},u,\boldsymbol{u_5},u_6,u_8,\boldsymbol{u_9}, \boldsymbol{u_{10}}, \boldsymbol{u_{11}})$ or ${\rm lk}(u_7) = C_9(u_1,\boldsymbol{u_2}$, $u,\boldsymbol{u_5},u_6,\boldsymbol{u_8},\boldsymbol{u_{9}}, \boldsymbol{u_{10}},u_{11})$. In case, ${\rm lk}(u_7) = C_9(u_1,\boldsymbol{u_2},u,\boldsymbol{u_5},u_6,u_8,\boldsymbol{u_9}, \boldsymbol{u_{10}}, \boldsymbol{u_{11}})$, $f_{seq}(u_6)$ is either $(3^3.4^2)$ or $(3.4^2.6)$. In both the 
possibilities, we get $f_{seq}(u_8) \neq f_1, f_2$.
On the other hand, when ${\rm lk}(u_7) = C_9(u_1,\boldsymbol{u_2},u,\boldsymbol{u_5},u_6$, $\boldsymbol{u_8},\boldsymbol{u_{9}}, \boldsymbol{u_{10}},u_{11})$, then $f_{seq}(u_{1}) = (3.4^2.6)$. This gives that $f_{seq}(u_{11}) \neq f_1,f_2$. So, $f_{seq}(u_7) \neq (3.4^2.6)$. Proceeding similarly, we see that $f_{seq}(u_{1})$,$f_{seq}(u_6) \neq (3.4^2.6)$. Thus for the pair $((3^3.4^2), (3.4^2.6))$, $M$ does not exist.


\smallskip

\noindent{\bf Subcase 1.3:} Suppose $f_2 = (3^4.6)$. Then at least one vertex in $\{u_3, u_4\}$ has the face-sequence $(3^4.6)$. Without loss of generality let $f_{seq}(u_3) = (3^4.6)$. Then ${\rm lk}(u_3) = C_8(u, u_2, u_8, \boldsymbol{u_9}$, $\boldsymbol{u_{10}}, \boldsymbol{u_{11}}, u_{12}$, $u_4)$ or 
${\rm lk}(u_3) = C_8(u, u_2, u_8, u_9, \boldsymbol{u_{10}}, \boldsymbol{u_{11}}, \boldsymbol{u_{12}}, u_4)$. In the first case of ${\rm lk}(u_3)$, completing successively ${\rm lk}(u_8)$, ${\rm lk}(u_9)$, ${\rm lk}(u_{10})$, ${\rm lk}(u_{11})$, ${\rm lk}(u_{12})$ and ${\rm lk}(u_4)$, we get $f_{seq}(u_5) \neq f_1, f_2$. On the other hand, when  ${\rm lk}(u_3) = C_8(u, u_2, u_8, u_9, \boldsymbol{u_{10}}, \boldsymbol{u_{11}}, \boldsymbol{u_{12}}, u_4)$, then ${\rm lk}(u_2) = C_7(\boldsymbol{u_{13}}, u_{14}, u_8, u_3, u, \boldsymbol{u_7})$. This gives that $f_{seq}(u_8)$ is either $(3^3.4^2)$ or $(3^4.6)$. But for both the cases, we see that $f_{seq}(u_9)$ or $f_{seq}(u_{14})$ is neither $f_1$ nor $f_2$. Thus, we do not get $M$ for the pair $((3^3.4^2), (3^4.6))$.

\smallskip

\noindent{\bf Subcase 1.4:} Now suppose $f_2 = (3^2.4.12)$. Then at least one vertex in $\{u_1, u_2, u_3, u_4, u_5,u_6\}$ has the face-sequence $(3^2.4.12)$. If $f_{seq}(u_3) = (3^2.4.12)$ (or $f_{seq}(u_4) = (3^2.4.12)$), then $f_{seq}(u_2) \neq f_1, f_2$ (resp. $f_{seq}(u_5) \neq f_1, f_2$).

If $f_{seq}(u_1) = (3^2.4.12)$ such that 
${\rm lk}(u_1) = C_{14}(u_2, \boldsymbol{u}, u_7, u_8,u_9, \boldsymbol{u_{10}}, \boldsymbol{u_{11}}, \boldsymbol{u_{12}}, \boldsymbol{u_{13}}, \boldsymbol{u_{14}}, \boldsymbol{u_{15}}, \boldsymbol{u_{16}}, \linebreak \boldsymbol{u_{17}}, \boldsymbol{u_{18}})$). Then $f_{seq}(u_2) = (3^2.4.12)$ and $f_{seq}(u_i) = (3^2.4.12)$, for $9 \leq i \leq 18$. This gives that $f_{seq}(u_3) \neq f_1, f_2$. So, $f_{seq}(u_1) \neq (3^2.4.12)$. Then we see that $f_{seq}(u_6) \neq (3^2.4.12)$.

If $f_{seq}(u_2) = (3^2.4.12)$, then we get $f_{seq}(u_3) \neq f_1, f_2$. So, $f_{seq}(u_2) \neq (3^2.4.12)$. Similarly we see that $f_{seq}(u_5) \neq (3^2.4.12)$.
Thus for the pair $((3^3.4^2), (3^2.4.12))$, $M$ does not exist.

\smallskip

\noindent{\bf Subcase 1.5:} Suppose $f_2 = (4.8^2)$. Then at least one vertex in $\{u_1, u_6\}$ has the face-sequence $(4.8^2)$. If $f_{seq}(u_1) = (4.8^2)$ (or $f_{seq}(u_6) = (4.8^2)$), then we see $f_{seq}(u_2) \neq f_1, f_2$ (resp. $f_{seq}(u_5) \neq f_1, f_2$). Thus, there is no map for the pair $((3^3.4^2), (4.8^2))$.

\smallskip

\noindent{\bf Subcase 1.6:} Suppose $f_2 = (4.5.20)$. Then at least one vertex in $\{u_1, u_6\}$ has the face-sequence $(4.5.20)$. This gives that $f_{seq}(u_7) \neq f_1, f_2$. Thus, $M$ does not exist for the pair $((3^3.4^2), (4.5.20))$.

\smallskip

\noindent{\bf Subcase 1.7:} Suppose $f_2 = (4.6.12)$. Then at least one vertex in $\{u_1, u_6\}$ has the face-sequence $(4.6.12)$. This gives that $f_{seq}(u_7) \neq f_1, f_2$. Thus, $M$ does not exist for the pair $((3^3.4^2), (4.6.12))$.

\smallskip

\noindent{\bf Subcase 1.8:} Suppose $f_2 = (3.4.3.12)$. Then at least one vertex in $\{u_1,u_2, u_5, u_6\}$ has the face-sequence $(3.4.3.12)$. If $f_{seq}(u_2) = (3.4.3.12)$ (or $f_{seq}(u_5) = (3.4.3.12)$), then $f_{seq}(u_3) \neq f_1, f_2$ (resp. $f_{seq}(u_4) \neq f_1, f_2$). If $f_{seq}(u_1) = (3.4.3.12)$, then $f_{seq}(u_{2})=(3.4.3.12)$). This gives that $f_{seq}(u_{3} \neq f_1,f_2)$. Similarly, we see that $f_{seq}(u_6) \neq (3.4.3.12)$. Thus, $M$ does not exist. 
  
\smallskip

\noindent{\bf Case 2:} Let $((3.4^2.6), f_2) \in A$ and $u$ be a vertex with $f_{seq}(u) = (3.4^2.6)$ such that ${\rm lk}(u) = C_9(u_1, \boldsymbol{u_2}, u_3, \boldsymbol{u_4}, u_5, u_6, \boldsymbol{u_7}, \boldsymbol{u_8}, \boldsymbol{u_9})$. Then we have the following cases.

\noindent{\bf Subcase 2.1:} If $f_2 = (3^2.6^2)$, then at least one vertex in $\{u_6, u_7, u_8, u_9\}$ has the face-sequence $(3^2.6^2)$. In case, $f_{seq}(u_6) = (3^2.6^2)$ or $f_{seq}(u_7) = (3^2.6^2)$, we see that $f_{seq}(u_5) \neq f_1, f_2$. In case, $f_{seq}(u_8) = (3^2.6^2)$, we see that either $f_{seq}(u_6) \neq f_1, f_2$ or $f_{seq}(u_1) \neq f_1, f_2$. Similarly, $f_{seq}(u_9) \neq f_1, f_2$. Thus, $M$ does not exist for $((3.4^2.6), (3^2.6^2))$. 

\noindent{\bf Subcase 2.2:} If $f_2 = (3^4.6)$, then at least one vertex in $\{u_6, u_7, u_8, u_9\}$ has the face-sequence $(3^4.6)$. If $f_{seq}(u_6) = (3^4.6)$ then $f_{seq}(u_5) \neq f_1, f_2$. If $f_{seq}(u_7) = (3^4.6)$ then $f_{seq}(u_6) \neq f_1, f_2$.  If $f_{seq}(u_8) = (3^4.6)$ then $f_{seq}(u_9) \neq f_1, f_2$. If $f_{seq}(u_9) = (3^4.6)$ then ${\rm lk}(u_9) = C_8(\boldsymbol{u}, \boldsymbol{u_6}, \boldsymbol{u_7}, u_8, u_{10}, u_{11}, u_{12}, u_1)$, which implies $f_{seq}(u_{12}) \neq f_1, f_2$. Thus, $M$ does not exist for $((3.4^2.6), (3^4.6))$.  

\noindent{\bf Subcase 2.3:} If $f_2 = (3^2.4.3.4)$, then at least one vertex in $\{u_2, u_4, u_5\}$ has the face-sequence $(3^2.4.3.4)$. If $f_{seq}(u_2) = (3^2.4.3.4)$ then $f_{seq}(u_1) \neq f_1, f_2$. If $f_{seq}(u_5) = (3^2.4.3.4)$ then $f_{seq}(u_6) \neq f_1, f_2$. If $f_{seq}(u_4) = (3^2.4.3.4)$ then ${\rm lk}(u_4) = C_7(\boldsymbol{u}, u_3, u_{10}, \boldsymbol{u_{11}}, u_{12}, u_{13}, u_5)$ or ${\rm lk}(u_4) = C_7(\boldsymbol{u}, u_3$, $u_{10}, u_{11},\boldsymbol{u_{12}}, u_{13}, u_5)$. In case, ${\rm lk}(u_4) = C_7(\boldsymbol{u}, u_3, u_{10}, \boldsymbol{u_{11}}, u_{12}, u_{13}, u_5)$, we see that $f_{seq}(u_{13}) \neq f_1, f_2$. On the other hand, for ${\rm lk}(u_4) = C_7(\boldsymbol{u}, u_3, u_{10}, u_{11},\boldsymbol{u_{12}}, u_{13}, u_5)$, completing successively, ${\rm lk}(u_5)$, ${\rm lk}(u_6)$, ${\rm lk}(u_7)$, ${\rm lk}(u_8)$, ${\rm lk}(u_9)$, ${\rm lk}(u_1)$, ${\rm lk}(u_2)$ and ${\rm lk}(u_3)$, we see that ${\rm lk}(u_{10}) \neq f_1, f_2$. 
Thus, $M$ does not exist for $((3.4^2.6), (3^2.4.3.4))$.  

\noindent{\bf Subcase 2.4:} If $f_2 = (4^4)$, then at least one vertex in $\{u_2, u_3, u_4\}$ has the face-sequence $(4^4)$. If $f_{seq}(u_2) = (4^4)$, then $f_{seq}(u_1) = (3.4^2.6)$. If ${\rm lk}(u_1) = C_9(u_9,u_{10},\boldsymbol{u_{11}},u_2,\boldsymbol{u_{3}},u,\boldsymbol{u_{6}},\boldsymbol{u_{7}},\boldsymbol{u_{8}})$, then $f_{seq}(u_{10}) = (3.4^2.6)$. This gives that $f_{seq}(u_9) \neq f_1,f_2$. Similarly for $f_{seq}(u_3) = (4^4)$ or $f_{seq}(u_4) = (4^4)$, we get a contradiction.  

\noindent{\bf Subcase 2.5:} If $f_2 = (3^2.4.12)$, then at least one vertex in $\{u_2, u_4, u_5\}$ has the face-sequence $(3^2.4.12)$. If $f_{seq}(u_2) = (3^2.4.12)$ then $f_{seq}(u_3) \neq f_1, f_2$. If $f_{seq}(u_4) = (3^2.4.12)$ then either $f_{seq}(u_5) \neq f_1, f_2$ or $f_{seq}(u_6) \neq f_1, f_2$. If $f_{seq}(u_5) = (3^2.4.12)$ then $f_{seq}(u_6) \neq f_1, f_2$.

\noindent{\bf Subcase 2.6:} If $f_2 \in \{(4.8^2), (4.5.20)\}$, then at least one of the vertex in the set $\{u_2, u_4\}$ has the face-sequence $f_2$. If $f_{seq}(u_2) = f_2$ or $f_{seq}(u_4) = f_2$, we see that $f_{seq}(u_3) \neq f_1, f_2$.

\noindent{\bf Subcase 2.7:} If $f_2 = (6^3)$, then at least one vertex in $\{u_7, u_8,u_9\}$ has the face-sequence $(6^3)$. If $f_{seq}(u_7) = (6^3)$ then $f_{seq}(u_6) \neq f_1, f_2$. If $f_{seq}(u_8) =(6^3)$ or $f_{seq}(u_9) = (6^3)$ then $f_{seq}(u_1) \neq f_1, f_2$. 


\noindent{\bf Subcase 2.8:} If $f_2 = (4.6.12)$, then at least one vertex in $\{u_1,u_2, u_4\}$ has the face-sequence $(4.6.12)$. If $f_{seq}(u_1)$ or $f_{seq}(u_2)$ or $f_{seq}(u_4)$ is $(4.6.12)$, then $f_{seq}(u_3) \neq f_1, f_2$.

\noindent{\bf Subcase 2.9:} If $f_2 = (3.4.3.12)$, then at least one vertex in $\{u_2, u_4, u_5\}$ has the face-sequence $(3.4.3.12)$. If $f_{seq}(u_4)$ or $f_{seq}(u_5)$ is $(3.4.3.12)$, then $f_{seq}(u_6) \neq f_1, f_2$. Similarly, if $f_{seq}(u_2)$ is $(3.4.3.12)$, then $f_{seq}(u_1) \neq f_1, f_2$.

\smallskip

\noindent{\bf Case 3:} Let $((3^2.6^2), f_2) \in A$ and $u$ be a vertex with $f_{seq}(u) = (3^2.6^2)$ such that ${\rm lk}(u) = C_{10}(u_1,\boldsymbol{u_2},\boldsymbol{u_3},\boldsymbol{u_4},u_5, u_6,u_7,\boldsymbol{u_8},\boldsymbol{u_9}, \boldsymbol{u_{10}})$. Then we have the following cases.

\noindent{\bf Subcase 3.1:} Suppose $f_2 \in \{ (3^2.4.3.4), (3^2.4.12)\}$. Then $f_{seq}(u_6) = f_2$. This gives that $f_{seq}(u_5) \neq f_1, f_2$.


\noindent{\bf Subcase 3.2:} Suppose $f_2 = (6^3)$. Then at least one of the vertex in the set $\{u_8, u_9,u_{10}, u_1,u_2, u_3,u_4\}$ has the face-sequence $(6^3)$. This implies $f_{seq}(u_4) = (6^3)$ or $(3^2.6^2)$. In case, $f_{seq}(u_4) = (6^3)$, then $f_{seq}(u_5) = (3^2.6^2)$. This gives that $f_{seq}(u_6) \neq f_1,f_2$. On the other hand, when $f_{seq}(u_4) = (3^2.6^2)$, then $f_{seq}(u_5) \neq f_1,f_2 $.

\smallskip

\noindent{\bf Case 4:} Let $((3^4.6), f_2) \in A$ and $u$ be a vertex with $f_{seq}(u) = (3^4.6)$ such that ${\rm lk}(u) = C_8(\boldsymbol{u_1}, \boldsymbol{u_2},\boldsymbol{u_3}$, $u_4, u_5,u_6,u_7,u_8)$. Then we have the following cases.

\noindent{\bf Subcase 4.1:} Suppose $f_2 = \{(3^2.4.3.4), (3^2.4.12)\}$. Then at least one vertex in $\{u_5,u_6, u_7\}$ has the face-sequence $f_2$. If $f_{seq}(u_i) = f_2$, where $5 \leq i \leq 7$, then $f_{seq}(u_4)$ or $f_{seq}(u_8)$ is neither $f_1$ nor $f_2$.



\noindent{\bf Subcase 4.2:} Suppose $f_2 \in \{(6^3), (4.6.12) \}$. Then at least one in $\{u_1,u_2, u_8\}$ has the face-sequence $f_2$. If $f_{seq}(u_1)$ or $f_{seq}(u_2)$ or $f_{seq}(u_8)$ is $f_2$, then $f_{seq}(u_3)$ or $f_{seq}(u_7)$ is neither $f_1$ nor $f_2$.

\smallskip

\noindent{\bf Case 5:} Let $\{f_1 = (3^2.4.3.4), f_2\} \in A$ and $u$ be a vertex with $f_{seq}(u) = (3^2.4.3.4)$ such that ${\rm lk}(u) = C_7(u_1, \boldsymbol{u_2}, u_3, u_4, u_5, \boldsymbol{u_6}, u_7)$. Then we have the following cases.

\noindent{\bf Subcase 5.1:} Suppose $f_2 \in \{ (4^4),  (4.8^2)\}$. Then at least one vertex in $\{u_2, u_6\}$ has the face-sequence $f_2$. If $f_{seq}(u_2)$ or $f_{seq}(u_6)$ is $f_2$, then $f_{seq}(u_3)$ or $f_{seq}(u_{8})$ is neither $f_1$ nor $f_2$.

\noindent{\bf Subcase 5.2:} Suppose $f_2 = (3^2.4.12)$. Then at least one vertex in $\{u_2, u_4, u_6\}$ has the face-sequence $f_2$. If $f_{seq}(u_{2i}) = f_2$, where $2 \leq i \leq 3$, then $f_{seq}(u_3)$ or $f_{seq}(u_5)$ is neither $f_1$ nor $f_2$.

\noindent{\bf Subcase 5.3:} Suppose $f_2 = (3.4.3.12)$. Then $u_3,u_4,u_5$ has face-sequence $(3^2.4.3.4)$ and at least one vertex in $\{u_1,u_2,u_7\}$ has the face-sequence $f_2$. If $f_{seq}(u_1)$ or $f_{seq}(u_2)$ or $f_{seq}(u_7)$ is $f_2$, the we get a vertex in ${\rm lk}(u_1)$ whose face-sequence is neither $f_1$ nor $f_2$.

\smallskip

\noindent{\bf Case 6:} Let $((3.6.3.6), f_2) \in A$ and $u$ be a vertex with $f_{seq}(u) = (3.6.3.6)$ such that ${\rm lk}(u) = C_{10}(\boldsymbol{u_1},\boldsymbol{u_2},\boldsymbol{u_3},u_4, u_5,\boldsymbol{u_6},\boldsymbol{u_7}, \boldsymbol{u_{8}},u_9,u_{10})$. Suppose $f_2 \in  \{(6^3), (4.6.12)\}$. Then at least one vertex $\{u_1,u_2,u_3,u_6, u_7,u_8\}$ has the face-sequence $f_2$. If $f_{seq}(u_1)$ or $f_{seq}(u_2)$ or $f_{seq}(u_3)$ is $f_2$, then $f_{seq}(u_4)$ or $f_{seq}(u_{10})$ is neither $f_1$ nor $f_2$. Similarly, if $f_{seq}(u_6)$ or $f_{seq}(u_7)$ or $f_{seq}(u_8)$ is $f_2$, then $f_{seq}(u_5)$ or $f_{seq}(u_{9})$ is neither $f_1$ nor $f_2$.


\smallskip

\noindent{\bf Case 7:} Let $((4^4), f_2) \in A$ and $u$ be a vertex with $f_{seq}(u) = (4^4)$ such that ${\rm lk}(u) = C_8(u_1, \boldsymbol{u_2}, u_3, \boldsymbol{u_4}$, $u_5, \boldsymbol{u_6}, u_7,\boldsymbol{u_8})$.  Suppose $f_2 \in \{(3.4.6.4), (3^2.4.12),  (4.5.20), (4.6.12), (3.4.3.12),  (4.8^2) \}$. Then at least one vertex in $\{u_2,u_4,u_6,u_8\}$ has the face-sequence $f_2$. If $f_{seq}(u_2) = f_2$, then $f_{seq}(u_1) \neq f_1,f_2$. Similarly, $f_{seq}(u_{2i}) \neq f_2$, where $2 \leq i \leq 4$.



\smallskip

\noindent{\bf Case 8:} Let $((3.4.6.4), f_2) \in A$ and $u$ be a vertex with $f_{seq}(u) = (3.4.6.4)$ such that ${\rm lk}(u) = C_9(u_1, \boldsymbol{u_2}, \boldsymbol{u_3}, \boldsymbol{u_4}, u_5, \boldsymbol{u_6}, u_7, u_8, \boldsymbol{u_9})$. Then we have the following cases.

\noindent{\bf Subcase 8.1:} If $f_2 \in \{(3^2.4.12), (3.4.3.12)\}$, then at least one vertex in $\{u_6, u_7, u_8, u_9\}$ has the face-sequence $f_2$. If $f_{seq}(u_6) = f_2$, then $f_{seq}(u_5)$ or $f_{seq}(u_7)$ is neither $f_1$ or $f_2$.  Following similar argument for $f_{seq}(u_i) = f_2$, where $7\leq i \leq 9$, we get a vertex in ${\rm lk}(u_i)$ whose face-sequence is neither $f_1$ nor $f_2$.

\noindent{\bf Subcase 8.2:} If $f_2 = (4.8^2)$, then at least one vertex in $\{u_6, u_9\}$ has the face-sequence $(4.8^2)$. If $f_{seq}(u_6) = (4.8^2)$ or $f_{seq}(u_9) = (4.8^2)$ then $f_{seq}(u_5) \neq f_1, f_2$ or $f_{seq}(u_1) \neq f_1, f_2$. 

\noindent{\bf Subcase 8.3:} If $f_2 = (6^3)$, then at least one vertex in $\{u_2, u_3, u_4\}$ has the face-sequence $(6^3)$. But for each $f_{seq}(u_i) = (6^3)$, where $2 \leq i \leq 4$, we see either $f_{seq}(u_1) \neq f_1, f_2$ or $f_{seq}(u_5) \neq f_1, f_2$.  

\noindent{\bf Subcase 8.4:} If $f_2 = (4.5.20)$, then at least one vertex in $\{u_6, u_9\}$ has the face-sequence $(4.5.20)$. But for each $f_{seq}(u_i) = (4.5.20)$, where $i \in \{6,9\}$, we see either $f_{seq}(u_1) \neq f_1, f_2$ or $f_{seq}(u_5) \neq f_1, f_2$.


\smallskip

\noindent{\bf Case 9:} Let $((3^2.4.12), f_2) \in A$ and $u$ be a vertex with $f_{seq}(u) = (3^2.4.12)$ such that ${\rm lk}(u) = C_{14}(u_1, \boldsymbol{u_2}, \boldsymbol{u_3}, \boldsymbol{u_4}, \boldsymbol{u_5}, \boldsymbol{u_6}, \boldsymbol{u_7}, \boldsymbol{u_8}, \boldsymbol{u_{10}}, u_{11}, \boldsymbol{u_{12}}, u_{13}, u_{14})$. Then we have the following cases.

\noindent{\bf Subcase 9.1:} If $f_2 \in \{(4.5.20),(4.8^2)\}$, then $f_{seq}(u_{12})=f_2$. This implies $f_{seq}(u_{11}) \neq f_1, f_2$.  

\noindent{\bf Subcase 9.2:} If $f_2 = (3.12^2)$, then at least one vertex in $\{u_1, u_2, \ldots, u_{10}\}$ has the face-sequence $(3.12^2)$ and $f_{seq}(u_{14}) = (3^2.4.12)$. This gives that $f_{seq}(u_{13})$ or $f_{seq}(u_{1})$ is neither $f_1$ nor $f_2$.


\noindent{\bf Subcase 9.3:} If $f_2 = (4.6.12)$, then at least one vertex in $\{u_2, u_3, \ldots, u_{12}\}$ has the face-sequence $(4.6.12)$ and $f_{seq}(u_1) = (3^2.4.12)$. This gives that $f_{seq}(u_{14}) \neq f_1, f_2$.   

\noindent{\bf Subcase 9.4:} If $f_2 \in  \{(3.4.3.12)$, then at least one vertex in $\{u_1, u_2, \ldots, u_{13}\}$ has the face-sequence $(3.4.3.12)$ and $f_{seq}(u_{14}) = (3^2.4.12)$.  This gives that $f_{seq}(u_{13})$ or $f_{seq}(u_{1})$ is neither $f_1$ nor $f_2$.  

\smallskip 

\noindent{\bf Case 10:} Let $((4.8^2), f_2) \in A$ and $u$ be a vertex with $f_{seq}(u) = (4.8^2)$ such that ${\rm lk}(u) = C_{14}(u_1, \boldsymbol{u_2}, \boldsymbol{u_3}, \boldsymbol{u_4}, \boldsymbol{u_5}, \boldsymbol{u_6}, u_7, \boldsymbol{u_8}, \boldsymbol{u_{9}}, \boldsymbol{u_{10}}, \boldsymbol{u_{11}}, \boldsymbol{u_{12}}, u_{13}, \boldsymbol{u_{14}})$. Then for each $f_2 \in \{(4.5.20)$, $(4.6.12)$, $(3.4.3.12)\}$, we see that $f_{seq}(u_{14})=f_2$ and this implies $f_{seq}(u_1) \neq f_1, f_2$ or $f_{seq}(u_{13})\neq f_1, f_2$. 

\smallskip  

Similarly, a small computation shows that $M$ does not exists for the pairs $\{((3.12^2), (3,4,3,12))$, $((6^3), (4.6.12))$, $((5^2.10), (4.5.20))\}$. Thus the claim \ref{c1}.


\smallskip

Now for the pairs $(f_1, f_2) \in B$, we have Claims \ref{c2}, \ref{c3}, \ref{c4}, \ref{c5},  \ref{c6},  \ref{c7}, \ref{c8}, \ref{c9}. 

\begin{claim}\label{c2}
	For the pairs $((3^3.4^2), f_2) \in B$, where $f_2 \in \{ (3^6), (4^4), (3.4.6.4), (3^2.4.3.4) \}$, we get DSEMs $N_1, N_2, \ldots, N_7$ of respective types $T_1, T_2, \ldots, T_7$ on the plane, as shown in example Sec. \ref{s3}.
\end{claim}

Let $u$ be a vertex with  $f_{seq}(u) =f_1= (3^3.4^2)$ and  ${\rm lk}(u) = C_7(\boldsymbol{u_1},u_2,u_3,u_4,u_5,\boldsymbol{u_6},u_7)$. Then we have the following cases.

\smallskip 

{\bf Case 1:} If $f_2 = (3^6)$, then we see that $f_{seq}(u_1) = f_{seq}(u_2) = f_{seq}(u_5) = f_{seq}(u_6) = f_{seq}(u_7) = (3^3.4^2)$ and therefore $u_{3},u_4$ can have face-sequence $(3^3.4^2)$ or $(3^6)$. Note that, at least one vertex in $\{u_3, u_4\}$ has face-sequence $(3^6)$. If $f_{seq}(u_3) = (3^6)$ (or $f_{seq}(u_4) = (3^6)$), then $f_{seq}(u_4) = (3^6)$ (resp. $f_{seq}(u_3) = (3^6)$). Thus $f_{seq}(u_3) = f_{seq}(u_4) = (3^6)$ and this implies the vertices with face-sequence $f_1 = (3^3.4^2)$ can have links of the face-sequence $({f_1}^5.{f_2}^2)$. Now we find the face-sequence of links of those vertices which have face-sequence $(3^6)$. For this, consider ${\rm lk}(u_3) = C_6(u,u_2,u_8,u_9,u_{10},u_4)$, then $f_{seq}(u_9)$ is either $(3^3.4^2)$ or $(3^6)$. In case, $f_{seq}(u_9) = (3^3.4^2)$, we have $f_{seq}(u_{10})=(3^3.4^2)$ and $f_{seq}(u_8)=(3^6)$. This gives the vertices with face-sequence $ f_2 = (3^6)$ can have links of the face-sequence $({f_1}^2.{f_2}. {f_1}^2.{f_2})$. On the other hand, when $f_{seq}(u_9)=(3^6)$, then $f_{seq}(u_8)=f_{seq}(u_{4})=f_{seq}(u_{10}) = (3^6)$ and this gives the vertices with face-sequence $(3^6)$ can have links of the face-sequence $({f_1}^2.{f_2}^4)$. Thus, there are two possibilities for the face-sequence of links of vertices with face-sequence $(3^6)$. By Lemma \ref{l1}, constructing maps for both the possibilities, we get $M=N_1$ of type $T_1 = [(3^6)^{({f_1}^2.{f_2}. {f_1}^2.{f_2})}:(3^3.4^2)^{({f_1}^5.{f_2}^2)}]$ and $M=N_2$ of type $T_2= [(3^6)^{({f_1}^2.{f_2}^4)}:(3^3.4^2)^{({f_1}^5.{f_2}^2)}]$. 

\smallskip

{\bf Case 2:} If $f_2 = (4^4)$, then we see that $f_{seq}(u_2) = f_{seq}(u_3) = f_{seq}(u_4) = f_{seq}(u_5) = (3^3.4^2)$ and $f_{seq}(u_7) = (4^4)$. This implies $f_{seq}(u_{6}) = f_{seq}(u_{1})=(4^4)$ and thus the vertices with face-sequence $(3^3.4^2)$ can have links of face-sequence $({f_1}^4.{f_2}^3)$. Now we find the face-sequence of links of those vertices with have face-sequence $(4^4)$. For this, consider ${\rm lk}(u_7) = C_8(\boldsymbol{u_{10}},u_1,\boldsymbol{u_{2}},u,\boldsymbol{u_5},u_6,\boldsymbol{u_8},u_9)$. Then $f_{seq}(u_9)$ is either $(3^3.4^2)$ or $(4^4)$. In the first case, when $f_{seq}(u_9) = (3^3.4^2)$, then $f_{seq}(u_8) = f_{seq}(u_{10}) = (3^3.4^2)$ and this shows that the vertices with face-sequence $(4^4)$ can have links of the face-sequence  $({f_1}^3.{f_2}.{f_1}^3.{f_2})$. 
On the other hand, when $f_{seq}(u_9) = (4^4)$, then $f_{seq}(u_8) = f_{seq}(u_{10}) = (4^4)$ and we see that the vertices with face-sequence $(4^4)$ can have links of the face-sequence  $({f_1}^3.{f_2}^5)$. Thus, we have two possibilities for the links of face-sequence $(4^4)$. By Lemma \ref{l1}, constructing $M$ on the plane for the pair $((3^3.4^2),(4^4))$, we get $N_3$ of type $T_3 = [(3^3.4^2)^{({f_1}^4.{f_2}^3)}:(4^4)^{({f_1}^3.{f_2}.{f_1}^3.{f_2} )}]$ and $N_4$ of type $T_4= [(3^3.4^2)^{({f_1}^4.{f_2}^3)}:(4^4)^{({f_1}^3.{f_2}^5)}]$. 

\smallskip

{\bf Case 3:}  If $f_2 = (3.4.6.4)$, then $f_{seq}(u_3) = f_{seq}(u_4) = f_{seq}(u_7) =(3^3.4^2)$. Now, we show:

\begin{itemize}
	\item The vertices $u_1,u_2,u_5,u_6$ can not have face-sequence 
	$(3^3.4^2)$.
\end{itemize} 

Suppose, if possible, $f_{seq}(u_5) = (3^3.4^2)$, then $f_{seq}(u_2) = (3^3.4^2)$.  This implies $f_{seq}(u_1) = f_{seq}(u_6) = (3^3.4^2)$, which shows that each vertex in ${\rm lk}(u)$ has the face-sequence $(3^3.4^2)$. But this can not be true. Hence, $f_{seq}(u_5) \neq (3^3.4^2)$. Similarly, we see that the vertices $u_1,u_2,u_6$ can not have face-sequence $(3^3.4^2)$ and therefore $f_{seq}(u_1) = f_{seq}(u_2) = f_{seq}(u_5) = f_{seq}(u_6) = (3.4.6.4)$. This shows, the vertices with face-sequence $(3^3.4^2)$ can have links of face-sequence $({f_2}^2.{f_1}.{f_2}^2.{f_1}^2)$. Now to find the face-sequence of links of those vertices which have face-sequence $(3.4.6.4)$, consider 
${\rm lk}(u_5) = C_9(\boldsymbol{u_8},u_9,\boldsymbol{u_{10}},\boldsymbol{u_{11}}, \boldsymbol{u_{12}},u_6,\boldsymbol{u_{7}},u,u_4)$. As $f_{seq}(u_4) = (3^3.4^2)$, we get $f_{seq}(u_8) = (3^3.4^2)$ and therefore $f_{seq}(u_9) = f_{seq}(u_{10}) = f_{seq}(u_{11}) = f_{seq}(u_{12}) = f_{seq}(u_{6}) = (3.4.6.4)$. This shows, the vertices with face-sequence $(3.4.6.4)$ can have links of face-sequence  $({f_2}^5.{f_1}^4)$. By Lemma \ref{l1}, for the pair $((3^3.4^2),(3.4.6.4))$, we get $N_5$ of type $T_5$ = $[(3^3.4^2)^{({f_2}^2. {f_1}.{f_2}^2.{f_1}^2)}:(3.4.6.4)^{({f_2}^5.{f_1}^4)}]$. 

\smallskip

{\bf Case 4:} If $f_2= (3^2.4.3.4)$, then $f_{seq}(u_7) = (3^3.4^2)$. For the exitence of $M$, we need some $u_i$, for $( 1 \leq i \leq 6)$, with face-sequence $f_2$. Without loss of generality, let $f_{seq}(u_2) = (3^2.4.3.4)$. Then ${\rm lk}(u_2) = C_{7}(u,\boldsymbol{u_7},u_1,u_8,\boldsymbol{u_9},u_{10},u_3)$ or ${\rm lk}(u_2) = C_{7}(u,\boldsymbol{u_7},u_1,u_8,u_9,\boldsymbol{u_{10}},u_3)$. 

{\bf Subcase 4.1:} If ${\rm lk}(u_2) = C_{7}(u,\boldsymbol{u_7},u_1,u_8,\boldsymbol{u_9},u_{10},u_3)$, then $f_{seq}(u_3)=(3^3.4^2)$, $f_{seq}(u_{10}) = f_{seq}(u_1) = (3^2.4.3.4)$. Now, we show:

\begin{itemize}
	\item The vertices $u_4,u_5,u_6$ can not have face-sequence 
	$(3^3.4^2)$.
\end{itemize} 

Without loss of generality, let $f_{seq}(u_4) = (3^3.4^2)$. Then $f_{seq}(u_5)$ is $(3^3.4^2)$ or $(3^2.4.3.4)$. If $f_{seq}(u_5) = (3^3.4^2)$, then $f_{seq}(u_6) = (3^3.4^2)$. This shows that ${\rm lk}(u)$ has five consecutive vertices with face-sequence $(3^3.4^2)$ while ${\rm lk}(u_3)$ has only four. This can not be true, as $f_{seq}(u)=f_{seq}(u_3)$. Similarly, if $f_{seq}(u_5)=(3^2.4.3.4)$, then $f_{seq}(u_6)=(3^2.4.3.4)$ and we see that ${\rm lk}(u)$ has two consecutive vertices with face-sequence $(3^3.4^2)$ while ${\rm lk}(u_3)$ has four. Hence $f_{seq}(u_4) \neq (3^3.4^2)$. A similar calculation shows that the vertices $u_5,u_6$ can not have face-sequence $(3^3.4^2)$. Therefore, $f_{seq}(u_4) = f_{seq}(u_5) = f_{seq}(u_6) = (3^2.4.3.4)$. This shows the vertices with face-sequence $(3^3.4^2)$ can have links of face-sequence $(f_1.{f_2}^3.f_1.{f_2}^2)$. Now, we need to find the face-sequence of links of those vertices which have face-sequence $(3^2.4.3.4)$. For this we consider ${\rm lk}(u_1)$. This gives ${\rm lk}(u_1) = C_7(\boldsymbol{u},u_7,u_{11},u_{12},\boldsymbol{u_{13}},u_8,u_2)$ or ${\rm lk}(u_1) = C_7(\boldsymbol{u},u_7,u_{11},\boldsymbol{u_{12}},u_{13},u_8,u_2)$. If ${\rm lk}(u_1) = C_7(\boldsymbol{u},u_7,u_{11},u_{12},\boldsymbol{u_{13}},u_8,u_2)$, then $f_{seq}(u_2) = f_{seq}(u_8) = (3^2.4.3.4)$ but ${\rm lk}(u_2)$ has four consecutive vertices with face-sequence $(3^2.4.3.4)$ while ${\rm lk}(u_8)$ has at least five. So, ${\rm lk}(u_1) = C_7(\boldsymbol{u},u_7,u_{11},\boldsymbol{u_{12}}$, $u_{13},u_8,u_2)$. Now proceeding similarly, as above, we see that $f_{seq}(u_{13}) = f_{seq}(u_8) = f_{seq}(u_9) = (3^3.4^2)$. This shows, the vertices with face-sequence $f_2=(3^2.4.3.4)$ can have links of face-sequence $( {f_1}^3.f_2.{f_1}^2.f_2)$. 

{\bf Subcase 4.2:} If ${\rm lk}(u_2) = C_{7}(u,\boldsymbol{u_7},u_1,u_8,u_9,\boldsymbol{u_{10}},u_3)$, then proceeding similarly, as in Subcase 4.1, we see that the vertices with face-sequences $(3^3.4^2)$ can have links of face-sequence $(f_1.{f_2}^6)$ and the vertices with face-sequences $(3^2.4.3.4)$ can have links of face-sequence  $( {f_1}^2.{f_2}^3.{f_1}.f_2)$. 

Now, by Lemma \ref{l1}, constructing $M$ for the pair $((3^3.4^2),(3^2.4.3.4))$ we get $N_6$ of type $T_6 = [(3^3.4^2)^{( {f_1}.{f_2}^6)}:(3^2.4.3.4)^{( {f_1}^2.{f_2}^3.{f_1}.f_2)}]$ and $N_7$ of type $T_7 = [(3^3.4^2)^{(f_1.{f_2}^3.f_1.{f_2}^2)}:(3^2.4.3.4)^{( {f_1}^3.f_2.{f_1}^2.f_2)}]$. 

The Claim \ref{c2} follows form Cases 1-4.

\smallskip

\begin{claim}\label{c3}
	For the pairs $((3^6), f_2) \in B$, $f_2 \in \{ (3^2.6^2),  (3^2.4.12), (3^2.4.3.4), (3^4.6) \}$, we get DSEMs $N_8, N_9, \ldots, N_{13}$ of respective types $T_8, \ldots T_{13}$ on the plane, as shown in example Sec. \ref{s3}.
\end{claim}

Let $u$ be a vertex with $f_{seq}(u)=f_1 = (3^6)$ and ${\rm lk}(u) = C_6(u_1,u_2,u_3,u_4,u_5,u_6)$. Then we have the following cases. 

\smallskip

{\bf Case 1:} If $f_2 = (3^2.6^2)$, then at least one vertex in ${\rm lk}(u)$ has the face-sequence $(3^2.6^2)$. Without loss of generality, let $f_{seq}(u_1) = (3^2.6^2)$. This implies, $f_{seq}(u_2) = f_{seq}(u_3) = f_{seq}(u_4) = f_{seq}(u_5) = f_{seq}(u_6) = (3^2.6^2)$. This shows that the vertices with face-sequence $(3^6)$ can have links of face-sequence  $({f_2}^6)$.  Observe that ${\rm lk}(u_1) = C_{10}(u,u_6,\boldsymbol{u_7},\boldsymbol{u_8},\boldsymbol{u_9},u_{10},\boldsymbol{u_{11}},\boldsymbol{u_{12}},\boldsymbol{u_{13}},u_2)$, this gives that the vertices $u_7,u_8,u_9,u_{10},u_{11},u_{12},u_{13}$ have face-sequence $(3^2.6^2)$. Thus, the vertices with face-sequence $(3^2.6^2)$ can have links of face-sequence $({f_1}.{f_2}^9)$. By Lemma \ref{l1}, we get $N_8$ of type $T_8 = [(3^6)^{({f_2}^6)}:(3^2.6^2)^{(f_1.{f_2}^9)}]$. 

\smallskip

{\bf Case 2:} If $f_2 = (3^2.4.12)$, then proceeding similarly, as above, we see that the vertices with face-sequence $(3^6)$ can have links of face-sequence $({f_2}^6)$ and the vertices with face-sequence $(3^2.4.12)$ can have links of face-sequence $(f_1.{f_2}^{13})$. By Lemma \ref{l1}, we get $N_9$ of type $T_9 = [(3^6)^{({f_2}^6)}:(3^2.4.12)^{({f_1}.{f_2}^{13})}]$. 

\smallskip

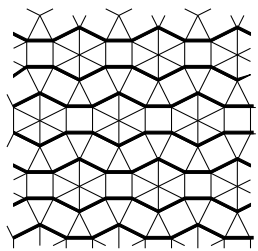
{\bf Case 3:} If $f_2 = (3^2.4.3.4)$, then proceeding similarly, as in Case 1, the vertices with face-sequence $(3^6)$ can have links of face-sequence $(f_2^6)$ and the vertices with face-sequence $(3^2.4.3.4)$ can have links of face-sequence  $(f_1.{f_2}^{6})$. By Lemma \ref{l1}, we get $N_{10}$ of type $T_{10} = [(3^6)^{({f_2}^6)}:(3^2.4.3.4)^{(f_1.{f_2}^{6})}]$. 

\smallskip

{\bf Case 4:} Let  $v$ be a vertex with face-sequence $f_2 = {(3^4.6)}$ and ${\rm lk}(v) = C_{8}(\boldsymbol{v_1},\boldsymbol{v_2},\boldsymbol{v_3},v_4,v_5,v_6,v_7$, $v_8)$. Then $v_1,v_2,v_3,v_4,v_8$ have face-sequence ${(3^4.6)}$. This implies at least one vertex in $\{v_5,v_6,v_7\}$ has the face-sequence $(3^6)$. Without loss of generality, let $f_{seq}(v_5) = (3^6)$. Then ${\rm lk}(v_5) = C_{6}(v,v_4,v_9$, $v_{10}, v_{11},v_6)$, which gives $f_{seq}(v_6)$ is ${(3^4.6)}$ or $(3^6)$. 

\smallskip

{\bf Subcase 4.1:} If $f_{seq}(v_6) = {(3^4.6)}$, then either ${\rm lk}(v_6) = C_{8}(\boldsymbol{v_{12}},\boldsymbol{v_{13}},\boldsymbol{v_{14}},v_{15},v_7,v,v_5, v_{11})$ or  ${\rm lk}(v_6) = C_{8}({v_{12}},\boldsymbol{v_{13}},\boldsymbol{v_{14}},\boldsymbol{v_{15}},v_7,v,v_5,v_{11})$.
In case ${\rm lk}(v_6) = C_{8}(\boldsymbol{v_{12}},\boldsymbol{v_{13}},\boldsymbol{v_{14}},v_{15},v_7,v,v_5, v_{11})$, we get $f_{seq}(v_{11}) = f_{seq}(v_{12}) = f_{seq}(v_{13}) = f_{seq}(v_{14}) = f_{seq}(v_{15}) = {(3^4.6)}$ and hence $f_{seq}(v_7) = (3^6)$.  This shows the vertices with face-sequence ${(3^4.6)}$ can have links of face-sequence  $({{f_2}}^5.{f_1}.{{f_2}}.{f_1})$. Similarly, we see that the vertices with face-sequence $(3^6)$ can have links of face-sequence $({f_2}^6)$. On the other hand, when ${\rm lk}(v_6) = C_{8}({v_{12}},\boldsymbol{v_{13}},\boldsymbol{v_{14}},\boldsymbol{v_{15}},v_7,v,v_5,v_{11})$, then $f_{seq}(v_7) = (3^4.6)$. This shows that the vertices with face-sequence $(3^4.6)$ can have link of face-sequence $({{f_2}}^7.{f_1})$ and the vertices with face-sequence $(3^6)$ can have links of face-sequence  $({{f_2}}^6)$.  
	
	\smallskip
	
{\bf Subcase 4.2:}  If $f_{seq}(v_6) = (3^6)$, then $f_{seq}(v_7)$ is ${(3^4.6)}$ or $(3^6)$. If $f_{seq}(v_7) = {(3^4.6)}$, then $v_6$ is the only vertex with face-sequence $(3^6)$ in  ${\rm lk}(v_7)$ while there are two vertices having face-sequence $(3^6)$ in ${\rm lk}(v)$. This can not be possible. So $f_{seq}(v_7) = (3^6)$. This shows that the vertices with face-sequence ${(3^4.6)}$ can have links of face-sequence $({{f_2}}^5.{f_1}^3)$. 
	
Note that $f_{seq}(v_{10})$ is ${(3^4.6)}$ or $(3^6)$. In case $f_{seq}(v_{10}) = (3^6)$, then $f_{seq}(v_{11})$ is $(3^6)$ or ${(3^4.6)}$. If, $f_{seq}(v_{11}) = (3^6)$, we get $f_{seq}(u_5) = f_{seq}(u_{6}) = (3^6)$ but $u_5$ has face-sequence of its link as $({{f_1}}^4.{f_2}^2)$ while $u_6$ have face-sequence of its link as $({{f_1}}^5.{f_2})$ or ${(3^6)}^{({{f_1}}^2.{f_2}^2.{f_1}.{f_2})}$. This is not possible. So $f_{seq}(v_{11})={(3^4.6)}$. Then we see that the vertices with face-sequence $(3^6)$ can have links of face-sequence $({{f_1}}^2.{f_2}^2.{f_1}.{f_2})$. On the other hand, when $f_{seq}(v_{10}) = (3^4.6)$, then either ${\rm lk}(v_{10}) = C_{8}({v_{12}},v_{13},\boldsymbol{v_{14}},\boldsymbol{v_{15}},\boldsymbol{v_{16}},v_{11},v_5,v_{9})$ or  ${\rm lk}(v_{10}) = C_{8}(\boldsymbol{v_{13}},\boldsymbol{v_{14}},\boldsymbol{v_{15}},v_{16},v_{11},v_5,v_{9},v_{12})$.
	
If ${\rm lk}(v_{10}) = C_{8}({v_{12}},v_{13},\boldsymbol{v_{14}},\boldsymbol{v_{15}},\boldsymbol{v_{16}},v_{11},v_5,v_{9})$, then $f_{seq}(v_5) = f_{seq}(v_9) = (3^6)$ but $v_5$ has face-sequence of its link as $({{f_1}}.{f_2}^2.f_1.{f_2}^2)$ while $u_9$ has face-sequence of its link as  $({{f_1}}^3.{f_2}.f_1.f_2)
$ or $({f_2}^3.{f_1}.{f_2}_1)$. This is not possible. So ${\rm lk}(v_{10}) = C_{8}(\boldsymbol{v_{13}},\boldsymbol{v_{14}},\boldsymbol{v_{15}},v_{16},v_{11},v_5,v_{9},v_{12})$. Now, proceeding similarly, we see easily that the vertices with face-sequence $(3^6)$ can have links of face-sequence $({{f_1}}^2.{f_2}^2.f_1.{f_2})$. 
	
By Lemma \ref{l1}, constructing $M$ for the pair $\{3^6,3^4.6\}$, we get $N_{11}$ type $T_{11} = [{(3^6)}^{({{f_1}}^2.{f_2}^2.{f_1}.{f_2})}:{(3^4.6)}^{({{f_2}}^5.{f_1}^3)}]$, $N_{12}$ of type $T_{12} = [(3^6)^{{(f_2)}^6}:{(3^4.6)}^{({{f_2}}^5.{f_1}.{{f_2}}.{f_1})}]$,  and $N_{13}$ of type $T_{13} = [(3^6)^{{{(f_2)}}^6}:{(3^4.6)}^{({{f_2}}^7.{f_1})}]$. 

Now, the Claim \ref{c3} follows from Cases 1-4.


\begin{claim}\label{c4}
	For the pairs $((3.4.6.4), f_2) \in B$, where $f_2 \in \{ (3.4^2.6), (3^2.4.3.4), (4.6.12)\}$, we get DSEMs $N_{14}, N_{15}, N_{16}$ of respective types $T_{14}, T_{15}, T_{16}$ on the plane,  as shown in example Sec \ref{s3}.
\end{claim}

Let $v$ be a vertex with $f_{seq}(v)=f_1 = (3.4.6.4)$ and let ${\rm lk}(v)=C_{9}(\boldsymbol{v_1},v_2,\boldsymbol{v_3},v_4,v_5,\boldsymbol{v_6},v_7,\boldsymbol{v_8}, \boldsymbol{v_9})$. Then we have the following cases. 

{\bf Case 1:} If $f_2 = (3.4^2.6)$, then at least one vertex in ${\rm lk}(v)$ has the face-sequence $(3.4^2.6)$. Now we show:

\begin{itemize}
	\item The vertices  $v_1,v_2,v_7,v_8,v_9$ can not have the face-sequence $(3.4^2.6).$ 
\end{itemize}

Suppose, if possible, $f_{seq}(v_2) = (3.4^2.6)$. Then $f_{seq}(v_3) = f_{seq}(v_1) = (3.4^2.6)$. This gives ${\rm lk}(v_2)=C_{9}(v_1,v_{10},\boldsymbol{v_{11}},v_3,\boldsymbol{v_4},v,\boldsymbol{v_7},\boldsymbol{v_8},\boldsymbol{v_9})$. Since $f_{seq}(v_2) = f_{seq}(v_3) = f_{seq}(v_1) = (3.4^2.6)$, we get $f_{seq}(v_4) = f_{seq}(v_5) = f_{seq}(v_{10}) = f_{seq}(v_{11}) = (3.4.6.4)$. This implies $f_{seq}(v_7)$ is $(3.4^2.6)$ or $(3.4.6.4)$

Let $f_{seq}(v_7) = (3.4^2.6)$ and ${\rm lk}(v_7)=C_{9}(\boldsymbol{v_1},\boldsymbol{v_2},v,\boldsymbol{v_5},v_6,\boldsymbol{v_{12}},v_{13},v_8, \boldsymbol{v_9})$. Then $f_{seq}(v_6) = f_{seq}(v_8) = (3.4^2.6)$.
This gives $f_{seq}(v_9) = f_{seq}(v_{12}) = f_{seq}(v_{13}) = (3.4.6.4)$. 
Note that $v_1,v_2,v_3$ or $v_6,v_7,v_8$  are consecutive vertices  in  ${\rm lk}(v)$ with face-sequence $(3.4^2.6)$. While constructing ${\rm lk}(v_{13})$, we get ${\rm lk}(v_{13})=C_{9}(\boldsymbol{v_6},v_{12}, \boldsymbol{v_{14}},\boldsymbol{v_{15}}, \boldsymbol{v_{16}},v_{17},\boldsymbol{v_{18}},v_8,v_7)$, and we see four consecutive vertices $u_{18},u_8,u_7,u_6$ with face-sequence $(3.4^2.6)$. This is not possible.

Let $f_{seq}(v_7) = (3.4.6.4)$ and ${\rm lk}(v_7)=C_{9}(\boldsymbol{v_1},\boldsymbol{v_2},v,\boldsymbol{v_5},v_6,v_{12},\boldsymbol{v_{13}},v_8, \boldsymbol{v_9})$. Then $f_{seq}(v_8)$ is either $(3.4.6.4)$ or $(3.4^2.6)$. 
In case,  $f_{seq}(v_8) = (3.4.6.4)$, we get $f_{seq}(v_6) = f_{seq}(v_{12}) = f_{seq}(v_{13}) = f_{seq}(v_9) = (3.4.6.4)$. This implies the vertices $v$ and $v_7$ have links of its face-sequences $({f_2}^{3}.{f_1}^6)$ and $({f_2}^{2}.{f_1}^7)$. This is not possible. Similarly, for $f_{seq}(v_8) = (3.4^2.6)$, we get a contradiction.

A small computation shows that the vertices $v_1,v_7,v_8,v_9$ can not have the face-sequence $(3.4^2.6)$. So, the vertices $v_1,v_2,v_7,v_8,v_9$ have face-sequence $(3.4.6.4)$ and therefore at least one vertex in $\{v_3,v_4,v_5,v_6\}$ has face-sequence $(3.4^2.6)$. Without loss of generality, let $f_{seq}(v_4) = (3.4^2.6)$ with ${\rm lk}(v_4)=C_{9}(v,\boldsymbol{v_2},v_3,\boldsymbol{v_{19}},v_{20},\boldsymbol{v_{21}}, \boldsymbol{v_{22}},\boldsymbol{v_{23},v_5})$. Then $f_{seq}(v_3) = f_{seq}(v_5) = f_{seq}(v_6) = (3.4^2.6)$. This shows the vertices with face-sequence $(3.4.6.4)$ can have links of face-sequence $({f_1}^5.{f_2}^{4})$.

Also, $f_{seq}(v_{20})$ is $(3.4.6.4)$ or $(3.4^2.6)$. In case $f_{seq}(v_{20}) = (3.4.6.4)$, the face-sequence of the link of $v_{20}$ is $({f_1}^6.{f_2}^{3})$ or  $({f_1}^2.{f_2}^{3}.f_1.{f_2}^3)$. Note that $f_{seq}(v) = f_{seq}(v_{20}) = (3.4.6.4)$ but their respective face-sequence of links are different. So $f_{seq}(v_{20}) \neq (3.4.6.4)$. On the other hand, if $f_{seq}(v_{20}) = (3.4^2.6)$, then $f_{seq}(v_{19}) = f_{seq}(v_{27}) = (3.4^2.6)$. Observe that, the vertices $v_{25},v_{26}$ have face-sequence $(3.4.6.4)$ or $(3.4^2.6)$. In case, any vertex in $\{v_{25}, v_{26}\}$ has face-sequence $(3.4.6.4)$ then the other has also $(3.4.6.4)$. 
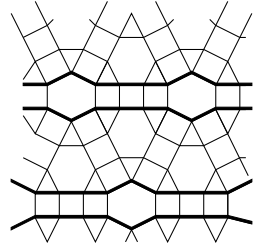
This gives ${\rm lk}(v_{25})$ has five consecutive vertices with face-sequence $(3.4^2.6)$, while ${\rm lk}(v)$ has four, which is not possible. So, $f_{seq}(v_{25}) = (3.4^2.6)$ and therefore $f_{seq}(v_{26}) = (3.4^2.6)$. This shows that the vertices with face-sequence $(3.4^2.6)$ can have links of face-sequence $({f_1}^2.{f_2}^{7})$. By Lemma \ref{l1}, we get $N_{14}$ of type $T_{14} = [(3.4.6.4)^{({f_1}^5.{f_2}^{4})}:(3.4^2.6)^{({f_1}^2.{f_2}^{7})}]$. 

\smallskip
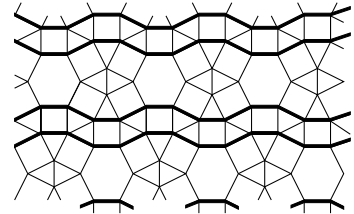
{\bf Case 2:} If $f_2 = (3^2.4.3.4)$ then ${\rm lk}(u) = C_7(u_1,\boldsymbol{u_2},u_3,u_4,u_5,\boldsymbol{u_6},u_{7})$. This implies $f_{seq}(u_4) = (3^2.4.3.4)$, $f_{seq}(u_3) = f_{seq}(u_5) = (3.4.6.4)$ and ${\rm lk}(u_4) = C_7(u_3,\boldsymbol{u_8},u_9,u_{10},\boldsymbol{u_{11}},u_{5}, u)$ and ${\rm lk}(u_5) = C_9(u_4,\boldsymbol{u_{10}},u_{11},  \boldsymbol{u_{12}}, \boldsymbol{u_{13}},\boldsymbol{u_{14}},u_{6},\boldsymbol{u_{7}},u)$. Now, we note that $u_6,u_{14}, u_{13}, u_{12},u_{11}$ have face-sequence $(3.4.6.4)$ and $u_1, u_7,u_{10},u_9$ have face-sequence $(3^2.4.3.4)$. This implies $f_{seq}(u_8) = (3.4.6.4)$. So, the vertices with face-sequence $f_2 =  (3^2.4.3.4)$ can have links of face-sequence $({f_1}^2.{f_2}.{f_1}^2.{f_2}^2)$ and the vertices with face-sequence $(3.4.6.4)$ can have links of face-sequence $({f_1}^5.{f_2}^4)$. By Lemma \ref{l1}, we get $N_{15}$ of type $T_{15} = [(3.4.6.4)^{({f_1}^5.{f_2}^4)}:(3^2.4.3.4)^{({f_1}^2.{f_2}.{f_1}^2.{f_2}^2)}]$. 

\smallskip

{\bf Case 3:} If  $f_2 = (4.6.12)$ then ${\rm lk}(u)=C_{9}(\boldsymbol{u_1},u_2, \boldsymbol{u_3},u_4,u_5,\boldsymbol{u_6},u_7,\boldsymbol{u_8},\boldsymbol{u_9})$ and we see that the vertices $u_4,u_5$ have face-sequence $(3.4.6.4)$ and at least one vertex in $\{u_1, u_2,u_3,u_6,u_7,u_8,u_9\}$ has the face-sequence $(4.6.12)$. Now we show:

\begin{itemize}
	\item The vertex $u_9$ can not have the face-sequence $(4.6.12)$.
\end{itemize}

Suppose, if possible, $f_{seq}(u_9) = (4.6.12)$, then ${\rm lk}(u_9) = C_{16}(u_1,\boldsymbol{u_2},\boldsymbol{u},\boldsymbol{u_7},u_8,\boldsymbol{u_{10}},\boldsymbol{u_{11}},\boldsymbol{u_{12}},\boldsymbol{u_{13}}$, $\boldsymbol{u_{14}},  \boldsymbol{u_{15}},\boldsymbol{u_{16}},\boldsymbol{u_{17}},\boldsymbol{u_{18}},u_{19}, \boldsymbol{u_{20}})$. This implies the vertices $u_8, u_{10},u_{11},u_{12},u_{13},u_{14},u_{15},u_{16},u_{17}$, $u_{18}, u_{19}$ have face-sequence $(4.6.12)$ and $u_7$ has face-sequence $(3.4.6.4)$. 

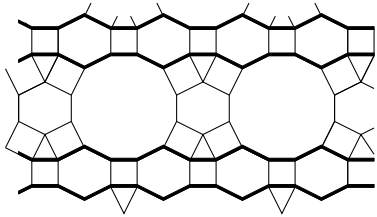
Let ${\rm lk}(u_7)=C_{9}(\boldsymbol{u_1},\boldsymbol{u_2},u,\boldsymbol{u_5},u_6,\boldsymbol{u_{21}},\boldsymbol{u_{10}}, u_8,\boldsymbol{u_{9}})$. Then $f_{seq}(u_{21}) = (3.4.6.4)$ and $f_{seq}(u_3)$ is $(3.4.6.4)$ or $(4.6.12)$. In case, $f_{seq}(u_3) = (3.4.6.4)$, we see that the vertices $u_1,u_2,u_{20}$ have face-sequence $(3.4.6.4)$. This implies the face-sequence of links of $u$ and $u_7$ are $({f_2}^{2}.{f_1}^{7})$ and $({f_1}^6.{f_2}^{3})$ respectively. This is not possible. On the other hand, if $f_{seq}(u_3) = (4.6.12)$, then the face-sequence of the links of $u_2$ and $u_9$ are $({f_2}^{13}.{f_1}^{3})$ and $({f_1}^4.{f_2}^{12})$ respectively. This is not possible. Hence, $f_{seq}(u_9) \neq (4.6.12)$. So, $f_{seq}(u_9) = (3.4.6.4)$. Then at least one vertex in $\{u_1, u_2,u_3,u_4,u_6, u_7,u_8\}$ has the face-sequence $(4.6.12)$. Without loss of generality, suppose that $f_{seq}(u_1) = (4.6.12)$. This gives $f_{seq}(u_2) = f_{seq}( u_3) = (4.6.12)$. Note that if any vertex in $\{u_6,u_7,u_8\}$ has the face-sequence $(3.4.6.4)$, then the remaining two vertices have also the face-sequence $(3.4.6.4)$. In this case, the face-sequence of links of $u$ and $u_7$ are $({f_2}^{3}.{f_1}^{6})$ and $({f_1}^7.{f_2}^{2})$ respectively. This is not possible. So, the vertices $u_6,u_7,u_8$ have face-sequence $(4.6.12)$. This shows the vertices with face-sequence $(3.4.6.4)$ and $(4.6.12)$ can have links of face-sequence $({f_2}^{3}.{f_1}^{2}.{f_2}^{3}.{f_1})$ and $({f_2}^{11}.{f_1}.{f_2}^{2}.{f_1}^2)$ respectively.  By Lemma \ref{l1}, we get $N_{16}$ of type $T_{16}$ = $[(3.4.6.4)^{({f_2}^{3}.{f_1}^{2}.{f_2}^{3}.{f_1})}:(4.6.12)^{({f_2}^{11}.{f_1}.{f_2}^{2}.{f_1}^2)}]$. 

The Claim \ref{c4} follows from Cases 1-4.


\begin{claim}\label{c5}
	The pair $((3.4.3.12),(3.12^2)) \in B$ gives DSEM $N_{17}$  of type $T_{17}$ on the plane, as shown in example Sec \ref{s3}.
\end{claim}

Let  $v$ be a vertex with face-sequence $(3.4.3.12)$. Let ${\rm lk}(v) = C_{14}(v_1,\boldsymbol{v_2},\boldsymbol{v_3},\boldsymbol{v_4},\boldsymbol{v_5},\boldsymbol{v_6},\boldsymbol{v_7},\boldsymbol{v_8},  \boldsymbol{v_9}$, $\boldsymbol{v_{10}},v_{11},v_{12},\boldsymbol{v_{13}},v_{14})$. Then $f_{seq}(v_{11}) = f_{seq}(v_{12}) = f_{seq}(v_{13}) = f_{seq}(v_{14}) = (3.4.3.12)$. This implies $f_{seq}(v_{10}) = f_{seq}(v_{2}) = (3.12^2)$ and ${\rm lk}(v_{10}) = C_{20}(v_1,\boldsymbol{v_2},\boldsymbol{v_3},  \boldsymbol{v_4},\boldsymbol{v_5},  \boldsymbol{v_6},\boldsymbol{v_7},\boldsymbol{v_8},v_9,v_{15},\boldsymbol{v_{16}},\boldsymbol{v_{17}},\boldsymbol{v_{18}}$, $\boldsymbol{v_{19}}, \boldsymbol{v_{20}}, \boldsymbol{v_{21}},\boldsymbol{v_{22}},\boldsymbol{v_{23}},\boldsymbol{v_{12}},v_{11})$. Then, we show: 

\begin{itemize}
	\item The vertices $v_3,v_9$ can not have face-sequence $(3.12^2)$.
\end{itemize}

Suppose, if possible, $f_{seq}(v_9) = (3.12^2)$. Then $f_{seq}(v_{15}) = (3.12^2)$. Now observe that $f_{seq}(v_{10}) = f_{seq}(v_{11}) = (3.12^2)$ but $v_{10}$ has three adjacent vertices with face-sequence $(3.12^2)$ while $v_{11}$ has only one.
This can not be possible. Similarly, we see that $v_3$ can not have face-sequence $(3.12^2)$. 
So, $f_{seq}(v_3) = f_{seq}(v_9) = (3.4.3.12)$. This implies $f_{seq}(v_7) = f_{seq}(v_4) =(3.12^2)$ and $f_{seq}(v_6) = (3.4.3.12)$. Thus the vertices with face-sequence $(3.4.3.12)$ can have links of face-sequence  $F_1 = (f_2^2.f_1.f_2^2.f_1.f_2^2.f_1.f_2^2.f_1^3)$ and the vertices with face-sequence  $(3.12^2)$ can have links of face-sequence $F_2 = (f_1.f_2^2.f_1.f_2^2.f_1.f_2^2.f_1.f_2.f_1.f_2^2.f_1.f_2^2.f_1.f_2^2.f_1)$.
Now, by Lemma \ref{l1}, constructing $M$ on the plane for the pair $\{(3.12^2)$, $(3.4.3.12)\}$, we get $N_{17}$ of type $T_{17}$ = $[{(3.4.3.12)}^{F_1}:(3.12^2)^{F_2}]$. This proves the claim.



\begin{claim}\label{c6}
	The pair $((3^4.6), (3^2.6^2)) \in B$ gives DSEM $N_{18}$ of type $T_{18}$ on the plane, as shown in example Sec \ref{s3}.
\end{claim}


Let  $u$ be a vertex with $ (3^2.6^2)$ and ${\rm lk}(u) = C_{10}(\boldsymbol{u_1},\boldsymbol{u_2},u_3,\boldsymbol{u_4},\boldsymbol{u_5},\boldsymbol{u_6},u_7,u_{8},u_9,\boldsymbol{u_{10}})$. Then $f_{seq}(u_3) = (3^2.6^2)$. Now, we show: 

\begin{itemize}
	\item The vertices $u_2,u_4,u_7,u_{9}$ can not have face-sequence $(3^2.6^2)$.
\end{itemize}

Suppose, the vertex $u_4$  has face-sequence $(3^2.6^2)$. Then ${\rm lk}(u_4) = C_{10}(\boldsymbol{u_{13}},\boldsymbol{u_{14}},\boldsymbol{u_{15}},u_5,\boldsymbol{u_6},\boldsymbol{u_7},\boldsymbol{u}, \linebreak u_3, u_{11},u_{12})$, and $f_{seq}(u_{11}) = (3^4.6)$. This implies either $f_{seq}(u_2) = (3.6.y_1^{l_1}.y_2^{l_2}. \ldots y_k^{l_k}.6)$ or $f_{seq}(u_{12}) = (3.6.y_1^{l_1}.y_2^{l_2}. \ldots y_k^{l_k}.6)$. This is not be possible. Thus, $f_{seq}(u_4) \neq (3^2.6^2)$. Similarly, we see that the vertices $u_2,u_7,u_{9}$ can not have face-sequence $(3^2.6^2)$. Hence $f_{seq}(u_2) = f_{seq}(u_4) = f_{seq}(u_7) = f_{seq}(u_{9}) = (3^4.6)$. Now, the vertices $u_1,u_5,u_6,u_{10}$ can have face-sequence $(3^4.6)$ or $(3^2.6^2)$. Then we show:

\begin{itemize}
	\item The vertices $u_1,u_5,u_6,u_{10}$ can not have face-sequence $(3^4.6)$.
\end{itemize}

Suppose, if possible, $f_{seq}(u_1) = (3^4.6)$. Then $f_{seq}(u_{10}) = (3^4.6)$ and ${\rm lk}(u_9) = C_{8}(\boldsymbol{u_{1}},\boldsymbol{u_{2}},\boldsymbol{u_{3}},u$, $u_8,u_{11},u_{12},u_{10})$, ${\rm lk}(u_{10}) = C_{8}({u_{1}},\boldsymbol{u_{2}},\boldsymbol{u_{3}},\boldsymbol{u},u_9,u_{12},  u_{13},u_{14})$, ${\rm lk}(u_{1}) = C_{8}({u_{2}},\boldsymbol{u_{3}},\boldsymbol{u},\boldsymbol{u_9},u_{10},u_{14}$, $u_{15},u_{16})$, ${\rm lk}(u_{2}) = C_{8}(u_1,{u_{16}},{u_{17}},{u_{18}},u_3,\boldsymbol{u},\boldsymbol{u_9}, \boldsymbol{u_{10}})$. This implies $f_{seq}(u_{12}) = f_{seq}(u_{14}) = f_{seq}(u_{16})$ $= (3^4.6)$. Now, observe that $f_{seq}(u_{11})$ is $(3^4.6)$ or $(3^2.6^2)$. If $f_{seq}(u_{11}) = (3^2.6^2)$, then ${\rm lk}(u_{11}) = C_{10}(\boldsymbol{u_{20}},\boldsymbol{u_{21}},\boldsymbol{u_{22}},u_{23}, \boldsymbol{u_{24}},\boldsymbol{u_{25}},\boldsymbol{u_{26}},u_{12},u_9,u_{19})$. This implies 
$f_{seq}(u_{23}) = (3^2.6^2)$, $f_{seq}(u_{13}) = f_{seq}(u_{15}) = (3^4.6)$, and $f_{seq}(u_{19}) = f_{seq}(u_1) = (3^4.6)$. Now, observe that ${\rm lk}(u_{19})$ has at least three vertices with face-sequence $(3^2.6^2)$ while ${\rm lk}(u_{1})$ has only two. This is not possible. 
Similarly, we see that $f_{seq}(u_{11}) \neq (3^4,6)$.  This implies $f_{seq}(u_5) = f_{seq}(u_6) = f_{seq}(u_{10}) = (3^2.6^2)$. Thus vertices with face-sequence  $(3^2.6^2)$ can have links of face-sequence $(f_2^2.f_1^3.f_2^2.f_1.f_2.f_1)$ and vertices with face-sequence $(3^4.6)$ can have links of face-sequence $(f_2^2.f_1.f_2.f_1.f_2^2.f_1)$. By Lemma \ref{l1}, constructing $M$ on the plane for the pair $\{(3^2.6^2),(3^4.6)\}$, we get $N_{18}$ of type $T_{18}$ = $[(3^4.6)^{(f_2^2.f_1.f_2.f_1.f_2^2.f_1)}:(3^2.6^2)^{(f_2^2.f_1^3.f_2^2.f_1.f_2.f_1)}]$. This proves the claim.



\begin{claim}\label{c7}
	For the pairs $((3^4.6), (3.6.3.6)) \in B$, we get DSEMs $N_{19}, N_{20}$ of respective types $T_{19},T_{20}$ on the plane, as shown in example Sec. \ref{s3}.
\end{claim}

Let $u$ be a vertex with face-sequence $(3^4.6)$ and ${\rm lk}(u)=C_{8}(\boldsymbol{u_1},\boldsymbol{u_2},\boldsymbol{u_3},u_4,u_5,u_6,u_7,u_8)$. Then $u_5,u_6$, and $u_7$
have face-sequence $(3^4.6)$ and at least one of the vertex in the set $\{u_1,u_2,u_3,u_4,u_8\}$ has face-sequence $(3.6.3.6)$. Without loss of generality, let $f_{seq}(u_5)=(3^4.6)$. Then ${\rm lk}(u_5)=C_{8}(\boldsymbol{u_{11}},\boldsymbol{u_{12}},\boldsymbol{u_{13}},u_6,u,u_4,u_9,u_{10})$ or ${\rm lk}(u_5)=C_{8}(\boldsymbol{u_{10}},\boldsymbol{u_{11}},\boldsymbol{u_{12}}, u_{13}, u_6,u,u_{4},u_9)$ or ${\rm lk}(u_5)=C_{8}(\boldsymbol{u_9}$, $\boldsymbol{u_{10}},\boldsymbol{u_{11}},u_{12},u_{13},u_6,u,u_4)$.

{\bf Case 1:} If ${\rm lk}(u_5)=C_{8}(\boldsymbol{u_{11}},\boldsymbol{u_{12}},\boldsymbol{u_{13}},u_6,u,u_4,u_9,u_{10})$, then ${\rm lk}(u_6)=C_{8}(\boldsymbol{u_{10}},\boldsymbol{u_{11}},\boldsymbol{u_{12}},u_{13},u_{14}, \linebreak u_7,u,u_{5})$ and ${\rm lk}(u_7)=C_{8}(\boldsymbol{u_{16}},\boldsymbol{u_{17}},\boldsymbol{u_{18}},u_{8},u,u_6,u_{14},u_{15})$ or ${\rm lk}(u_7)=C_{8}(\boldsymbol{u_{15}},\boldsymbol{u_{16}},\boldsymbol{u_{17}},u_{18},u_8,u, \linebreak u_6,u_{14})$.

{\bf Subcase 1.1:} If ${\rm lk}(u_7)=C_{8}(\boldsymbol{u_{15}},\boldsymbol{u_{16}},\boldsymbol{u_{17}},u_{18},u_8,u,u_6,u_{14})$, then  $f_{seq}(u_8) = (3.6.3.6)$ and $f_{seq}(u_1)$ is either $(3^4.6)$ or $(3.6.3.6)$. 

	If  $f_{seq}(u_1)=(3^4.6)$, then $f_{seq}(u_2)$ is either $(3^4.6)$ or $(3.6.3.6)$. If $f_{seq}(u_2) = (3.6.3.6)$, then $f_{seq}(u) = f_{seq}(u_1) =(3^4.6)$ but $u$ has at least four consecutive vertex having face-sequence $(3^4.6)$ while $u_1$ has only three. This can not be possible. On the other hand, if  $f_{seq}(u_2)=(3^4.6)$, then $f_{seq}(u_3)$ is either $(3^4.6)$ or $(3.6.3.6)$. If $f_{seq}(u_3)$ is $(3^4.6)$ (or $(3.6.3.6)$), then $u$ has the face-sequence of links $(f^7_1.f_2)$ (resp. $(f^4_1.f_2.f^2_1.f_2)$) and hence the vertices having the face-sequence $(3.6.3.6)$ can have links of face-sequence $(f^{10}_1)$ (resp. $(f^4_1.f_2.f^4_1.f_2)$). By Lemma \ref{l1} constructing $M$ on the plane for the pair $\{(3^4.6),(3.6.3.6)\}$, we get DSEM $N_{19}$ of type $T_{19} = [(3^4.6)^{(f^4_1.f_2.f_1^2.f_2)}:(3.6.3.6)^{(f^{4}_1.f_2.f_1^4.f_2)}]$ and DSEM $N_{20}$ of type $T_{20} = [(3^4.6)^{(f^7_1.f_2)}:(3.6.3.6)^{(f^{10}_1)}]$. 
		
    If $f_{seq}(u_1)=(3.6.3.6)$, then $f_{seq}(u_2)$ is either $(3^4.6)$ or $(3.6.3.6)$. 
	
	$(i)$ If $f_{seq}(u_2) = (3^4.6)$, then $f_{seq}(u_3)$ is $(3^4.6)$ or $(3.6.3.6)$. If $f_{seq}(u_3) = (3.6.3.6)$, then $f_{seq}(u) = f_{seq}(u_2) =(3^4.6)$ but $u$ has four consecutive vertices having face-sequence $(3^4.6)$ while $u_2$ has only three. This can not be possible. Similarly
		$f_{seq}(u_3) = (3^4.6)$ is not possible.
	
	$(ii)$ If $f_{seq}(u_2) = (3.6.3.6)$, then $f_{seq}(u_3)$ is $(3^4.6)$ or $(3.6.3.6)$. 
	If $f_{seq}(u_3) = (3^4.6)$, then the vertices having the face-sequence $(3^4.6)$  have links of face-sequence $(f^{5}_1, f^{3}_2)$. This gives that the vertex $u_2$ has face-sequence of link  $(f^{6}_1, f^{4}_2)$. Note that, both $u_1$ and $u_2$ has face-sequence $(3.6.3.6)$ but $u_1$ has three consecutive vertices having face-sequence $(3^4.6)$ while $u_2$ has six. This can not be possible. On the other hand, if $f_{seq}(u_3) = (3.6.3.6)$, then the vertices having the face-sequence $(3^4.6)$ can have links of face-sequence $(f^{4}_1, f^{4}_2)$.  Note that, both $u_8$ and $u_2$ has face-sequence $(3.6.3.6)$ but $u_8$ has four consecutive vertices having face-sequence $(3^4.6)$ while $u_2$ has two. This can not be possible.

{\bf Subcase 1.2:} If ${\rm lk}(u_7)=C_{8}(\boldsymbol{u_{16}},\boldsymbol{u_{17}},\boldsymbol{u_{18}},u_{8},u,u_6,u_{14},u_{15})$, then proceeding similarly as in Subcase 1.1, we get $N_{19}$ and $N_{20}$ of types $T_{19}$ and $T_{20}$ respectively.

{\bf Case 2:} If ${\rm lk}(u_5)=C_{8}(\boldsymbol{u_{10}},\boldsymbol{u_{11}},\boldsymbol{u_{12}},  u_{13}, u_6,u,u_{4},u_9)$ or ${\rm lk}(u_5)=C_{8}(\boldsymbol{u_9},\boldsymbol{u_{10}},\boldsymbol{u_{11}},u_{12},u_{13},u_6, \linebreak u,u_4)$, then proceeding similarly as in Case 1, we get again $N_{19}$ and $N_{20}$ of types $T_{19}$ and $T_{20}$ respectively. Thus the claim.


\begin{claim}\label{c8}
	For the pairs $((3.6.3.6), (3^2.6^2)) \in B$, we get a unique DSEM $N_{21}$ of type $T_{21}$ on the plane, as shown in example Sec. \ref{s3}. 
\end{claim}

Let $u$ be a vertex with $f_{seq}(u) = (3^2.6^2)$ and ${\rm lk}(u) = C_{10}(\boldsymbol{u_1},\boldsymbol{u_2}, u_3,\boldsymbol{u_4},\boldsymbol{u_5},\boldsymbol{u_6},u_7,u_{8}$, $u_9,\boldsymbol{u_{10}})$. 
Then $f_{seq}(u_3) = f_{seq}(u_8) = (3^2.6^2)$ and $f_{seq}(u_7) = f_{seq}(u_9) = (3.6.3.6)$. 
Let ${\rm lk}(u_3) = C_{10}(\boldsymbol{u_1},u_2$, $u_{11},u_4,\boldsymbol{u_5},\boldsymbol{u_6},\boldsymbol{u_7},u,\boldsymbol{u_9},\boldsymbol{u_{10}})$. Then $f_{seq}(u_{11}) = (3^2.6^2)$ and this implies $f_{seq}(u_2) = f_{seq}(u_4) = (3.6.3.6)$. Now we show: 


\begin{itemize}
	\item The vertices $u_1,u_5,u_6,u_{10}$ can not have face-sequence $(3.6.3.6)$.
\end{itemize}

Suppose, if possible, $f_{seq}(u_5) = (3.6.3.6)$, then $f_{seq}(u_6) = (3.6.3.6)$. If ${\rm lk}(u_7) = C_{10}(u,\boldsymbol{u_3},\boldsymbol{u_4}$, $\boldsymbol{u_5},u_6,u_{12},\boldsymbol{u_{13}},\boldsymbol{u_{14}},\boldsymbol{u_{15}},u_{8})$ and ${\rm lk}(u_6) = C_{10}(\boldsymbol{u},\boldsymbol{u_3},\boldsymbol{u_4}, u_5,u_{16},\boldsymbol{u_{17}},\boldsymbol{u_{18}},\boldsymbol{u_{19}},u_{12},u_7)$, then \linebreak
$f_{seq}(u_{16}) = f_{seq}(u_{19}) = f_{seq}(u_{13}) = f_{seq}(u_{14}) = (3.6.3.6)$ and $f_{seq}(u_{15}) = (3^2.6^2)$. Now, observe that, there are four consecutive vertices in ${\rm lk}(u_7)$ with face-sequence $(3^2.6^2)$ while there are at most two consecutive vertex in ${\rm lk}(u_6)$ with face-sequence $(3^2.6^2)$. This can not be possible. Therefore, $f_{seq}(u_5)\neq (3.6.3.6)$. Similarly, we see that $u_1,u_6,u_{10}$ can not have face-sequence $(3.6.3.6)$. 
Thus $f_{seq}(u_1) = f_{seq}(u_5) = f_{seq}(u_6) = f_{seq}(u_{10}) = (3^2.6^2)$. This implies that the vertices with face-sequences $(3^2.6^2)$ and $(3.6.3.6)$ can have links of face-sequences $({f_2}^2.{f_1}.{f_2}.{f_1}.{f_2}^2.{f_1}.{f_2}.{f_1})$ and $({f_2}^4.{f_1}.{f_2}^4.{f_1})$ respectively. 
By Lemma \ref{l1}, we get $N_{21}$ of type $T_{21}$ = $[(3.6.3.6)^{({f_2}^4.{f_1}.{f_2}^4.{f_1})}:(3^2.6^2)^{({f_2}^2.{f_1}.{f_2}.{f_1}.{f_2}^2.{f_1}.{f_2}.{f_1})}]$. 

\begin{claim}\label{c9}
	For the pairs $((3.6.3.6), (3.4^2.6)) \in B$, we get infinitely many non-isomorphic DSEMs of types $T_{22}$ on the plane.
\end{claim}

Let $u$ be a vertex with face-sequence $(3.4^2.6)$ and ${\rm lk}(u) = C_{9}(\boldsymbol{u_1},u_2,\boldsymbol{u_3},\boldsymbol{u_4},\boldsymbol{u_5},u_6,u_7,\boldsymbol{u_8},u_9)$. Then $f_{seq}(u_1) = f_{seq}(u_2) = f_{seq}(u_7) = f_{seq}(u_8) = f_{seq}(u_{9}) = (3.4^2.6)$. If ${\rm lk}(u_2)  = C_{9}(\boldsymbol{u_{10}},u_{11},u_3$, $\boldsymbol{u_{4}},\boldsymbol{u_{5}},\boldsymbol{u_6},u,\boldsymbol{u_9},u_1)$ and ${\rm lk}(u_7)  = C_{9}(u,u_6,\boldsymbol{u_{12}},\boldsymbol{u_{13}}, \boldsymbol{u_{14}},u_{15}$, $\boldsymbol{u_{16}},u_8,\boldsymbol{u_{9}})$, then $f_{seq}(u_3) = f_{seq}(u_6)$ $= f_{seq}(u_{14}) = (3.6.3.6)$. Now we show: 


\begin{itemize}
	\item The vertices $u_5,u_{12},u_{13},u_4$ can not have face-sequence $(3.6.3.6)$.
\end{itemize}

Suppose, if possible, $f_{seq}(u_5) = (3.6.3.6)$. Then $f_{seq}(u_{12}) = f_{seq} (u_{13}) = (3.6.3.6)$. Note that $f_{seq}(u_6) = f_{seq}(u_{12}) = (3.6.3.6)$ but $u_6$ has four consecutive vertices with face-sequence $(3.4^2.6)$ in its link while $u_{12}$ has at most two. This is not possible. So, 
$f_{seq}(u_5) = (3.4^2.6)$. Similarly, we see that $f_{seq}(u_{12}) = f_{seq}(u_{13}) = f_{seq}(u_4) = (3.4^2.6)$. This gives that the vertices with face-sequence $(3.6.3.6)$ can have links of face-sequence $({f_1}^4.{f_2}.{f_1}^4.{f_2})$ and the vertices with face-sequence $(3.4^2.6)$ can have links of face-sequence $({f_1}^4.{f_2}.{f_1}^2.{f_2})$.


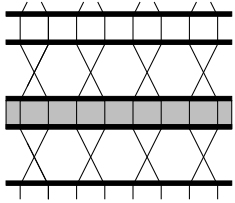
Observe that, there are two possibilities for the link of vertices having face-sequence $(3.4^2.6)$. For example, if we consider  ${\rm lk}(u_9)$, then either ${\rm lk}(u_9)  = C_{9}(u,\boldsymbol{u_7},u_8,u_{17},\boldsymbol{u_{18}},\boldsymbol{u_{19}},\boldsymbol{u_{20}},u_1,\boldsymbol{u_2})$ or ${\rm lk}(u_9)  = C_{9}(u,\boldsymbol{u_7},u_8,\boldsymbol{u_{17}},\boldsymbol{u_{18}},\boldsymbol{u_{19}},u_{20},u_1,\boldsymbol{u_2})$. Using these two possibilities of the links of vertices having face-sequence $(3.4^2.6)$ we can construct infinitely many DSEMs of type $T_{22}$ as follows. Let $H_1$ be an strip in $M$ containing only quadrangular faces such that exactly two triangles are adjacent with each quadrangular face in $H_1$ (depicted by gray color in Figure 22(i)), and let $H_2$ be an strip in $M$ containing quadrangular faces such that exactly one triangular face is adjacent with each quadrangular face in $H_2$ (depicted by brown color in Figure 22(ii)). Now if we construct $M$ which has only $H_1$ type strips, we get $M_{1}$ of type $T_{22}$ as shown in Figure 22(i). If we construct $M$ which has only $H_2$ type strips, we get $M_{2}$ of type $T_{22}$. Further, using both $H_1$ and $H_2$ types strips, we can construct infinitely many non-isomorphic DSEMs, see for example $M_3$ and $M_4$ of type $T_{22}$ in Figures 22(iii) and 22(iv) respectively. 

This proves the claim and hence Theorem \ref{t1}. \hfill$\Box$




\section{Proof of Theorem \ref{t2} and Corollary \ref{co2}}\label{s4}

The enumeration presented here for DSEMs of types $T_i$, for $1 \leq i \leq 4$, is motivated by the construction given in \cite{alt(1973)} and further extended in \cite{MaUp(2018)}. We give this enumeration in a detailed way for the type $T_1, T_2, T_3$ and $T_4$. However, similarly, one can enumerate DSEMs of the remaining types $T_i$ ($5 \leq i \leq 22$). Throughout this section, by a DSEM, we mean a doubly semi-equivelar map on the torus. The  labeling of vertices in the link of a vertex is considered here anti-clock wise.


\subsection{DSEMs of type $T_1 = [f_1^{({f_1}^2.{f_2}. {f_1}^2.{f_2})}:f_2^{({f_1}^5.{f_2}^2)}]$, where $(f_1, f_2)=((3^6), (3^3.4^2))$} \label{s4.1}

Let $M_1$ be a DSEM of type $T_1$ with the vertex set $V(M_1)$. Let $V_{(3^6)}$ and $V_{(3^3.4^2)}$ denote the vertex sets with face-sequences $(3^6)$ and $(3^3.4^2)$ respectively. It is easy to see that the number of triangular faces in $M_1$ is $4|V(3^6)|$ or $2|V(3^3.4^2)|$, where $|V(3^6)|$ and $|V(3^3.4^2)|$ denote the cardinality of sets $V(3^6)$ and $V(3^3.4^2)$ respectively. Thus for the existence of $M_1$, we have $2|V(3^6)| = |V(3^3.4^2)|$. Now we define following three types of paths in $M_1$ (Figure 5.1.1) as follows.


\begin{deff} \normalfont
	
	A path $P_{1} = P( \ldots, y_{i-1},y_{i},y_{i+1}, \ldots)$ in $M_1$ is of type $A_{1}$ if either of the conditions follows: 
	$(i)$ every vertex of $P_1$ has the face-sequence $(3^6)$ or
	$(ii)$ each vertex of $P_1$ has the face-sequence $(3^3.4^2)$ such that all the triangles (or squares) that incident on its inner vertices lie on the one side  of the path $P_1$. For example, see thick balck colored paths in Figure 5.1.1.
	
\end{deff}

\begin{deff}\label{d4.1.2} \normalfont A path $P_{2} = P(\ldots, z_{i-1},z_{i},z_{i+1}, \ldots)$ in $M_1$, such that $z_{i-1}, z_{i}, z_{i+1}$ are inner vertices of $P_{2}$ or an extended path of $P_{2}$, is of type $A_{2}$ (indicated by red colored  paths in Figure 5.1.1) if either of the following three conditions follows for each vertex of the path:

	\begin{enumerate}

		\item  If ${\rm lk}(z_{i})=C_7(\boldsymbol{m},z_{i-1},\boldsymbol{n},o,z_{i+1},p,q)$ and  ${\rm lk}(z_{i-1})= C_7(\boldsymbol{q},m,z_{i-2},r,n,\boldsymbol{o},z_{i})$,  then ${\rm lk}(z_{i+1}) \linebreak =C_6(s,z_{i+2},t,p,z_i,o)$.
		
		\item If ${\rm lk}(z_{i})=C_7(\boldsymbol{m},z_{i+1},\boldsymbol{n},o,z_{i-1},p,q)$ and  ${\rm lk}(z_{i-1})=C_6(r,z_{i-2},s,p,z_{i},o)$, then ${\rm lk}(z_{i+1})=C_7(\boldsymbol{o},z_{i},\boldsymbol{q},m,z_{i+2},t,n)$. 
		
		\item  If ${\rm lk}(z_{i})=C_6(z_{i-1},m,n,z_{i+1},o,p)$ and ${\rm lk}(z_{i-1})=C_7(\boldsymbol{r},z_{i-2},\boldsymbol{s},m,z_{i},p,q)$, then ${\rm lk}(z_{i+1})=C_7(\boldsymbol{t},z_{i+2},\boldsymbol{u},o,z_{i},n,v)$.
		
	\end{enumerate}
	
\end{deff}

\begin{deff} \label{d4.1.3} \normalfont A path $P_{3} = P(\ldots, w_{i-1},w_{i},w_{i+1}, \ldots)$ in $M_1$, such that $w_{i-1}, w_i, w_{i+1}$ are inner vertices of $P_3$ or an extended path of $P_3$, is of type $A_{3}$ (shown by green colored paths in Figure 5.1.1), if either of the following three conditions follows for each vertex of the path.

	\begin{enumerate}

		\item  If ${\rm lk}(w_{i})=C_7(\boldsymbol{m},w_{i-1},\boldsymbol{n},o,p,w_{i+1},q)$ and  ${\rm lk}(w_{i-1})=C_7(\boldsymbol{q},m,r,w_{i-2},n,\boldsymbol{o},w_{i})$,  then ${\rm lk}(w_{i+1}) \linebreak =C_6(s,q,w_i,p,t,w_{i+2})$.
		
		\item If ${\rm lk}(w_{i})=C_7(\boldsymbol{m},w_{i+1},\boldsymbol{n},o,p,w_{i-1},q)$ and  ${\rm lk}(w_{i-1})=C_6(r,w_{i-2},s,q,w_{i},p)$, then ${\rm lk}(w_{i+1})=C_7(\boldsymbol{o},w_{i},\boldsymbol{q},m,t,w_{i+2},n)$.

		\item  If ${\rm lk}(w_{i})=C_6(w_{i-1},m,n,w_{i+1},o,p)$ and ${\rm lk}(w_{i-1})=C_7(\boldsymbol{r},w_{i-2},\boldsymbol{s},t,m,w_{i},p)$, then ${\rm lk}(w_{i+1})=C_7(\boldsymbol{u},w_{i+2},\boldsymbol{v},w,o,w_{i},n)$.
		
	\end{enumerate}
	
\end{deff}

\vspace{-.5cm}
\begin{picture}(0,0)(-26,23)
\setlength{\unitlength}{5mm}

\drawpolygon(1,0)(.5,1)(-.5,1)(-1,0)(-.5,-1)(.5,-1)
\drawpolygon(3,0)(2.5,1)(1.5,1)(1,0)(1.5,-1)(2.5,-1)
\drawpolygon(5,0)(4.5,1)(3.5,1)(3,0)(3.5,-1)(4.5,-1)

\drawline[AHnb=0](4.5,1)(4.5,2)
\drawline[AHnb=0](4.5,-1)(4.5,-2)

\drawline[AHnb=0](-.5,1)(-.5,2)
\drawline[AHnb=0](-.5,-1)(-.5,-2)

\drawline[AHnb=0](.5,-1)(.5,-2)
\drawline[AHnb=0](1.5,-1)(1.5,-2)
\drawline[AHnb=0](2.5,-1)(2.5,-2)
\drawline[AHnb=0](3.5,-1)(3.5,-2)

\drawline[AHnb=0](.5,1)(.5,2)
\drawline[AHnb=0](1.5,1)(1.5,2)

\drawline[AHnb=0](2.5,1)(2.5,2)
\drawline[AHnb=0](3.5,1)(3.5,2)


\drawline[AHnb=0](-1.11,0)(5.11,0)
\drawline[AHnb=0](-1,-2)(5,-2)
\drawline[AHnb=0](-1,2)(5,2)
\drawline[AHnb=0](-1,1)(5,1)


\drawline[AHnb=0](.5,1)(-.5,-1)
\drawline[AHnb=0](-.5,1)(.5,-1)

\drawline[AHnb=0](1.5,-1)(2.5,1)
\drawline[AHnb=0](2.5,-1)(1.5,1)

\drawline[AHnb=0](3.5,-1)(4.5,1)
\drawline[AHnb=0](4.5,-1)(3.5,1)


\drawline[AHnb=0](.5,2)(.75,2.5)
\drawline[AHnb=0](.5,2)(.25,2.5)

\drawline[AHnb=0](-.5,2)(-.75,2.5)
\drawline[AHnb=0](-.5,2)(-.25,2.5)

\drawline[AHnb=0](1.5,2)(1.75,2.5)
\drawline[AHnb=0](1.5,2)(1.25,2.5)

\drawline[AHnb=0](3.5,2)(3.75,2.5)
\drawline[AHnb=0](3.5,2)(3.25,2.5)

\drawline[AHnb=0](4.5,2)(4.75,2.5)
\drawline[AHnb=0](4.5,2)(4.25,2.5)

\drawline[AHnb=0](2.5,2)(2.25,2.5)
\drawline[AHnb=0](2.5,2)(2.75,2.5)

\drawline[AHnb=0](.5,-2)(.75,-2.5)
\drawline[AHnb=0](.5,-2)(.25,-2.5)

\drawline[AHnb=0](-.5,-2)(-.75,-2.5)
\drawline[AHnb=0](-.5,-2)(-.25,-2.5)

\drawline[AHnb=0](1.5,-2)(1.75,-2.5)
\drawline[AHnb=0](1.5,-2)(1.25,-2.5)

\drawline[AHnb=0](3.5,-2)(3.75,-2.5)
\drawline[AHnb=0](3.5,-2)(3.25,-2.5)

\drawline[AHnb=0](4.5,-2)(4.75,-2.5)
\drawline[AHnb=0](4.5,-2)(4.25,-2.5)

\drawline[AHnb=0](2.5,-2)(2.25,-2.5)
\drawline[AHnb=0](2.5,-2)(2.75,-2.5)


\drawpolygon[fillcolor=black](-1,-2)(5,-2)(5,-2.1)(-1,-2.1)

\drawpolygon[fillcolor=black](-1,-1)(5,-1)(5,-1.1)(-1,-1.1)

\drawpolygon[fillcolor=black](-1.12,0)(5.1,0)(5.1,.1)(-1.12,.1)

\drawpolygon[fillcolor=black](-1,1)(5,1)(5,1.1)(-1,1.1)

\drawpolygon[fillcolor=black](-1,2)(5,2)(5,2.1)(-1,2.1)

\drawpolygon[fillcolor=green](.75,-2.5)(.5,-2)(.5,-1)(-.5,1)(-.5,2)(-.75,2.5)(-.65,2.5)(-.4,2)(-.4,1)(.6,-1)(.6,-2)(.85,-2.5)

\drawpolygon[fillcolor=green](1.75,-2.5)(1.5,-2)(1.5,-1)(.5,1)(.5,2)(.25,2.5)(.35,2.5)(.6,2)(.6,1)(1.6,-1)(1.6,-2)(1.85,-2.5)

\drawpolygon[fillcolor=green](2.75,-2.5)(2.5,-2)(2.5,-1)(1.5,1)(1.5,2)(1.25,2.5)(1.35,2.5)(1.6,2)(1.6,1)(2.6,-1)(2.6,-2)(2.85,-2.5)

\drawpolygon[fillcolor=green](3.75,-2.5)(3.5,-2)(3.5,-1)(2.5,1)(2.5,2)(2.25,2.5)(2.35,2.5)(2.6,2)(2.6,1)(3.6,-1)(3.6,-2)(3.85,-2.5)

\drawpolygon[fillcolor=green](4.75,-2.5)(4.5,-2)(4.5,-1)(3.5,1)(3.5,2)(3.25,2.5)(3.35,2.5)(3.6,2)(3.6,1)(4.6,-1)(4.6,-2)(4.85,-2.5)

\drawpolygon[fillcolor=red](-.75,-2.5)(-.5,-2)(-.5,-1)(.5,1)(.5,2)(.75,2.5)(.65,2.5)(.4,2)(.4,1)(-.6,-1)(-.6,-2)(-.85,-2.5)

\drawpolygon[fillcolor=red](.25,-2.5)(.5,-2)(.5,-1)(1.5,1)(1.5,2)(1.75,2.5)(1.65,2.5)(1.4,2)(1.4,1)(.4,-1)(.4,-2)(.15,-2.5)

\drawpolygon[fillcolor=red](1.25,-2.5)(1.5,-2)(1.5,-1)(2.5,1)(2.5,2)(2.75,2.5)(2.65,2.5)(2.4,2)(2.4,1)(1.4,-1)(1.4,-2)(1.15,-2.5)

\drawpolygon[fillcolor=red](2.25,-2.5)(2.5,-2)(2.5,-1)(3.5,1)(3.5,2)(3.75,2.5)(3.65,2.5)(3.4,2)(3.4,1)(2.4,-1)(2.4,-2)(2.15,-2.5)

\drawpolygon[fillcolor=red](3.25,-2.5)(3.5,-2)(3.5,-1)(4.5,1)(4.5,2)(4.75,2.5)(4.65,2.5)(4.4,2)(4.4,1)(3.4,-1)(3.4,-2)(3.15,-2.5)


\put(1,2.65) {\scriptsize {\tiny $A_3$}}
\put(2.5,2.65) {\scriptsize {\tiny $A_2$}}
\put(5.5,-1.1) {\scriptsize {\tiny $A_1$}}
\put(5.5,-.1) {\scriptsize {\tiny $A_1$}}

\put(-2.5,-3.5){\scriptsize {\tiny {\bf Figure 5.1.1:} Paths of types $A_1, A_2, A_3$}}
\end{picture}

\begin{picture}(0,0)(-80,30)
\setlength{\unitlength}{4.25mm}

\drawpolygon(0,0)(3,0)(3,3)(0,3)
\drawpolygon(0,5)(3,5)(3,8)(0,8)

\drawline[AHnb=0](1,0)(1,3)
\drawline[AHnb=0](2,0)(2,3)

\drawline[AHnb=0](1,5)(1,8)
\drawline[AHnb=0](2,5)(2,8)


\drawline[AHnb=0](0,1)(3,1)
\drawline[AHnb=0](0,2)(3,2)

\drawline[AHnb=0](0,6)(3,6)
\drawline[AHnb=0](0,7)(3,7)


\drawline[AHnb=0](0,1)(2,3)
\drawline[AHnb=0](1,1)(3,3)
\drawline[AHnb=0](2,1)(3,2)
\drawline[AHnb=0](0,2)(1,3)

\drawline[AHnb=0](0,6)(2,8)
\drawline[AHnb=0](1,6)(3,8)
\drawline[AHnb=0](2,6)(3,7)
\drawline[AHnb=0](0,7)(1,8)

\put(3.5,.6){\scriptsize $\ldots$} 
\put(3.5,.6){\scriptsize $\ldots$} 

\put(3.5,2.1){\scriptsize $\ldots$} 
\put(3.5,2.1){\scriptsize $\ldots$} 

\put(3.5,5.6){\scriptsize $\ldots$} 
\put(3.5,5.6){\scriptsize $\ldots$} 

\put(3.5,7.1){\scriptsize $\ldots$} 
\put(3.5,7.1){\scriptsize $\ldots$}
\put(8.1,1){\scriptsize $\ldots$} 
\put(8.1,1){\scriptsize $\ldots$}

\put(8.1,2.1){\scriptsize $\ldots$} 
\put(8.1,2.1){\scriptsize $\ldots$} 

\put(8.1,7.1){\scriptsize $\ldots$} 
\put(8.1,7.1){\scriptsize $\ldots$} 

\put(8.1,5.6){\scriptsize $\ldots$} 
\put(8.1,5.6){\scriptsize $\ldots$} 

\put(1.5,3.7){\scriptsize $\vdots$}
\put(6,3.7){\scriptsize $\vdots$}
\put(10.5,3.7){\scriptsize $\vdots$}

\put(1.5,3.7){\scriptsize $\vdots$}
\put(6,3.7){\scriptsize $\vdots$}
\put(10.5,3.7){\scriptsize $\vdots$}

\drawpolygon(4.5,0)(7.5,0)(7.5,3)(4.5,3)

\drawpolygon(4.5,5)(7.5,5)(7.5,8)(4.5,8)


\drawline[AHnb=0](5.5,0)(5.5,3)
\drawline[AHnb=0](6.5,0)(6.5,3)

\drawline[AHnb=0](5.5,5)(5.5,8)
\drawline[AHnb=0](6.5,5)(6.5,8)

\drawline[AHnb=0](4.5,1)(7.5,1)
\drawline[AHnb=0](4.5,2)(7.5,2)

\drawline[AHnb=0](4.5,6)(7.5,6)
\drawline[AHnb=0](4.5,7)(7.5,7)

\drawline[AHnb=0](4.5,1)(6.5,3)
\drawline[AHnb=0](5.5,1)(7.5,3)
\drawline[AHnb=0](6.5,1)(7.5,2)
\drawline[AHnb=0](4.5,2)(5.5,3)

\drawline[AHnb=0](4.5,6)(6.5,8)
\drawline[AHnb=0](5.5,6)(7.5,8)
\drawline[AHnb=0](6.5,6)(7.5,7)
\drawline[AHnb=0](4.5,7)(5.5,8)

\drawpolygon(9,0)(12,0)(12,3)(9,3)
\drawpolygon(9,5)(12,5)(12,8)(9,8)

\drawline[AHnb=0](10,0)(10,3)
\drawline[AHnb=0](11,0)(11,3)

\drawline[AHnb=0](10,5)(10,8)
\drawline[AHnb=0](11,5)(11,8)


\drawline[AHnb=0](9,1)(12,1)
\drawline[AHnb=0](9,2)(12,2)

\drawline[AHnb=0](9,6)(12,6)
\drawline[AHnb=0](9,7)(12,7)

\drawline[AHnb=0](9,1)(11,3)
\drawline[AHnb=0](10,1)(12,3)
\drawline[AHnb=0](11,1)(12,2)
\drawline[AHnb=0](9,2)(10,3)

\drawline[AHnb=0](9,6)(11,8)
\drawline[AHnb=0](10,6)(12,8)
\drawline[AHnb=0](11,6)(12,7)
\drawline[AHnb=0](9,7)(10,8)

\put(-.2,-.5){\scriptsize {\tiny $1$}}
\put(.8,-.5){\scriptsize {\tiny $2$}}
\put(1.8,-.5){\scriptsize {\tiny $3$}}
\put(2.8,-.5){\scriptsize {\tiny $4$}}
\put(4.4,-.5){\scriptsize {\tiny $k$}}
\put(5,-.5){\scriptsize {\tiny $k$+$1$}}
\put(6.1,-.5){\scriptsize {\tiny $k$+$2$}}
\put(7.2,-.5){\scriptsize {\tiny $k$+3}}
\put(8.8,-.5){\scriptsize {\tiny $i$-2}}
\put(9.9,-.5){\scriptsize {\tiny $i$-1}}
\put(11,-.5){\scriptsize {\tiny $i$}}
\put(12,-.5){\scriptsize {\tiny $1$}}

\put(-.75,8.2){\scriptsize {\tiny $k$+1}}
\put(.4,8.2){\scriptsize {\tiny $k$+2}}
\put(1.6,8.2){\scriptsize {\tiny $k$+3}}
\put(2.7,8.2){\scriptsize {\tiny $k$+4}}

\put(8.7,8.2){\scriptsize {\tiny $k$-2}}
\put(9.7,8.2){\scriptsize {\tiny $k$-1}}
\put(10.8,8.2){\scriptsize {\tiny $k$}}
\put(11.6,8.2){\scriptsize {\tiny $k+1$}}

\put(-1.6,-.2){\scriptsize {\tiny $Q_1$}}
\put(-1.6,.8){\scriptsize {\tiny $Q_2$}}
\put(-1.6,1.8){\scriptsize {\tiny $Q_3$}}
\put(-1.6,2.8){\scriptsize {\tiny $Q_4$}}
\put(-1.6,6.8){\scriptsize {\tiny $Q_{j}$}}
\put(-1.6,7.8){\scriptsize {\tiny $Q_{1}$}}
\put(-1.6,5.8){\scriptsize {\tiny $Q_{j-1}$}}
\put(-1.6,4.8){\scriptsize {\tiny $Q_{j-2}$}}

\put(3.6,-1.4){\scriptsize {\tiny {\bf Figure 5.1.2:} $M(i,j,k)$ }} 

\end{picture}

\newpage

\begin{picture}(0,0)(-20,18.25)
\setlength{\unitlength}{6mm}

\drawpolygon(0,0)(3,0)(3,2)(0,2)


\drawline[AHnb=0](1,0)(1,2)
\drawline[AHnb=0](2,0)(2,2)
\drawline[AHnb=0](3,0)(3,2)


\drawline[AHnb=0](0,1)(3,1)
\drawline[AHnb=0](0,2)(3,2)


\drawline[AHnb=0](0,1)(1,2)
\drawline[AHnb=0](1,1)(2,2)
\drawline[AHnb=0](2,1)(3,2)

\put(-.2,-.3){\scriptsize {\tiny $x_1$}}
\put(.8,-.3){\scriptsize {\tiny $x_2$}}
\put(1.8,-.3){\scriptsize {\tiny $x_3$}}
\put(2.8,-.3){\scriptsize {\tiny $x_4$}}

\put(3.7,1){\scriptsize $\ldots$} 
\put(3.7,1){\scriptsize $\ldots$} 
\put(3.7,1){\scriptsize $\ldots$}

\put(3.7,1.8){\scriptsize $\ldots$} 
\put(3.7,1.8){\scriptsize $\ldots$} 
\put(3.7,1.8){\scriptsize $\ldots$} 

\put(3.7,0){\scriptsize $\ldots$} 
\put(3.7,0){\scriptsize $\ldots$} 
\put(3.7,0){\scriptsize $\ldots$}

\put(-.6,.7){\scriptsize {\tiny $y_1$}}
\put(1.05,.7){\scriptsize {\tiny $y_2$}}
\put(2.05,.7){\scriptsize {\tiny $y_3$}}
\put(3.05,.7){\scriptsize {\tiny $y_4$}}

\put(-.6,1.7){\scriptsize {\tiny $z_1$}}
\put(1.05,1.7){\scriptsize {\tiny $z_2$}}
\put(2.05,1.7){\scriptsize {\tiny $z_3$}}
\put(3.05,1.7){\scriptsize {\tiny $z_4$}}

\drawpolygon(4.5,0)(5.5,0)(5.5,2)(4.5,2)


\drawline[AHnb=0](4.5,1)(5.5,1)


\drawline[AHnb=0](4.5,1)(5.5,2)

\put(4.5,-.3){\scriptsize {\tiny $x_i$}}  \put(5.5,-.3){\scriptsize {\tiny $x_1$}}

\put(4.6,.7){\scriptsize {\tiny $y_i$}} 
\put(5.6,.7){\scriptsize {\tiny $y_1$}}

\put(4.55,1.7){\scriptsize {\tiny $z_i$}}  \put(5.6,1.7){\scriptsize {\tiny $z_1$}}

\put(.6,-1){\scriptsize {\tiny {\bf Figure 5.1.3:} Cylinder-I}}

\end{picture}

\begin{picture}(0,0)(-88,15) 
\setlength{\unitlength}{6.5mm}

\drawpolygon(0,0)(3,0)(3,2)(0,2)


\drawline[AHnb=0](1,0)(1,2)
\drawline[AHnb=0](2,0)(2,2)


\drawline[AHnb=0](0,1)(3,1)


\drawline[AHnb=0](0,1)(1,2)
\drawline[AHnb=0](1,0)(3,2)
\drawline[AHnb=0](0,0)(2,2)
\drawline[AHnb=0](2,0)(3,1)

\put(-.2,-.3){\scriptsize {\tiny $u_1$}}
\put(.8,-.3){\scriptsize {\tiny $u_2$}}
\put(1.8,-.3){\scriptsize {\tiny $u_3$}}
\put(2.9,-.3){\scriptsize {\tiny $u_4$}}

\put(-.5,.7){\scriptsize {\tiny $v_1$}}
\put(1.05,.7){\scriptsize {\tiny $v_2$}}
\put(2.05,.7){\scriptsize {\tiny $v_3$}}
\put(3.05,.7){\scriptsize {\tiny $v_4$}}

\put(-.6,1.7){\scriptsize {\tiny $w_1$}}
\put(1.05,1.7){\scriptsize {\tiny $w_2$}}
\put(2.05,1.7){\scriptsize {\tiny $w_3$}}
\put(3.05,1.7){\scriptsize {\tiny $w_4$}}

\put(3.5,1){\scriptsize $\ldots$} 
\put(3.5,1){\scriptsize $\ldots$}
\put(3.5,1){\scriptsize $\ldots$}

\put(3.5,1.8){\scriptsize $\ldots$} 
\put(3.5,1.8){\scriptsize $\ldots$} 
\put(3.5,1.8){\scriptsize $\ldots$} 

\put(3.5,0){\scriptsize $\ldots$} 
\put(3.5,0){\scriptsize $\ldots$} 
\put(3.5,0){\scriptsize $\ldots$}

\drawpolygon(4.2,0)(5.2,0)(5.2,2)(4.2,2)


\drawline[AHnb=0](4.2,1)(5.2,2)
\drawline[AHnb=0](4.2,0)(5.2,1)

\drawline[AHnb=0](4.2,1)(5.2,1)

\put(4.2,-.3){\scriptsize {\tiny $u_i$}}  \put(5.3,-.3){\scriptsize {\tiny $u_1$}}

\put(4.3,.7){\scriptsize {\tiny $v_i$}} \put(5.3,.7){\scriptsize {\tiny $v_1$}}

\put(4.3,1.7){\scriptsize {\tiny $w_i$}}  \put(5.3,1.7){\scriptsize {\tiny $w_1$}}

\put(.6,-1){\scriptsize {\tiny {\bf Figure 5.1.4:} Cylinder-II}}

\end{picture}

\vspace{2cm}

\begin{rem}\label{r1}
	Note that $M_1$ is a map on finite vertex set. Thus, if $P(v_1, v_2, \ldots, v_l)$ is a maximal path (a path of maximum length) of type $A_{\alpha}$ for $\alpha \in \{1,2,3\}$ in $M_1$, then the path gives a cycle $Q = C_l(v_1, v_2, \ldots, v_l)$. We refer to the corresponding cycle as a cycle of type $A_{\alpha}$. 
	
\end{rem}

For a cycle $Q$ of type $A_{\alpha}$ with the vertex set $V(Q)$, where ${\alpha} \in \{1,2,3\}$, we have following.





\begin{lem}\label{l4.1.3} If $Q$ is a cycle of type $A_{\alpha}$, for $\alpha \in$ $\{1, 2, 3\} $, in $M_1$ then $Q$ is non-contractible.
\end{lem}

\noindent{\bf Proof.}  Suppose $Q$ is a contractible cycle of type $A_{1}$. Then $Q$ is the boundary of a 2-disk $\mathbb{D}^2$. It is disused in Lemma \ref{l1} that this is not possible. Thus $Q$ is non-contractible. Similarly, we see, if $Q$ is a cycle of type $A_2$ or $A_{3}$ then it is also non-contractible. Thus the proof. \hfill $\Box$

\vspace{.2cm}
Let $S$ denote the set of all the faces incident at $v$ for all $v \in V(Q)$. Then the geometric carrier $S_{Q}$ (union of all the faces in $S$) is a cylinder, as $Q$ is non-contractible (for example see Figures 5.1.3 and 5.1.4). 
Let $\partial S_{Q} = \{Q_1, Q_2\}$, where $\partial S_Q$ denote the boundary of $S_{Q}$ and $Q_{1}$, $Q_{2}$ be the boundary cycles. Then we show:

\begin{lem}\label{l4.1.4} Let $Q$ be a cycle of type $A_{\alpha}$, $\alpha \in \{1, 2, 3\}$, such that $\partial S_{Q} = \{Q_{1}, Q_{2}\}$. Then $(i)$ $Q$, $Q_1$ and $Q_2$ are of same type, $(ii) $ length$(Q)$ = length$(Q_{1})$ = length$(Q_{2})$. 
\end{lem}

\noindent{\bf Proof.}  Let $Q$ be a cycle of type $A_{1}$ such that $S_{Q}$ is a cylinder with the boundary $\partial S_{Q} = \{Q_{1},Q_{2}\}$. Consider the faces which are incident with both $Q$ and $Q_n$ for a fixed $n \in \{1,2\}$. Without loss of generality let $Q_n = Q_1$. Now, depending on $Q$, if the faces incident with both $Q$ and $Q_1$ are quadrangles then the faces lying on the other side of $Q_1$ must be triangles and $(ii)$ if the faces incident with both $Q$ and $Q_1$ are triangles then the faces on the other side of $Q_1$ must be either triangles or quadrangles. From the definition of type $A_1$, we see for both the cases, $Q_1$ is of type $A_1$ and hence $Q_2$ is also of type $A_1$.        

Now let $Q = C_l(v_1, \ldots, v_l)$, $Q_{1} = C_{l_1}(u_1, \ldots, u_{l_1})$ and $Q_2 = C_{l_2}(w_1, \ldots, w_{l_2})$. We show that $l = l_1 = l_2$. Suppose $l \neq l_1 \neq l_2$. Without loss of generality let $ l < l_1$. By the definition of type $A_1$, the face-sequence of $v_{1},v_{2}, \ldots, v_{l-1},v_{l}$ will be same through out the cycle. Since $l < l_1$, so, the ${\rm lk}(v_{l})$ contains the vertices $u_{l},u_{l+i}$, and $w_{l}$ for some $i > 0$. This shows that the face-sequences of $v_{l}$ and $v_{l-1}$ are not same, a contradiction. Therefore $l = l_1 = l_2$.  
Proceeding similarly, we get the above result for a cycles of type $A_2$ and $A_3$.  \hfill$\Box$

\vspace{.2cm}

Let $Q_1$ and $Q_2$ be two same type cycles in a DSEM $M_1$ on the torus. We say that cycles $Q_{1}$ and $Q_{2}$ are homologous if there is a cylinder whose boundary is $\{Q_1, Q_2\}$. Thus in Lemma \ref{l4.1.4}, the cycles $Q$, $Q_1$ and $Q_2$ are homologous to each other. 

\smallskip

Now we give the notion of a planar representation, denoted as $M(i,j,k)$ representation, for a DSEM $M_1$. This is obtained by cutting $M_1$ along any two non-homologous cycles.

\smallskip

 Let $u \in V(M_1)$ and $Q_{\alpha}$ be cycles of type $A_\alpha$ through $u$, where $ \alpha \in \{1,2,3\}$. Let $Q_{1} = C_i(u_1,u_2, \ldots, u_i)$. We cut $M_1$ first along the cycle $Q_{1}$. We get a cylinder, say $R_1$, bounded by identical cycle $Q_{1}$. We say that a cycle is  horizontal if it is $Q_{1}$ or homologous to $Q_{1}$. Then we say that a cycle is vertical if it is $Q_{\alpha}$ or homologous to $Q_{\alpha}$ for $\alpha\in\{2,3\}$. In $R_{1}$, starting from the vertex $u$, make another cut along the path $P \subset Q_{3}$, until it reaches $Q_{1}$ again for the first time. As a result, we unfold the torus into a planer representation, say $R_{2}$.

\smallskip
\noindent{Claim.} The representation $R_{2}$ is connected.

Since $Q_1$ is a non-contractible cycle, $R_{1}$ is connected. Suppose that $R_{2}$ is disconnected. This means there is a 2-disk $D_{Q}$ with boundary cycle $Q = P_1 \cup P_1' = P(a_{m} , \ldots, a_{n})$ $\cup$ $P(b_{s} , \ldots, b_{t})$ where $P_1 \subset Q_{3}$, $P_1' \subset Q_{1}$, $a_{n} = b_{s}$ and $a_{m} = b_{t}$. Let $\bigtriangleup$ and $\Box$ denote triangular and quadrangular faces respectively. Now in the first case when quadrangular faces are incident with $P_1'$ in $D_{Q}$ and $\Box_{m}$, $\bigtriangleup_{m,1}$, $\bigtriangleup_{m,2}$, $\bigtriangleup_{m,3}$, $\bigtriangleup_{m,4}$, $\Box_{m+1}$, $\bigtriangleup_{m+1,1}$, $\bigtriangleup_{m+1,2}$, $\bigtriangleup_{m+1,3}$, $\bigtriangleup_{m+1,4}, \ldots, \bigtriangleup_{n-1,1}$, $\bigtriangleup_{n-1,2}$,   $\bigtriangleup_{n-1,3}$, $\bigtriangleup_{n-1,4}$, $\Box_{n}$ are incident with $P_1$ in $D_{Q}$ then as in Lemma \ref{l1}, we calculate $v$, $e$ and $f$ of $D_{Q}$ and we get  $v-e+f=0$. On the other hand, if triangular faces are incident with $P_1'$ in $D_{Q}$ then again calculating $v$, $e$ and $f$ of $D_{Q}$, we get $v-e+f=0$. This is not possible as the Euler characteristic of the 2-disk $D_{Q}$ is 1. So, $R_{2}$ is connected.

Without loss of generality, suppose that the quadrangular faces are incident on $Q_{1}$. In $R_{2}$, let there are $j$ cycles which are homologous to $Q_{1}$ along the path $P \subset Q_{3}$. Since length($Q_{1}$) = $i$ and the number of horizontal cycles along $P$ is $j$, as shown in Figure 5.1.2. So, we denote $R_{2}$ by $(i, j)$ representation. Now to obtain $M_1$ from its $(i, j)$ representation one has to go in a reverse way, i.e., identify both the vertical sides and horizontal sides along the vertices, but observe that identification of the horizontal side may need some shifting, as shown in Figure 5.1.2. Let $k$ be the shifting, i.e., $u_{k+1}$ be the first vertex in the upper horizontal cycle. This gives another planar representation of $M_1$ (corresponding to the $(i,j)$ representation) also called $M(i,j,k)$ representation of $M_1$.

\begin{lem}\label{l4.1.6} In $M(i,j,k)$, $A_{1}$ type cycles have unique length and $A_{2}$ type (or $A_3$ type) cycles have at most two different lengths.
\end{lem}

\noindent{\bf Proof.}  Consider an $M(i,j,k)$ representation of DSEM $M_1$. Let $Q_{1}$ be an $A_1$ type cycle in $A_1$. Consider a cylinder $S_{Q_{1}}$ which is defined by $Q_{1}$. Let $\partial S_{Q_{1}}= \{Q_{0}, Q_{2}\}$. By Lemma \ref{l4.1.4}, the cycles $Q_{0}$, $Q_{1}$ and $Q_{2}$ are homologous and length$(Q_{0})$ = length$(Q_{1})$ = length$(Q_{2})$. Now consider the cycle $Q_{2}$ and repeat the above procedure. In this process, let $Q_{m}$ indicate a cycle at $m^{th}$ step such that $\partial S_{Q_{m}}= \{ Q_{m-1},Q_{m+1} \}$ and length$(Q_{m-1})$ = length$(Q_{m})$ = length$(Q_{m+1})$. Since $V(M_1)$ is finite, this process terminates after finite number of steps. By the construction of $M(i,j,k)$, shown in Figure 5.1.2, the process stops
after $j+1$ number of steps, i.e., $Q_{0}$ appears again. Thus, the cycles  $Q_{r}, Q_{s}$  are homologous for every $1 \leq r,s \leq j$ and $\bigcup\limits_{m=1}^{j} V(Q_m)=V(M_1)$. Note that in $M(i,j,k)$ there is only one cycle of type $A_{1}$  through any vertex. As a result, $Q_{1}, Q_{2}, \ldots, Q_{j}$ are the only type $A_{1}$ cycles in $M_1$. Since these cycles are homologous and length$(Q_{1})$ = length$(Q_{r})$ for all $1 \leq r \leq j$. Therefore, $M_1$ has $A_{1}$ type cycles with unique length.

Let $Q_{1}'$, $Q_{1}''$ be the cycles of type $A_2$ and $A_3$ respectively. Now repeating the above process for the cycles $Q_{1}'$ and $Q_{1}''$, we see that all the $A_2$ type cycles have same length say $l_1$ and all the $A_3$ type cycles have same length say $l_{2}$. Observe that $Q_{1}'$ and $Q_{1}''$ define same type cycles as these are mirror image of each other. So, the map $M_1$ contains the cycles of type $A_{2}$ (or type $A_3$) of lengths $l_{1}$ and $l_{2}$. Therefore, $M_1$ has $A_{2}$ type (or $A_3$ type) cycles with at most two different lengths. \hfill$\Box$ 

\smallskip
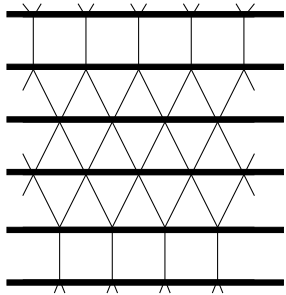
Now we define a cycle of new type (other than $A_1$, $A_2$ and $A_3$) as follows: \\ In $M(i,j,k)$-representation of DSEM $M_1$, Let $Q_{lh} = C_i(x_{1},x_{2}, \ldots, x_{i})$ and $Q_{uh} = C_i(x_{k+1}, \linebreak x_{k+2}$, $\ldots, x_{k})$ be the lower and upper horizontal cycles of type $A_1$ respectively. We define two paths $P_{1} = P(x_{k+1}, y_1,y_2, \ldots, y_{\alpha}, x_{k_{1}})$  and $P_{1}'=P(x_{k+1}, z_1,z_2, \ldots, z_{\beta},x_{k_{2}})$ of type $A_{2}$ and $A_{3}$ respectively through $x_{k+1}$, where $x_{k_{1}}, x_{k_{2}} \in V(Q_{uh})$. Note that, the paths $P_{1}$ and $P_{1}'$ are not parts of horizontal cycles. Consider the paths $P_{2 }= P(x_{k_{1}}, \ldots, x_{k+1})$ and $P_{2}'=P(x_{k_{2}}, \ldots, x_{k+1})$ in $Q_{uh}$ such that $P_{2} \cup P_{2}' \subseteq Q_{uh}$,  shown by the blue colored path in Figure 5.1.5. Let $Q^1_{4} = P_{1} \cup P_{2} = C_{\gamma}(x_{k+1}, y_1,y_2, \ldots, y_{\alpha}, x_{k_1}, \ldots, x_{k+1})$ and $Q^2_{4} = P_{1}'\cup P_{2}' = C_{\gamma^{'}}(x_{k+1}, z_1, z_2, \ldots, z_{\beta},x_{k_2}, \ldots, x_{k+1})$, where $\gamma$ and $\gamma^{'}$ are lengths of $Q^1_{4}$ and $Q^2_{4}$ respectively. Define a cycle $Q_{4}$ of new type, say $A_4$, as: \\ \begin{equation}
Q_{4}= \begin{cases} 
Q^1_{4}, & \text{if $\gamma(Q^1_{4}) \leq \gamma^{'} (Q^2_{4})$  }\\
Q^2_{4}, & \text{if $\gamma(Q^1_{4}) > \gamma^{'} (Q^2_{4})$}
\end{cases}
\end{equation}
From (1), clearly the length$(Q_{4})$ = min$ \{$length$(P_{1})$ + length$(P_{2})$, length $(P_{1}')$ + length$(P_{2}')\}$ = min$\{ k + j,(i - k-2j/3)(mod\,i) + j \}$. Here we denote $(i - k-2j/3)$ for $(i - k-2j/3)(mod\,i)$.

\begin{picture}(0,0)(-45,25)
\setlength{\unitlength}{3.75mm}

\drawpolygon(0,0)(12,0)(12,6)(0,6)


\drawline[AHnb=0](1,0)(1,6)
\drawline[AHnb=0](2,0)(2,6)
\drawline[AHnb=0](3,0)(3,6)
\drawline[AHnb=0](4,0)(4,6)
\drawline[AHnb=0](5,0)(5,6)
\drawline[AHnb=0](6,0)(6,6)
\drawline[AHnb=0](7,0)(7,6)
\drawline[AHnb=0](8,0)(8,6)
\drawline[AHnb=0](9,0)(9,6)
\drawline[AHnb=0](10,0)(10,6)
\drawline[AHnb=0](11,0)(11,6)


\drawline[AHnb=0](0,1)(12,1)
\drawline[AHnb=0](0,2)(12,2)
\drawline[AHnb=0](0,3)(12,3)
\drawline[AHnb=0](0,4)(12,4)
\drawline[AHnb=0](0,5)(12,5)

\drawline[AHnb=0](0,2)(1,3)
\drawline[AHnb=0](0,1)(2,3)
\drawline[AHnb=0](1,1)(3,3)
\drawline[AHnb=0](2,1)(4,3)

\drawline[AHnb=0](3,1)(5,3)
\drawline[AHnb=0](4,1)(6,3)
\drawline[AHnb=0](5,1)(7,3)
\drawline[AHnb=0](6,1)(8,3)
\drawline[AHnb=0](7,1)(9,3)
\drawline[AHnb=0](8,1)(10,3)
\drawline[AHnb=0](9,1)(11,3)
\drawline[AHnb=0](10,1)(12,3)
\drawline[AHnb=0](11,1)(12,2)

\drawline[AHnb=0](0,5)(1,6)
\drawline[AHnb=0](0,4)(2,6)
\drawline[AHnb=0](1,4)(3,6)
\drawline[AHnb=0](2,4)(4,6)
\drawline[AHnb=0](3,4)(5,6)
\drawline[AHnb=0](4,4)(6,6)
\drawline[AHnb=0](5,4)(7,6)
\drawline[AHnb=0](6,4)(8,6)
\drawline[AHnb=0](7,4)(9,6)
\drawline[AHnb=0](8,4)(10,6)
\drawline[AHnb=0](9,4)(11,6)
\drawline[AHnb=0](10,4)(12,6)
\drawline[AHnb=0](11,4)(12,5)

\put(-.2,-.5){\scriptsize {\tiny $x_1$}}
\put(.8,-.5){\scriptsize {\tiny $x_2$}}
\put(1.8,-.5){\scriptsize {\tiny $x_3$}}
\put(3.4,-.5){\scriptsize {\tiny $x_k$}}
\put(4.2,-.5){\scriptsize {\tiny $x_{k+1}$}}
\put(5.8,-.5){\scriptsize {\tiny $x_{k+2}$}}
\put(9.4,-.5){\scriptsize {\tiny $x_{i-1}$}}
\put(10.9,-.5){\scriptsize {\tiny $x_i$}}
\put(12,-.5){\scriptsize {\tiny $x_1$}}

\put(-1.5,6.4){\scriptsize {\tiny $x_{k+1}$}}
\put(.1,6.4){\scriptsize {\tiny $x_{k+2}$}}
\put(1.75,6.4){\scriptsize {\tiny $x_{k+3}$}}
\put(4.75,6.4){\scriptsize {\tiny $x_{k_1}$}}

\put(8.5,6.4){\scriptsize {\tiny $x_{k_2}$}}
\put(10.7,6.4){\scriptsize {\tiny $x_k$}}
\put(11.7,6.4){\scriptsize {\tiny $x_{k+1}$}}

\drawpolygon[fillcolor=black](5,0)(5,1)(7,3)(7,4)(9,6)(9,5.85)(7.15,4)(7.15,3)(5.15,1)(5.15,0)

\drawpolygon[fillcolor=black](5,0)(5,6)(5.15,6)(5.15,0)

\drawpolygon[fillcolor=blue](0,6)(5,6)(5,5.85)(0,5.85)

\drawpolygon[fillcolor=blue](9,6)(12,6)(12,5.85)(9,5.85)

\put(3.6,-1.6){\scriptsize {\tiny {\bf Figure 5.1.5:} $M(i,j,k)$ }} 

\end{picture}

\vspace{3.25cm}

From the notion of $M(i,j,k)$ representation, we prove the following lemma.

\begin{lem}\label{l4.1.7} Let $M_1$ be a DSEM of type $T_1$. Then $M_1$ admits an $M(i,j,k)$ representation iff the following holds: $(i)$ $i \geq 3$ and $j=3m$, where $m \in \mathbb{N}$, $(ii)$ $ij \geq 9 $, $(iii)$ $ 0 \leq k \leq i-1 $.
\end{lem}

\noindent{\bf Proof.} Let $M_1$ be a DSEM of type $T_1$ with $n$ vertices. By definition $M(i, j, k)$ of $M_1$ has $j$ number of $A_{1}$ type disjoint horizontal cycles of length $i$. Since all the vertices of $M_1$ lie in these cycles, the number of vertices in $M_1$ is $n = ij$. Clearly if $i \leq 2$, $M_1$ is not a map. So $i \geq 3$. If $j=1$ then $M_1$ is not a map and if $j =2$ then $M_1$ has no vertices with face-sequence $(3^6)$. If $j=3m+1$ or $3m+2$, where $ m \in \mathbb{N}$, then $2|V_{(3^6)}| \neq |V_{(3^3.4^2)}|$. So $j = 3m$ for $m \in \mathbb{N}$. Thus $n = ij \geq 9$. Since the length of the horizontal cycle is $i$, we get $ 0 \leq k \leq i-1 $. The converse part follows directly by constructing
$M(i, j, k)$ representation for the given values of $i, j$, and $k$.   \hfill $\Box$

\smallskip

Let $M_{t}$, $t = 1,2$, be DSEMs of type  $T_{1}$  on $n_t$ number of vertices with $n_1 = n_2$. Let $M_t(i_{t}, j_{t}, k_{t})$ be representation of $M_t$. Let $Q^t_{\alpha}$ be cycles of type $A_{\alpha}$ and $l^t_{\alpha}$ = length of the cycle of type $A_{\alpha}$, $\alpha \in \{1,2,3,4\}$, in $M_t(i_{t}, j_{t}, k_{t})$. We say that $M_t(i_{t}, j_{t}, k_{t})$ has cycle-type $(l^t_{1}, l^t_{2}, l^t_{3}, l^t_{4})$ if $l^t_{2} \leq l^t_{3}$ or $(l^t_{2}, l^t_{3}$, $l^t_{2}, l^t_{4})$ if $l^t_{3} < l^t_{2}$. Now, we show the following.

\begin{lem}\label{l4.1.8} The DSEMs $M_{1} \cong  M_{2}$ iff they have same cycle-type.
\end{lem}

\noindent{\bf Proof.} Let $M_{1}$ and $M_{2}$ be two DSEMs. Suppose the maps have same cycle-type. Then $l^1_{1} = l^2_{1}$, $\{l^1_{2}, l^1_{3}\} = \{l^2_{2}, l^2_{3}\}$ and $l^1_{4} = l^2_{4}$. To show that $M_1 \cong M_2$, it is enough to show that:

\smallskip

\noindent {\bf Claim.}  $M_1(i_1,j_1,k_1) \cong M_2(i_2,j_2,k_2)$.

\smallskip

By the definition, $M_t(i_{t},j_{t},k_{t})$ has $j_{t}$ horizontal cycles of type $A_{1}$, say $Q_{0}=C_{i_1}(w_{0,0},w_{0,1},\linebreak \ldots,w_{0,i_{1}-1}),Q_{1}=C_{i_1}(w_{1,0},w_{1,1},\ldots,w_{1,i_{1}-1}),\ldots,Q_{j_{1}-1}=C_{i_1}(w_{{j}_1-1,0},w_{{j}_1-1,1},\ldots,w_{j_{1}-1,i_{1}-1})$ in $M_{1}$ $(i_1,j_1,k_1)$ and $Q_{0}'=C_{i_2}(x_{0,0},x_{0,1},\ldots ,x_{0,i_{2}-1}), Q_{1}'=C_{i_2}(x_{1,0},x_{1,1},\ldots ,x_{1,i_{2}-1}),\ldots,Q_{j_{2}-1}'=C_{i_2}(x_{j_{2}-1,0}$, $x_{j_{2}-1,1},\ldots,x_{j_{2}-1,i_{2}-1})$ in $M_{2}(i_2,j_2,k_2)$. Then we have the following cases.

\smallskip
\noindent{\bf Case 1:} If $(i_{1}, j_{1}, k_{1})=(i_{2}, j_{2}, k_{2})$, $i_{1}=i_{2}, j_{1}=j_{2}, k_{1}=k_{2}$. Define an isomorphism  $f: V(M_1(i_1,j_1$, $k_1)) \to V(M_{2}(i_2,j_2,k_2))$ such that $f(w_{u,v})=x_{u,v}$ for $ 0 \leq u \leq j_{1}-1 $ and $ 0 \leq v \leq i_{1}-1 $. So, $M_1(i_1,j_1,k_1)$ $\cong M_2(i_2,j_2,k_2)$.

\smallskip
\noindent{\bf Case 2:} If $ i_{1} \neq i_{2} $, then it contradicts the fact that $l^1_{1}=l^2_{1}$. Thus $i_{1} = i_{2} $.

\smallskip
\noindent{\bf Case 3:} If $ j_{1} \neq j_{2} $, then $n_{1} = i_{1}j_{1} \neq i_{2}j_{2} = n_2$ as $ i_{1}=i_{2} $. A contradiction, as $n_{1} = n_{2}$. So, $j_{1}=j_{2}$.

\smallskip
\noindent{\bf Case 4:} If $k_{1} \neq k_{2} $. Since 
$l^1_{4} = l^2_{4}$, length$(Q^1_{4})$ = length$(Q^2_{4})$. This means min$\{k_1+j_1,i_1-k_1+j_1/3\}$ =  min$\{k_2+j_2,i_2-k_2+j_2/3\}$. Since $i_{1} = i_{2}$, $j_{1} = j_{2}$ and $k_{1} \neq k_{2}$, we get $k_1+j_1 \neq k_2+j_2 $ and $i_1-k_1+j_1/3 \neq i_2-k_2+j_2/3$. This gives that  $ k_1+j_1=i_2-k_2+j_2/3=i_1-k_2+j_1/3$, as $ i_{1}=i_{2}$ and  $j_{1}=j_{2}$. That is, $k_{2}=i_{1}-k_{1}-2j_1/3$. We identify $M_{2}(i_2,j_2,k_2)$ along the vertical boundary $P(x_{0,0},x_{1,0},\ldots , x_{{j_{2}-1},0},x_{0,k_{2}})$ and then cut along the path $P(x_{0,0},x_{1,0},x_{2,1},x_{3,2}, \ldots, x_{{j_{2}-1},2j_{2}/3-1}$, $x_{0,k_{2}+2j_{2}/3})$ of type $A_{2}$ through vertex $x_{0,0}$. This gives another representation of $M_2$, say $R$, with a map $f_1: V(M_{2}(i_2,j_2,k_2)) \to V(R)$ such that $f_1({x_{u,v}})=x_{u,(i_{2}-v+ \lfloor 2u/3 \rfloor)(mod\,i_{2})}$  for $0 \leq u \leq j_{2}-1$ and $0 \leq v \leq i_{2}-1 $. In $R$ the lower and upper horizontal cycles are $Q'=C_{i_2}(x_{0,0},x_{0,i_{2}-1},x_{0,i_{2}-2},\ldots , x_{0,1})$ and $Q''=C_{i_2}(x_{0,k_{2}+2j_{2}/3},x_{0,k_{2}+2j_{2}/3-1},\ldots, x_{0,k_{2}+2j_{2}/3+1)}$ respectively. The path $P(x_{0,0},x_{0,i_{2}-1}$, $x_{0,i_{2}-2}, \ldots, x_{0,k_{2}+2j_{2}/3})$ in $Q'$ has length $i_{2}-k_{2}-2j_{2}/3$.	Clearly $R$ has $j_{2}$ number of horizontal cycles of length $i_{2}$. So, $R = M_2(i_{2},j_{2},i_{2}-k_{2}-2j_{2}/3)$. Note that  $i_{2}-k_{2}-2j_{2}/3=i_{2}-(i_{1}-k_{1}-2j_{1}/3)-2j_{2}/3=k_{1}$, as $i_{1}=i_{2}$ and $j_{1}=j_{2}$. So,  $M_2(i_{2},j_{2},i_{2}-k_{2}-2j_{2}/3) = M_1(i_{1}, j_{1}, k_{1})$. Therefore by $f$, $M_1(i_1,j_1,k_1) \cong M_2(i_2,j_2,k_2)$. So, by Cases 1-4, the claim follows. Hence, $M_{1} \cong M_{2}$. 

Conversely, let  $M_{1} \cong M_{2}$ by an isomorphism $f$. Let $Q^1_{\alpha}$ and $Q^2_{\alpha}$ be  cycles of type $A_{\alpha}$ for $\alpha \in \{ 1, 2, 3, 4\}$ in $M_{1}$ and $M_{2}$ respectively. Let $f :V(M_{1}) \to V(M_{2})$ be such that $f(Q^1_{\alpha})$ = $Q^2_{\alpha}$. Since $f$ is an isomorphism, length$(Q^1_{\alpha})$ = length($f(Q^1_{\alpha})$) = length($Q^2_{\alpha}$).  Hence, $M_{1}$ and $M_{2}$ have the same cycle-type. This completes the proof. \hfill $\Box$

\smallskip

By Lemmas \ref{l4.1.7} -
\ref{l4.1.8}, one can compute and classify DSEMs of type $T_1$ for any $|V(M_1)|$. A tabular list of the DSEMs for first four admissible values of $|V(M_1)|$, i.e., $|V(M_1)| = 9, 12,15, 18$ is given in Table \ref{table:1}. For $|V(M_1)| = 12$, the computation is illustrated as follows. 


\begin{eg} \normalfont Let $M_1$ be a DSEM of type $T_1$ with $n$ vertices. Then $M_1$ has an $M(i,j,k)$ representation, with $n = ij \geq 9$, $j=3m$, where $m \in \mathbb{N}$ and $0 \leq k \leq i-1 $. If $n = 12$ and $j = 3$, we have  $i = 4$ and $k = 0,1,2,3 $ by Lemma \ref{l4.1.7}. So, $M(i, j, k) = M(4, 3, 0), M(4, 3, 1), M(4, 3, 2)$ and $M(4, 3, 3)$, see Figures 5.1.6, 5.1.7, 5.1.8, and 5.1.9, respectively. In $M(4, 3, 0)$, $Q^1_1 = C_4(v_{1}, v_{2}, v_{3},v_{4})$ is $A_{1}$ type cycle, $Q^1_2 = C_6(v_{1},v_{5},u_{2},v_{3},v_{7},u_{4})$ and $Q^1_3 = C_3(v_{1},v_{5},u_{1})$ are two $A_{2}$ type cycles and $Q^1_4 = C_3(v_{1},v_{5},u_{1})$ is $A_{4}$ type cycle. In $M(4, 3, 1)$, $Q^2_1 = C_4(w_{1}, w_{2}, w_{3},w_{4})$ is $A_{1}$ type cycle, $Q^2_2 = C_{12}(w_{1},w_{5},u_{2},w_{4},w_{8},u_{1},w_{3},w_{7},u_{4},w_{2},w_{6},	u_{3})$ and $Q^2_3 = C_{12}(w_{1}, w_{5}, u_{1}, w_{2}, w_{6}, u_{2}, w_{3}, w_{7}$, $u_{3}, w_{4}, w_{8}, u_{4})$ are two $A_{2}$ type cycle and $Q^2_4 = C_4(w_{2},w_{6}, u_{3},w_{1})$ is $A_{4}$ type cycle. In $M(4, 3, 2)$, $Q^3_1 = C_4(x_{1}, x_{2}, x_{3},x_{4})$ is a cycle of type $A_{1}$, $Q^3_2 = C_3(x_{1},x_{5},u_{2})$ and $Q^3_3 = C_6(x_{1},x_{5},u_{1},x_{3},x_{7}$, $u_{3})$ are two $A_{2}$ type cycles and $Q^3_4 = C_3(x_{3},x_{7},u_{4})$ is $A_{4}$ type cycle. In $M(4, 3, 3)$, $Q^4_1 = C_4(z_{1}, z_{2}, z_{3},z_{4})$ is $A_{1}$ type cycle, $Q^4_2 = C_{12}(z_{1},z_{5},u_{2},z_{2},z_{6}, u_{3},z_{3},z_{7},u_{4},z_{4}$, $z_{8}, u_{1})$ and $Q^4_{4} = C_{12}(z_{1},z_{5},u_{1},z_{4},z_{8},u_{4}, z_{3},z_{7}$, $u_{3},z_{2},z_{6},	u_{2})$ are two $A_{2}$ type cycles and $Q^4_{4} =C_6(z_{4},z_{8},u_{1},z_{1},z_{2},z_{3})$ is $A_{4}$ type cycle.	
	
By Lemma \ref{l4.1.6}, $M_1$ has $A_{1}$ type cycles with unique length and $A_{2}$ type cycles with at most two different lengths. Since length$(Q^1_4) \neq$ length$(Q^r_{4})$ for $r \in \{2,4\}$, $M(4,3,0) \ncong M(4,3,1), M(4,3,3)$. Also, $M(4,3,1) \ncong M(4,3,2)$ as length$(Q^2_4)$ $\neq$ length$(Q^3_4)$ and $M(4,3,2) \ncong M(4,3,3)$ as  length$\linebreak(Q^3_4)$ $\neq$ length$(Q^4_{4})$. Observe that, length$(Q^1_1)$ = length$(Q^3_1)$, \{length$(Q^1_2)$, length$(Q^1_3)$\} =  \{length$(Q^3_2)$, length$(Q^3_3)$\} and length$(Q^1_4)$ = length$(Q^4_{4})$. Now identifying $M(4,3,2)$ along the vertical boundary and cutting along the path $P(x_{3},x_{7},u_{4},x_{3})$ leads to Figure 5.1.10, i.e., $M(4, 3, 0)$. So, by Lemma \ref{l4.1.8}, $M(4,3,0) \cong M(4,3,2)$. Thus, there are three DSEMs, up to isomorphism, of type $T_1$ with $12$ vertices on the torus. These are $M(4,3,0)$, $M(4,3,1)$, $M(4,3,3)$. 	
	

	\begin{picture}(0,0)(2.4,15.5)
	\setlength{\unitlength}{5mm}
	
	\drawpolygon(0,0)(4,0)(4,3)(0,3)
	
	
	\drawline[AHnb=0](1,0)(1,3)
	\drawline[AHnb=0](2,0)(2,3)
	\drawline[AHnb=0](3,0)(3,3)
	
	
	\drawline[AHnb=0](0,1)(4,1)
	\drawline[AHnb=0](0,2)(4,2)
	
	
	\drawline[AHnb=0](0,1)(2,3)
	\drawline[AHnb=0](1,1)(3,3)
	\drawline[AHnb=0](2,1)(4,3)
	\drawline[AHnb=0](3,1)(4,2)
	\drawline[AHnb=0](0,2)(1,3)
	
	\put(-.2,-.4){\scriptsize {\tiny $v_1$}}
	\put(.8,-.4){\scriptsize {\tiny $v_2$}}
	\put(1.8,-.4){\scriptsize {\tiny $v_3$}}
	\put(2.8,-.4){\scriptsize {\tiny $v_4$}}
	\put(3.8,-.4){\scriptsize {\tiny $v_1$}}
	
	\put(-.5,.7){\scriptsize {\tiny $v_5$}}
	\put(1.05,.7){\scriptsize {\tiny $v_6$}}
	\put(2.05,.7){\scriptsize {\tiny $v_7$}}
	\put(3.05,.7){\scriptsize {\tiny $v_8$}}
	\put(4.1,.7){\scriptsize {\tiny $v_5$}}
	
	\put(-.5,1.7){\scriptsize {\tiny $u_1$}}
	\put(1.05,1.7){\scriptsize {\tiny $u_2$}}
	\put(2.05,1.7){\scriptsize {\tiny $u_3$}}
	\put(3.05,1.7){\scriptsize {\tiny $u_4$}}
	\put(4.15,1.7){\scriptsize {\tiny $u_1$}}

	\put(-.2,3.15){\scriptsize {\tiny $v_1$}}
	\put(.85,3.15){\scriptsize {\tiny $v_2$}}
	\put(1.85,3.15){\scriptsize {\tiny $v_3$}}
	\put(2.85,3.15){\scriptsize {\tiny $v_4$}}
	\put(3.8,3.15){\scriptsize {\tiny $v_1$}}
	
	\put(-.5,-1.2){\scriptsize {\tiny Figure 5.1.6: $M(4,3,0)$}}

	\end{picture}

	\begin{picture}(0,0)(-30,11)
	\setlength{\unitlength}{5mm}
	
	\drawpolygon(0,0)(4,0)(4,3)(0,3)
	
	
	\drawline[AHnb=0](1,0)(1,3)
	\drawline[AHnb=0](2,0)(2,3)
	\drawline[AHnb=0](3,0)(3,3)
	
	
	\drawline[AHnb=0](0,1)(4,1)
	\drawline[AHnb=0](0,2)(4,2)
	
	
	\drawline[AHnb=0](0,1)(2,3)
	\drawline[AHnb=0](1,1)(3,3)
	\drawline[AHnb=0](2,1)(4,3)
	\drawline[AHnb=0](3,1)(4,2)
	\drawline[AHnb=0](0,2)(1,3)
	
	\put(-.2,-.4){\scriptsize {\tiny $w_1$}}
	\put(.8,-.4){\scriptsize {\tiny $w_2$}}
	\put(1.8,-.4){\scriptsize {\tiny $w_3$}}
	\put(2.8,-.4){\scriptsize {\tiny $w_4$}}
	\put(3.8,-.4){\scriptsize {\tiny $w_1$}}
	
	\put(-.6,.7){\scriptsize {\tiny $w_5$}}
	\put(1.05,.7){\scriptsize {\tiny $w_6$}}
	\put(2.05,.7){\scriptsize {\tiny $w_7$}}
	\put(3.05,.7){\scriptsize {\tiny $w_8$}}
	\put(4.1,.7){\scriptsize {\tiny $w_5$}}
	
	\put(-.6,1.7){\scriptsize {\tiny $u_1$}}
	\put(1.05,1.7){\scriptsize {\tiny $u_2$}}
	\put(2.05,1.7){\scriptsize {\tiny $u_3$}}
	\put(3.05,1.7){\scriptsize {\tiny $u_4$}}
	\put(4.1,1.7){\scriptsize {\tiny $u_1$}}

	\put(-.2,3.15){\scriptsize {\tiny $w_2$}}
	\put(.85,3.15){\scriptsize {\tiny $w_3$}}
	\put(1.85,3.15){\scriptsize {\tiny $w_4$}}
	\put(2.85,3.15){\scriptsize {\tiny $w_1$}}
	\put(3.8,3.15){\scriptsize {\tiny $w_2$}}
	
	\put(-.2,-1.2){\scriptsize {\tiny Figure 5.1.7: $M(4,3,1)$}}

	\end{picture}
	
	\begin{picture}(0,0)(-62,6.5)
	\setlength{\unitlength}{5mm}
	
	\drawpolygon(0,0)(4,0)(4,3)(0,3)
	
	
	\drawline[AHnb=0](1,0)(1,3)
	\drawline[AHnb=0](2,0)(2,3)
	\drawline[AHnb=0](3,0)(3,3)
	
	
	\drawline[AHnb=0](0,1)(4,1)
	\drawline[AHnb=0](0,2)(4,2)
	
	
	\drawline[AHnb=0](0,1)(2,3)
	\drawline[AHnb=0](1,1)(3,3)
	\drawline[AHnb=0](2,1)(4,3)
	\drawline[AHnb=0](3,1)(4,2)
	\drawline[AHnb=0](0,2)(1,3)
	
	\put(-.2,-.4){\scriptsize {\tiny $x_1$}}
	\put(.8,-.4){\scriptsize {\tiny $x_2$}}
	\put(1.8,-.4){\scriptsize {\tiny $x_3$}}
	\put(2.8,-.4){\scriptsize {\tiny $x_4$}}
	\put(3.8,-.4){\scriptsize {\tiny $x_1$}}
	
	\put(-.6,.7){\scriptsize {\tiny $x_5$}}
	\put(1.05,.7){\scriptsize {\tiny $x_6$}}
	\put(2.05,.7){\scriptsize {\tiny $x_7$}}
	\put(3.05,.7){\scriptsize {\tiny $x_8$}}
	\put(4.1,.7){\scriptsize {\tiny $x_5$}}
	
	\put(-.6,1.7){\scriptsize {\tiny $u_1$}}
	\put(1.05,1.7){\scriptsize {\tiny $u_2$}}
	\put(2.05,1.7){\scriptsize {\tiny $u_3$}}
	\put(3.05,1.7){\scriptsize {\tiny $u_4$}}
	\put(4.1,1.7){\scriptsize {\tiny $u_1$}}
	
	\put(-.2,3.15){\scriptsize {\tiny $x_3$}}
	\put(.85,3.15){\scriptsize {\tiny $x_4$}}
	\put(1.85,3.15){\scriptsize {\tiny $x_1$}}
	\put(2.85,3.15){\scriptsize {\tiny $x_2$}}
	\put(3.8,3.15){\scriptsize {\tiny $x_3$}}
	
	\put(-.4,-1.2){\scriptsize {\tiny Figure 5.1.8: $M(4,3,2)$}}

	\end{picture}

	\begin{picture}(0,0)(-93,1.5)
	\setlength{\unitlength}{5mm}
	
	\drawpolygon(0,0)(4,0)(4,3)(0,3)
	
	
	\drawline[AHnb=0](1,0)(1,3)
	\drawline[AHnb=0](2,0)(2,3)
	\drawline[AHnb=0](3,0)(3,3)
	
	
	\drawline[AHnb=0](0,1)(4,1)
	\drawline[AHnb=0](0,2)(4,2)
	
	
	\drawline[AHnb=0](0,1)(2,3)
	\drawline[AHnb=0](1,1)(3,3)
	\drawline[AHnb=0](2,1)(4,3)
	\drawline[AHnb=0](3,1)(4,2)
	\drawline[AHnb=0](0,2)(1,3)
	
	\put(-.2,-.4){\scriptsize {\tiny $z_1$}}
	\put(.8,-.4){\scriptsize {\tiny $z_2$}}
	\put(1.8,-.4){\scriptsize {\tiny $z_3$}}
	\put(2.8,-.4){\scriptsize {\tiny $z_4$}}
	\put(3.8,-.4){\scriptsize {\tiny $z_1$}}
	
	\put(-.6,.7){\scriptsize {\tiny $z_5$}}
	\put(1.05,.7){\scriptsize {\tiny $z_6$}}
	\put(2.05,.7){\scriptsize {\tiny $z_7$}}
	\put(3.05,.7){\scriptsize {\tiny $z_8$}}
	\put(4.1,.7){\scriptsize {\tiny $z_5$}}
	
	\put(-.6,1.7){\scriptsize {\tiny $u_1$}}
	\put(1.05,1.7){\scriptsize {\tiny $u_2$}}
	\put(2.05,1.7){\scriptsize {\tiny $u_3$}}
	\put(3.05,1.7){\scriptsize {\tiny $u_4$}}
	\put(4.1,1.7){\scriptsize {\tiny $u_1$}}
	
	\put(-.2,3.15){\scriptsize {\tiny $z_4$}}
	\put(.85,3.15){\scriptsize {\tiny $z_1$}}
	\put(1.85,3.15){\scriptsize {\tiny $z_2$}}
	\put(2.85,3.15){\scriptsize {\tiny $z_3$}}
	\put(3.8,3.15){\scriptsize {\tiny $z_4$}}
	
	\put(-.2,-1.2){\scriptsize {\tiny Figure 5.1.9: $(4,3,3)$}}

	\end{picture}
	
	\begin{picture}(0,0)(-124,-2.85)
	\setlength{\unitlength}{5mm}
	
	\drawpolygon(0,0)(4,0)(4,3)(0,3)
	
	
	\drawline[AHnb=0](1,0)(1,3)
	\drawline[AHnb=0](2,0)(2,3)
	\drawline[AHnb=0](3,0)(3,3)
	
	
	\drawline[AHnb=0](0,1)(4,1)
	\drawline[AHnb=0](0,2)(4,2)
	
	
	\drawline[AHnb=0](0,1)(2,3)
	\drawline[AHnb=0](1,1)(3,3)
	\drawline[AHnb=0](2,1)(4,3)
	\drawline[AHnb=0](3,1)(4,2)
	\drawline[AHnb=0](0,2)(1,3)
	
	\put(-.2,-.4){\scriptsize {\tiny $x_3$}}
	\put(.8,-.4){\scriptsize {\tiny $x_2$}}
	\put(1.8,-.4){\scriptsize {\tiny $x_1$}}
	\put(2.8,-.4){\scriptsize {\tiny $x_4$}}
	\put(3.8,-.4){\scriptsize {\tiny $x_3$}}
	
	\put(-.6,.7){\scriptsize {\tiny $x_7$}}
	\put(1.05,.7){\scriptsize {\tiny $x_6$}}
	\put(2.05,.7){\scriptsize {\tiny $x_5$}}
	\put(3.05,.7){\scriptsize {\tiny $x_8$}}
	\put(4.1,.7){\scriptsize {\tiny $x_7$}}
	
	\put(-.6,1.7){\scriptsize {\tiny $u_4$}}
	\put(1.05,1.7){\scriptsize {\tiny $u_3$}}
	\put(2.05,1.7){\scriptsize {\tiny $u_2$}}
	\put(3.05,1.7){\scriptsize {\tiny $u_1$}}
	\put(4.1,1.7){\scriptsize {\tiny $u_4$}}
	
	\put(-.2,3.15){\scriptsize {\tiny $x_3$}}
	\put(.85,3.15){\scriptsize {\tiny $x_2$}}
	\put(1.85,3.15){\scriptsize {\tiny $x_1$}}
	\put(2.85,3.15){\scriptsize {\tiny $x_4$}}
	\put(3.8,3.15){\scriptsize {\tiny $x_3$}}
	
	\put(-.3,-1.22){\scriptsize {\tiny Figure 5.1.10: $M(4,3,0)$}}

	\end{picture}
	
\end{eg}

\vspace{.1cm}

\begin{center}
	
	\noindent \textbf{Table \ref{table:1}} : DSEMs of type $T_1$  on the torus for $|V(M)| \leq 18$
	
	\renewcommand{\arraystretch}{1.1}
	\begin{tabular}{ |p{1.2cm}|p{4cm}|p{3.8cm}|p{3cm}| }
		
		\hline
		$|V(M)|$ & Isomorphic classes & Length of cycles & No of maps\\
		\hline
		9 &  $M(3,3,0)$, $M(3,3,1)$ & $( 3, \{ 3, 9 \}, 3 )$ & 2\\
		
		\cline{2-3}
		& $M(3,3,2)$ & $( 3, \{ 9, 9 \}, 5 )$ & \\	
		
		\hline

		\hline
		12 &  $M(4,3,0)$, $M(4,3,2)$ & $( 4, \{ 3, 6 \}, 3 )$& 3 \\
		
		\cline{2-3}
		
		& $M(4,3,1)$ & $( 4, \{ 12,12 \}, 4 )$ & \\
		
		\cline{2-3}
		
		&$M(4,3,3)$ & $( 4, \{ 12, 12 \}, 6 )$ & \\
		
		\hline

		15 & $M(5,3,0)$, $M(5,3,3)$ & $( 5, \{ 3, 15 \}, 3 )$ &3\\
		
		\cline{2-3}
		
		& $M(5,3,1)$, $M(5,3,2)$  &$( 5, \{ 15, 15 \}, 4 )$ &\\
		
		\cline{2-3}
		
		& $M(5,3,4)$  & $( 5, \{ 15, 15 \}, 7 )$  &\\
		
		\hline
		
	\end{tabular}

\end{center}

\begin{center}
\begin{tabular}{ |p{1.2cm}|p{4cm}|p{3.8cm}|p{3cm}| }

		\hline
		18 & $M(6,3,0)$, $M(6,3,4)$ & $( 6, \{ 3, 9 \}, 3 )$ & 6\\
		
		\cline{2-3}
		
		&$M(6,3,1)$, $M(6,3,3)$ &$( 6, \{ 6, 18 \}, 4 )$&\\
		
		\cline{2-3}
		
		& $M(6,3,2)$ &$( 6, \{ 9, 9 \}, 5)$&\\
		
		\cline{2-3}
		
		&$M(6,3,5)$ & $( 6, \{ 18, 18 \}, 8 )$&  \\
		
		\cline{2-3}
		
		& $M(3,6,0)$, $M(3,6,2)$
		&$(3, \{ 6, 18 \}, 6)$ & \\
		
		\cline{2-3}

		&$M(3,6,1)$ & $( 3, \{ 18, 18 \}, 7 )$ &  \\

		\hline

	\end{tabular}
	\label{table:1}
\end{center}

In the subsequent subsections, we proceed in a similar way. For each type DSEM $M$ with vertex set $|V(M)|$, we construct an $M(i,j,k)$ representation by considering suitable non-homologous cycles and finding admissible values of $i,j,k \in \mathbb{N} \cup \{0\}$. Further, we determine the classes of such representations by defining isomorphism between them, if exist. This gives the number of representations or maps, up to the isomorphism, on the $|V(M_1)|$.

\subsection{DSEMs of type  $T_2= [f_1^{({f_1}^2.{f_2}^4)}:f_2^{({f_1}^5.{f_2}^2)}]$, where $(f_1, f_2 )=( (3^6),(3^3.4^2))$ }\label{s4.2}
Let $M_2$ be a DSEM of type $T_2$ with the vertex set $V(M_2)$. Then for the existence of $M_2$ we see easily that $|V_{(3^6)}| = |V_{(3^3.4^2)}|$.   Now define following three types of paths in $M_2$ as follows.




\begin{deff} \normalfont
	
	A path $P_{1} = P( \ldots, y_{i-1},y_{i},y_{i+1}, \ldots)$ in $M_2$ is of type $B_{1}$ if either of the conditions follows: 
	$(i)$ every vertex of $P_1$ has the face-sequence $(3^6)$ or
	$(ii)$ each vertex of $P_1$ has the face-sequence $(3^3.4^2)$ such that all the triangles (squares) that incident on its inner vertices lie on the one side  of the path $P_1$. For example, see thick black colored paths in Figure 5.2.1.
	
\end{deff}


\begin{deff} \normalfont A path $P_{2} = P(\ldots, z_{i-1},z_{i},z_{i+1}, \ldots)$ in $M_2$, such that $z_{i-1}, z_{i}, z_{i+1}$ are inner vertices of $P_{2}$ or an extended path of $P_{2}$, is of type $B_{2}$ (shown by red colored  paths in Figure 5.2.1), if either of the following four conditions follows for each vertex of the path.

	\begin{enumerate} 
		\item If ${\rm lk}(z_{i})=C_7(\boldsymbol{m},z_{i-1},\boldsymbol{n},o,z_{i+1},p,q)$ and ${\rm lk}(z_{i-1})=C_7(m,z_{i-2},r,n,\boldsymbol{o},z_i,\boldsymbol{q})$, then ${\rm lk}(z_{i+1})=C_6(s,z_{i+2},t,p,z_i,o), {\rm lk}(z_{i+2})=C_6(z_{i+1},s,u,z_{i+3},v,t)$. 
		
		\item If ${\rm lk}(z_{i})=C_7(\boldsymbol{m},z_{i+1},\boldsymbol{n},o,z_{i-1},p,q)$ and ${\rm lk}(z_{i-1})=C_6(z_{i-2},s,p,z_{i},o,r)$, then ${\rm lk}(z_{i+1})=C_7(\boldsymbol{o},z_{i}, \boldsymbol{q},m,z_{i+2},t,n)$, ${\rm lk}(z_{i+2})=C_6(z_{i+1},m,u,z_{i+3},v,t)$.
		
		\item If  ${\rm lk}(z_{i})=C_6(z_{i-1},m,n,z_{i+1},o,p)$ and ${\rm lk}(z_{i-1})=C_7(\boldsymbol{r},z_{i-2}, \boldsymbol{s},m,z_{i},p,q)$, then ${\rm lk}(z_{i+1})=C_6(z_{i},n,t,z_{i+2},u,o)$, ${\rm lk}(z_{i+2})=C_7(\boldsymbol{v},z_{i+3}, \boldsymbol{w},u,z_{i+1},t,x)$.
		
		\item If ${\rm lk}(z_{i})=C_6(z_{i-1},m,n,z_{i+1},o,p)$ and ${\rm lk}(z_{i-1})=C_6(z_{i-2},r,m,z_{i},p,q)$, then  
		${\rm lk}(z_{i+1})=C_7(\boldsymbol{s},z_{i+2},\boldsymbol{t},o,z_{i},n,u)$, ${\rm lk}(z_{i+2})=C_7(\boldsymbol{o}, z_{i+1},\boldsymbol{u},s,z_{i+3},v,t)$.
		
	\end{enumerate}
\end{deff}

\begin{deff} \normalfont A path $P_{3} = P(\ldots, w_{i-1},w_{i},w_{i+1}, \ldots)$ in $M_2$, such that $w_{i-1}, w_{i}, w_{i+1}$ are inner vertices of $P_{3}$ or an extended path of $P_{3}$, is of type $B_{3}$ (shown by green colored  paths in Figure 5.2.1), if either of the following four conditions follows for each vertex of the path.
	
	\begin{enumerate} \item  If ${\rm lk}(w_{i})=C_7(\boldsymbol{m},w_{i-1},\boldsymbol{n},o,p,w_{i+1},q)$ and ${\rm lk}(w_{i-1})=C_7(r,w_{i-2},n,\boldsymbol{o},w_i,\boldsymbol{q},m)$, then \linebreak ${\rm lk}(w_{i+1})=C_6(w_{i},p,s,w_{i+2},t,q), {\rm lk}(w_{i+2})=C_6(w_{i+1},s,u,w_{i+3},v,t)$.

		\item If ${\rm lk}(w_{i})=C_7(\boldsymbol{m},w_{i+1},\boldsymbol{n},o,p,w_{i-1},q)$ and ${\rm lk}(w_{i-1})=C_6(w_{i-2},s,q,w_{i},p,r)$, then ${\rm lk}(w_{i+1})=C_7(\boldsymbol{o},w_{i}, \boldsymbol{q},m,t,w_{i+2},n)$, ${\rm lk}(w_{i+2})=C_6(w_{i+1},t,u,w_{i+3},v,n)$.
		
		\item If  ${\rm lk}(w_{i})=C_6(w_{i-1},m,n,w_{i+1},o,p)$ and ${\rm lk}(w_{i-1})=C_7(\boldsymbol{q},w_{i-2}, \boldsymbol{r},s,m,w_{i},p)$, then ${\rm lk}(w_{i+1})=C_6(w_{i},n,t,w_{i+2},u,o)$, ${\rm lk}(w_{i+2})=C_7(\boldsymbol{v},w_{i+3}, \boldsymbol{w},x,u,w_{i+1},t)$.
		
		\item  If ${\rm lk}(w_{i})=C_6(w_{i-1},m,n,w_{i+1},o,p)$ and ${\rm lk}(w_{i-1})=C_6(w_{i-2},r,m,w_{i},p,q)$, then  
		${\rm lk}(w_{i+1})=C_7(\boldsymbol{s},w_{i+2},\boldsymbol{t},u,o,w_{i},n)$, ${\rm lk}(w_{i+2})=C_7(\boldsymbol{u}, w_{i+1},\boldsymbol{n},s,v,w_{i+3},t)$.
	\end{enumerate}
\end{deff}

\begin{picture}(0,0)(-15,12.5)
\setlength{\unitlength}{4.5mm}


\drawline[AHnb=0](-1,-2)(-1,-3)
\drawline[AHnb=0](0,-2)(0,-3)
\drawline[AHnb=0](1,-2)(1,-3)
\drawline[AHnb=0](2,-2)(2,-3)
\drawline[AHnb=0](3,-2)(3,-3)
\drawline[AHnb=0](4,-2)(4,-3)
\drawline[AHnb=0](5,-2)(5,-3)

\drawline[AHnb=0](-.5,1)(-.5,2)
\drawline[AHnb=0](.5,1)(.5,2)
\drawline[AHnb=0](1.5,1)(1.5,2)
\drawline[AHnb=0](2.5,1)(2.5,2)
\drawline[AHnb=0](3.5,1)(3.5,2)
\drawline[AHnb=0](4.5,1)(4.5,2)


\drawline[AHnb=0](-1,-3)(5,-3)
\drawline[AHnb=0](-1,-2)(5,-2)
\drawline[AHnb=0](-1,1)(5,1)
\drawline[AHnb=0](-1,2)(5,2)
\drawline[AHnb=0](-1.11,0)(5.11,0)


\drawline[AHnb=0](-.75,.5)(-.5,1)
\drawline[AHnb=0](-1,-2)(.5,1)
\drawline[AHnb=0](0,-2)(1.5,1)
\drawline[AHnb=0](1,-2)(2.5,1)
\drawline[AHnb=0](2,-2)(3.5,1)
\drawline[AHnb=0](3,-2)(4.5,1)
\drawline[AHnb=0](4,-2)(4.75,-.5)

\drawline[AHnb=0](-.75,-.5)(0,-2)
\drawline[AHnb=0](-.5,1)(1,-2)
\drawline[AHnb=0](.5,1)(2,-2)
\drawline[AHnb=0](1.5,1)(3,-2)
\drawline[AHnb=0](2.5,1)(4,-2)
\drawline[AHnb=0](3.5,1)(5,-2)
\drawline[AHnb=0](4.5,1)(4.75,0.5)


\drawline[AHnb=0](.5,2)(.75,2.5)
\drawline[AHnb=0](.5,2)(.25,2.5)

\drawline[AHnb=0](-.5,2)(-.75,2.5)
\drawline[AHnb=0](-.5,2)(-.25,2.5)

\drawline[AHnb=0](1.5,2)(1.75,2.5)
\drawline[AHnb=0](1.5,2)(1.25,2.5)

\drawline[AHnb=0](3.5,2)(3.75,2.5)
\drawline[AHnb=0](3.5,2)(3.25,2.5)

\drawline[AHnb=0](4.5,2)(4.75,2.5)
\drawline[AHnb=0](4.5,2)(4.25,2.5)

\drawline[AHnb=0](2.5,2)(2.25,2.5)
\drawline[AHnb=0](2.5,2)(2.75,2.5)

\drawline[AHnb=0](-1,-3)(-.75,-3.5)
\drawline[AHnb=0](-1,-3)(-1.25,-3.5)

\drawline[AHnb=0](1,-3)(1.25,-3.5)
\drawline[AHnb=0](1,-3)(.75,-3.5)

\drawline[AHnb=0](0,-3)(-.25,-3.5)
\drawline[AHnb=0](0,-3)(.25,-3.5)

\drawline[AHnb=0](2,-3)(2.25,-3.5)
\drawline[AHnb=0](2,-3)(1.75,-3.5)

\drawline[AHnb=0](4,-3)(4.25,-3.5)
\drawline[AHnb=0](4,-3)(3.75,-3.5)

\drawline[AHnb=0](5,-3)(5.25,-3.5)
\drawline[AHnb=0](5,-3)(4.75,-3.5)

\drawline[AHnb=0](3,-3)(2.75,-3.5)
\drawline[AHnb=0](3,-3)(3.25,-3.5)


\drawpolygon[fillcolor=green](1.25,-3.5)(1,-3)(1,-2)(-.5,1)(-.5,2)(-.75,2.5)(-.65,2.5)(-.4,2)(-.4,1)(1.1,-2)(1.1,-3)(1.35,-3.5)

\drawpolygon[fillcolor=green](2.25,-3.5)(2,-3)(2,-2)(.5,1)(.5,2)(.25,2.5)(.35,2.5)(.6,2)(.6,1)(2.1,-2)(2.1,-3)(2.35,-3.5)

\drawpolygon[fillcolor=green](3.25,-3.5)(3,-3)(3,-2)(1.5,1)(1.5,2)(1.25,2.5)(1.35,2.5)(1.6,2)(1.6,1)(3.1,-2)(3.1,-3)(3.35,-3.5)

\drawpolygon[fillcolor=green](4.25,-3.5)(4,-3)(4,-2)(2.5,1)(2.5,2)(2.25,2.5)(2.35,2.5)(2.6,2)(2.6,1)(4.1,-2)(4.1,-3)(4.35,-3.5)

\drawpolygon[fillcolor=green](5.25,-3.5)(5,-3)(5,-2)(3.5,1)(3.5,2)(3.25,2.5)(3.35,2.5)(3.6,2)(3.6,1)(5.1,-2)(5.1,-3)(5.35,-3.5)

\drawpolygon[fillcolor=red](-1.25,-3.5)(-1,-3)(-1,-2)(.5,1)(.5,2)(.75,2.5)(.65,2.5)(.4,2)(.4,1)(-1.1,-2)(-1.1,-3)(-1.35,-3.5)

\drawpolygon[fillcolor=red](-.25,-3.5)(0,-3)(0,-2)(1.5,1)(1.5,2)(1.75,2.5)(1.65,2.5)(1.4,2)(1.4,1)(-.1,-2)(-.1,-3)(-.35,-3.5)

\drawpolygon[fillcolor=red](.75,-3.5)(1,-3)(1,-2)(2.5,1)(2.5,2)(2.75,2.5)(2.65,2.5)(2.4,2)(2.4,1)(.9,-2)(.9,-3)(.65,-3.5)

\drawpolygon[fillcolor=red](1.75,-3.5)(2,-3)(2,-2)(3.5,1)(3.5,2)(3.75,2.5)(3.65,2.5)(3.4,2)(3.4,1)(1.9,-2)(1.9,-3)(1.65,-3.5)

\drawpolygon[fillcolor=red](2.75,-3.5)(3,-3)(3,-2)(4.5,1)(4.5,2)(4.75,2.5)(4.65,2.5)(4.4,2)(4.4,1)(2.9,-2)(2.9,-3)(2.65,-3.5)

\drawpolygon[fillcolor=red](3.75,-3.5)(4,-3)(4,-2)(5.25,.5)(5.14,.5)(3.9,-2)(3.9,-3)(3.65,-3.5)

\drawpolygon[fillcolor=black](-1.5,-3)(5.25,-3)(5.25,-3.1)(-1.5,-3.1)

\drawpolygon[fillcolor=black](-1.5,-2)(5.25,-2)(5.25,-2.1)(-1.5,-2.1)

\drawpolygon[fillcolor=black](-1,-1)(5,-1)(5,-1.1)(-1,-1.1)

\drawpolygon[fillcolor=black](-1.12,0)(5.1,0)(5.1,.1)(-1.12,.1)

\drawpolygon[fillcolor=black](-1,1)(5,1)(5,1.1)(-1,1.1)

\drawpolygon[fillcolor=black](-1,2)(5,2)(5,2.1)(-1,2.1)

\put(2.35,2.65) {\scriptsize {\tiny $B_2$}}
\put(1,2.65) {\scriptsize {\tiny $B_3$}}
\put(5.5,-2.1) {\scriptsize {\tiny $B_1$}}
\put(5.5,-1.1) {\scriptsize {\tiny $B_1$}}

\put(-2.75,-4.4){\scriptsize {\tiny {\bf Figure 5.2.1:} Paths of types $B_1, B_2, B_3$ }} 
\end{picture}

\begin{picture}(0,0)(-75,19)
\setlength{\unitlength}{6.15mm}

\drawpolygon(0,0)(3,0)(3,4)(0,4)
\drawpolygon(4.5,0)(7.5,0)(7.5,4)(4.5,4)
\drawpolygon(9,0)(11,0)(11,4)(9,4)


\drawline[AHnb=0](1,0)(1,4)
\drawline[AHnb=0](2,0)(2,4)

\drawline[AHnb=0](5.5,0)(5.5,4)
\drawline[AHnb=0](6.5,0)(6.5,4)

\drawline[AHnb=0](10,0)(10,4)


\drawline[AHnb=0](0,1)(3,1)
\drawline[AHnb=0](0,2)(3,2)
\drawline[AHnb=0](0,3)(3,3)

\drawline[AHnb=0](4.5,1)(7.5,1)
\drawline[AHnb=0](4.5,2)(7.5,2)
\drawline[AHnb=0](4.5,3)(7.5,3)

\drawline[AHnb=0](9,1)(11,1)
\drawline[AHnb=0](9,2)(11,2)
\drawline[AHnb=0](9,3)(11,3)


\drawline[AHnb=0](0,1)(3,4)
\drawline[AHnb=0](1,1)(3,3)
\drawline[AHnb=0](2,1)(3,2)
\drawline[AHnb=0](0,2)(2,4)
\drawline[AHnb=0](0,3)(1,4)

\drawline[AHnb=0](4.5,1)(7.5,4)
\drawline[AHnb=0](4.5,2)(6.5,4)
\drawline[AHnb=0](5.5,1)(7.5,3)
\drawline[AHnb=0](6.5,1)(7.5,2)
\drawline[AHnb=0](4.5,3)(5.5,4)

\put(3.7,1.1){\scriptsize $\ldots$} 

\put(3.7,1.1){\scriptsize $\ldots$} 
\put(3.7,1.1){\scriptsize $\ldots$}

\put(3.7,1.9){\scriptsize $\ldots$} 
\put(3.7,1.9){\scriptsize $\ldots$} 
\put(3.7,1.9){\scriptsize $\ldots$} 

\put(3.7,0){\scriptsize $\ldots$} 
\put(3.7,0){\scriptsize $\ldots$} 
\put(3.7,0){\scriptsize $\ldots$}

\put(3.7,2.9){\scriptsize $\ldots$} 
\put(3.7,2.9){\scriptsize $\ldots$} 
\put(3.7,2.9){\scriptsize $\ldots$}

\put(8.1,1.1){\scriptsize $\ldots$} 
\put(8.1,1.1){\scriptsize $\ldots$} 
\put(8.1,1.1){\scriptsize $\ldots$} 

\put(8.1,1.9){\scriptsize $\ldots$} 
\put(8.1,1.9){\scriptsize $\ldots$} 
\put(8.1,1.9){\scriptsize $\ldots$} 

\put(8.1,0){\scriptsize $\ldots$} 
\put(8.1,0){\scriptsize $\ldots$} 
\put(8.1,0){\scriptsize $\ldots$} 

\put(8.1,2.9){\scriptsize $\ldots$} 
\put(8.1,2.9){\scriptsize $\ldots$} 
\put(8.1,2.9){\scriptsize $\ldots$} 

\put(-.2,-.4){\scriptsize {\tiny $x_1$}}
\put(.9,-.4){\scriptsize {\tiny $x_2$}}
\put(1.88,-.4){\scriptsize {\tiny $x_3$}}
\put(2.88,-.4){\scriptsize {\tiny $x_4$}}
\put(4.3,-.4){\scriptsize {\tiny $x_k$}}
\put(5.3,-.4){\scriptsize {\tiny $x_{k+1}$}}
\put(6.3,-.4){\scriptsize {\tiny $x_{k+2}$}}
\put(7.4,-.4){\scriptsize {\tiny $x_{k+3}$}}
\put(9,-.4){\scriptsize {\tiny $x_{i-1}$}}
\put(10,-.4){\scriptsize {\tiny $x_{i}$}}
\put(11,-.4){\scriptsize {\tiny $x_1$}}

\put(-.4,.7){\scriptsize {\tiny $y_1$}}
\put(1.12,.7){\scriptsize {\tiny $y_2$}}
\put(2.12,.7){\scriptsize {\tiny $y_3$}}
\put(3.12,.7){\scriptsize {\tiny $y_4$}}
\put(4.6,.7){\scriptsize {\tiny $y_k$}}
\put(5.6,.7){\scriptsize {\tiny $y_{k+1}$}}
\put(6.6,.7){\scriptsize {\tiny $y_{k+2}$}}
\put(7.6,.7){\scriptsize {\tiny $y_{k+3}$}}
\put(9.1,.7){\scriptsize {\tiny $y_{i-1}$}}
\put(10.1,.7){\scriptsize {\tiny $y_{i}$}}
\put(11.1,.7){\scriptsize {\tiny $y_{1}$}}

\put(-.4,1.7){\scriptsize {\tiny $z_1$}}
\put(1.12,1.7){\scriptsize {\tiny $z_2$}}
\put(2.12,1.7){\scriptsize {\tiny $z_3$}}
\put(3.12,1.7){\scriptsize {\tiny $z_4$}}
\put(4.6,1.7){\scriptsize {\tiny $z_k$}}
\put(5.6,1.74){\scriptsize {\tiny $z_{k+1}$}}
\put(6.6,1.74){\scriptsize {\tiny $z_{k+2}$}}
\put(7.6,1.76){\scriptsize {\tiny $z_{k+3}$}}
\put(9.1,1.77){\scriptsize {\tiny $z_{i-1}$}}
\put(10.1,1.77){\scriptsize {\tiny $z_{i}$}}
\put(11.1,1.77){\scriptsize {\tiny $z_{1}$}}

\put(-.5,2.7){\scriptsize {\tiny $w_1$}}
\put(1.12,2.7){\scriptsize {\tiny $w_2$}}
\put(2.12,2.7){\scriptsize {\tiny $w_3$}}
\put(3.12,2.7){\scriptsize {\tiny $w_4$}}
\put(4.6,2.8){\scriptsize {\tiny $w_k$}}
\put(5.5,2.77){\scriptsize {\tiny $w_{k+1}$}}
\put(6.5,2.77){\scriptsize {\tiny $w_{k+2}$}}
\put(7.6,2.7){\scriptsize {\tiny $w_{k+3}$}}
\put(9,2.77){\scriptsize {\tiny $w_{i-1}$}}
\put(10.1,2.77){\scriptsize {\tiny $w_{i}$}}
\put(11.1,2.77){\scriptsize {\tiny $w_{1}$}}

\put(-.2,4.15){\scriptsize {\tiny $x_{k+1}$}}
\put(.85,4.15){\scriptsize {\tiny $x_{k+2}$}}
\put(1.85,4.15){\scriptsize {\tiny $x_{k+3}$}}
\put(2.85,4.15){\scriptsize {\tiny $x_{k+4}$}}
\put(4.5,4.16){\scriptsize {\tiny $x_i$}}
\put(5.5,4.16){\scriptsize {\tiny $x_1$}}
\put(6.5,4.16){\scriptsize {\tiny $x_2$}}
\put(7.5,4.1){\scriptsize {\tiny $x_3$}}
\put(8.6,4.16){\scriptsize {\tiny $x_{k-1}$}}
\put(9.8,4.16){\scriptsize {\tiny $x_{k}$}}
\put(10.9,4.16){\scriptsize {\tiny $x_{k+1}$}}


\drawline[AHnb=0](9,1)(11,3)
\drawline[AHnb=0](10,1)(11,2)
\drawline[AHnb=0](9,2)(11,4)
\drawline[AHnb=0](9,3)(10,4)

\put(2.2,-1.4){\scriptsize {\tiny {\bf Figure 5.2.2:} $M(i,4,k)$ representation}} 

\end{picture}

\vspace{2.85cm}

Following Remark \ref{r1} and Lemma \ref{l4.1.3}, we see that a maximal path $P(v_1, v_2, \ldots, v_i)$ of the type $B_\alpha$, for $\alpha \in \{1,2,3\}$ provides a non-contractible cycle $Q=C_i(v_1, v_2, \ldots, v_i)$ of the same type $B_\alpha$. Note that the cycles of type $B_{2}$ and $B_{3}$ are same type cycles as they are mirror image of each other. By Equation (1) of Section \ref{s4.1}, there is a cycle of type $B_4$ in $M_2$. Let $u \in V(M_2)$ and $Q_{\alpha}$ be cycles of type $B_{\alpha}$ through $u$. As in Section \ref{s4.1}, we define an $M(i,j,k)$ representation of $M_2$ for suitable $i, j, k$. For this, we first cut $M_2$ along the cycle $Q_{1}$ and then cut it along the cycle $Q_{3}$. Without loss of generality, let quadrangular faces are incident with the horizontal base cycle $Q_{1}$, see for example $M(i,4,k)$ in Figure 5.2.2.

\begin{lem}\label{l4.2.1} The DSEM $M_2$ of type $T_2$ admits an $M(i,j,k)$ representation iff the following holds: $(i)$ $i \geq 3$ and $j=4m$, where $m \in \mathbb{N}$, $(ii)$ $ ij \geq 12 $, $(iii)$ $ 0 \leq k \leq i-1 $.
\end{lem}

\noindent{\bf Proof.} We consider a map of type $T_2$  in place of $T_1$ and different value of $i, j, k$. We proceed similarly as in the case of Lemma \ref{l4.1.7}. This gives all the cases.  \hfill$\Box$

Let $M_{t}$, $ t \in \{1,2\}$, be DSEMs of type $T_{2}$  on $n_t$ number of vertices and $n_1 = n_2$. Let $M_t(i_{t}, j_{t}, k_{t})$ be a representation of $M_{t}$. Let $Q^t_{\alpha}$ be cycles of type $B_{\alpha}$ and $l^t_{\alpha}$ = length of the cycle of type $B_{\alpha}$, $\alpha \in \{1,2,3,4\}$, in $M_{t}(i_t,j_t,k_t)$. We say $M_t(i_{t}, j_{t}, k_{t})$ has cycle-type $(l^t_1, l^t_2, l^t_3, l^t_4)$ if $l^t_2 \leq l^t_3$ or $(l^t_1, l^t_3, l^t_2, l^t_4)$  if $l^t_3 < l^t_2$. Now, we show the following.

\begin{lem}\label{l4.2.2} 
	The DSEMs $M_{1} \cong  M_{2}$ iff they have same cycle-type.
\end{lem}

\noindent{\bf Proof.} Suppose $M_{1}$ and  $M_{2}$ be two DSEMs of the type $[3^6:3^3.4^2]_2$ with same number of vertices such that they have same cycle-type. Then $l^1_{1} = l^2_{1}$, $\{l^1_{2}, l_{1,3}\} = \{l^2_{2}, l^2_{3}\}$ and $l^1_{4} = l^2_{4}$. To show $M_1 \cong M_2$, it is equivalent to show that $M_1(i_1,j_1,k_1) \cong M_2(i_2,j_2,k_2)$.

If $(i_{1}, j_{1}, k_{1})=(i_{2}, j_{2}, k_{2})$, then we consider horizontal cycles of type $B_1$ in $M_t(i_{t}, j_{t}, k_{t})$ and proceed similarly as in Lemma \ref{l4.1.8}. This gives that $M_1(i_1,j_1,k_1) \cong M_2(i_2,j_2,k_2)$. Similarly, if $(i_{1}, j_{1}, k_{1}) \neq (i_{2}, j_{2}, k_{2})$, then we get $M_1(i_1,j_1,k_1) \cong M_2(i_2,j_2,k_2)$.

\smallskip
Suppose $M_{1} \cong M_{2}$. Proceeding similarly as in  the converse part of
Lemma \ref{l4.1.8}, we get  $l^1_{1}=l^2_{1}$, $\{l^1_{2}, l^1_{3}\} = \{l^2_{2}, l^2_{3}\}$ and $l^1_{4} = l^2_{4} $. Hence, $M_{1}$ and $M_{2}$ have same cycle-type.  \hfill$\Box$

\smallskip


Now doing the computation for the first four admissible values of $|V(M_2)|$, we get Table \ref{table:2}. For $|V(M_2)| = 12$, we illustrate the computation as follows. 

\begin{eg}\normalfont Let $M_2$ be a DSEM of type $T_2$ with $12$ vertices on the torus. By Lemma \ref{l4.2.1}, $M$ has three $M(i,j,k)$ representation, namely, $M(3, 4, 0), M(3, 4, 1)$ and $M(3, 4, 2)$, see Figures 5.2.3, 5.2.4, and 5.2.5 respectively. In $M(3, 4, 0)$, $Q^1_1 = C_3(v_{1}, v_{2}, v_{3})$ is a $B_{1}$ type cycle, $Q^1_2 =C_4(v_{1},v_{4},u_{2},u_{6})$ and $Q^1_3 = C_4(v_{1},v_{4},u_{1},u_{4})$ are two $B_{2}$ type cycles and $Q^1_4 = C_4(v_{1},v_{4},u_{1},u_{4})$ is a $B_{4}$ type cycle. In $M( 3,4,1)$, $Q^2_1 = C_3(w_{1},w_{2},w_{3})$ is a $B_{1}$ type cycle, $Q^2_2 = C_{12}(w_{1},w_{4},u_{2},u_{6},w_{2},w_{5}, \linebreak u_{3},u_{4},w_{3},w_{6},u_{1}, u_{5})$ and $Q^2_3 =C_{12}(w_{1},w_{4},u_{1},u_{4},w_{2},w_{5},u_{2},u_{5},  w_{3},w_{6},u_{3},u_{6})$ are two $B_{2}$ type cycles and $Q^2_4 = C_5(w_{2},w_{5},u_{2},u_{5},w_{3})$ is a $B_{4}$ type cycles. In $M( 3,4,2)$, $Q^3_1 = C_{3}(x_{1},x_{2},x_{3})$ is a $B_{1}$ type cycle, $Q^3_2 = C_{12}(x_{1},x_{4},u_{2},u_{6},x_{3},x_{6},u_{1},u_{5},x_{2}, x_{5}, u_{3},u_{4})$ and $Q^3_3 = C_{12}(x_{1},x_{4},u_{1},u_{4},x_{3},x_{6},u_{3}, \linebreak u_{6},x_{2},x_{5},u_{2},u_{5})$ are $B_{2}$ type cycles and $Q^3_4 = C_5(x_{3},x_{6},u_{1},u_{5},x_{2})$ is a $B_{4}$ type cycle.
	
In $M(i,j,k)$, observe that $B_{1}$ type cycles have the same length and $B_{2}$ type cycles have at most two different lengths. Since length$(Q^1_4)$ $\neq$ length$(Q^r_{4})$ for $r \in \{2,3\}$, $M(3, 4, 0) \ncong  M(3, 4, 1),\linebreak M(3, 4, 2)$. Observe that, length $(Q^2_1)$ = length$(Q^3_1)$, \{length$(Q^2_2)$, length$(Q^2_3)$\} = \{length$(Q^3_2)$, length$ (Q^3_3)$\} and length$(Q^2_4)$ = length$(Q^3_4)$. Now cutting $M(3,4,1)$ along the path $P(w_{2},w_{5},u_{3}, \linebreak u_{4},w_{3})$ and identifying along the path $P(w_{1},w_{4},u_{1},u_{4},w_{2})$, leads to Figure 5.2.6, i.e., $M(3,4,2)$. By Lemma \ref{l4.2.2}, $M(3,4,1) \cong M(3,4,2)$. Therefore, there are two DSEMs of type $T_2$ with $12$ vertices on the torus upto isomorphism.
	
\end{eg}

\vspace{-.45cm}
\begin{picture}(0,0)(-10,24)
\setlength{\unitlength}{6mm}

\drawpolygon(0,0)(3,0)(3,4)(0,4)


\drawline[AHnb=0](1,0)(1,4)
\drawline[AHnb=0](2,0)(2,4)


\drawline[AHnb=0](0,1)(3,1)
\drawline[AHnb=0](0,2)(3,2)
\drawline[AHnb=0](0,3)(3,3)


\drawline[AHnb=0](0,3)(1,4)
\drawline[AHnb=0](0,2)(2,4)
\drawline[AHnb=0](0,1)(3,4)
\drawline[AHnb=0](1,1)(3,3)
\drawline[AHnb=0](2,1)(3,2)

\put(-.45,-.4){\scriptsize {\tiny $v_1$}}
\put(.8,-.4){\scriptsize {\tiny $v_2$}}
\put(1.8,-.4){\scriptsize {\tiny $v_3$}}
\put(2.8,-.4){\scriptsize {\tiny $v_1$}}

\put(-.6,.7){\scriptsize {\tiny $v_4$}}
\put(1.05,.7){\scriptsize {\tiny $v_5$}}
\put(2.05,.7){\scriptsize {\tiny $v_6$}}
\put(3.05,.7){\scriptsize {\tiny $v_4$}}

\put(-.6,1.7){\scriptsize {\tiny $u_1$}}
\put(1.05,1.7){\scriptsize {\tiny $u_2$}}
\put(2.05,1.7){\scriptsize {\tiny $u_3$}}
\put(3.05,1.7){\scriptsize {\tiny $u_1$}}

\put(-0.6,2.7){\scriptsize {\tiny $u_4$}}
\put(1.1,2.7){\scriptsize {\tiny $u_5$}}
\put(2.1,2.7){\scriptsize {\tiny $u_6$}}
\put(3.1,2.7){\scriptsize {\tiny $u_4$}}

\put(-.6,3.7){\scriptsize {\tiny $v_1$}}
\put(1.1,3.7){\scriptsize {\tiny $v_2$}}
\put(2.1,3.7){\scriptsize {\tiny $v_3$}}
\put(3.1,3.7){\scriptsize {\tiny $v_1$}}

\put(-.7,-1.2){\scriptsize {\tiny Figure 5.2.3 : $M(3,4,0)$}}

\end{picture}

\begin{picture}(0,0)(-16,19)
\setlength{\unitlength}{6mm}

\drawpolygon(5,0)(8,0)(8,4)(5,4)


\drawline[AHnb=0](5,1)(8,1)
\drawline[AHnb=0](5,2)(8,2)
\drawline[AHnb=0](5,3)(8,3)


\drawline[AHnb=0](6,0)(6,4)
\drawline[AHnb=0](7,0)(7,4)


\drawline[AHnb=0](5,3)(6,4)
\drawline[AHnb=0](5,2)(7,4)
\drawline[AHnb=0](5,1)(8,4)
\drawline[AHnb=0](6,1)(8,3)
\drawline[AHnb=0](7,1)(8,2)

\put(4.8,-.4){\scriptsize {\tiny $w_1$}}
\put(5.8,-.4){\scriptsize {\tiny $w_2$}}
\put(6.8,-.4){\scriptsize {\tiny $w_3$}}
\put(7.8,-.4){\scriptsize {\tiny $w_1$}}

\put(4.4,.7){\scriptsize {\tiny $w_4$}}
\put(6.1,.7){\scriptsize {\tiny $w_5$}}
\put(7.1,.7){\scriptsize {\tiny $w_6$}}
\put(8.1,.7){\scriptsize {\tiny $w_4$}}

\put(4.4,1.7){\scriptsize {\tiny $u_1$}}
\put(6.1,1.7){\scriptsize {\tiny $u_2$}}
\put(7.1,1.7){\scriptsize {\tiny $u_3$}}
\put(8.1,1.7){\scriptsize {\tiny $u_1$}}

\put(4.4,2.7){\scriptsize {\tiny $u_4$}}
\put(6.1,2.7){\scriptsize {\tiny $u_5$}}
\put(7.1,2.7){\scriptsize {\tiny $u_6$}}
\put(8.1,2.7){\scriptsize {\tiny $u_4$}}

\put(4.4,3.7){\scriptsize {\tiny $w_2$}}
\put(6.1,3.7){\scriptsize {\tiny $w_3$}}
\put(7.1,3.7){\scriptsize {\tiny $w_1$}}
\put(8.1,3.7){\scriptsize {\tiny $w_2$}}

\put(4.3,-1.2){\scriptsize {\tiny Figure 5.2.4 : $M(3,4,1)$}}

\end{picture}

\begin{picture}(0,0)(-20,14)
\setlength{\unitlength}{6mm}

\drawpolygon(10,0)(13,0)(13,4)(10,4)

\drawline[AHnb=0](11,0)(11,4)
\drawline[AHnb=0](12,0)(12,4)


\drawline[AHnb=0](10,1)(13,1)
\drawline[AHnb=0](10,2)(13,2)
\drawline[AHnb=0](10,3)(13,3)


\drawline[AHnb=0](10,3)(11,4)
\drawline[AHnb=0](10,2)(12,4)
\drawline[AHnb=0](10,1)(13,4)
\drawline[AHnb=0](11,1)(13,3)
\drawline[AHnb=0](12,1)(13,2)

\put(9.8,-.4){\scriptsize {\tiny $x_1$}}
\put(10.8,-.4){\scriptsize {\tiny $x_2$}}
\put(11.8,-.4){\scriptsize {\tiny $x_3$}}
\put(12.8,-.4){\scriptsize {\tiny $x_1$}}

\put(9.4,.7){\scriptsize {\tiny $x_4$}}
\put(11.1,.7){\scriptsize {\tiny $x_5$}}
\put(12.1,.7){\scriptsize {\tiny $x_6$}}
\put(13.1,.7){\scriptsize {\tiny $x_4$}}

\put(9.4,1.7){\scriptsize {\tiny $u_1$}}
\put(11.1,1.7){\scriptsize {\tiny $u_2$}}
\put(12.1,1.7){\scriptsize {\tiny $u_3$}}
\put(13.1,1.7){\scriptsize {\tiny $u_1$}}

\put(9.4,2.7){\scriptsize {\tiny $u_4$}}
\put(11.1,2.7){\scriptsize {\tiny $u_5$}}
\put(12.1,2.7){\scriptsize {\tiny $u_6$}}
\put(13.1,2.7){\scriptsize {\tiny $u_4$}}

\put(9.4,3.7){\scriptsize {\tiny $x_3$}}
\put(11.1,3.7){\scriptsize {\tiny $x_1$}}
\put(12.1,3.7){\scriptsize {\tiny $x_2$}}
\put(13.1,3.7){\scriptsize {\tiny $x_3$}}

\put(9.3,-1.2){\scriptsize {\tiny Figure 5.2.5 : $M(3,4,2)$}}

\end{picture}

\begin{picture}(0,0)(-24,9)
\setlength{\unitlength}{6mm}

\drawpolygon(15,0)(18,0)(18,4)(15,4)

\drawline[AHnb=0](16,0)(16,4)
\drawline[AHnb=0](17,0)(17,4)


\drawline[AHnb=0](15,1)(18,1)
\drawline[AHnb=0](15,2)(18,2)
\drawline[AHnb=0](15,3)(18,3)


\drawline[AHnb=0](15,3)(16,4)
\drawline[AHnb=0](15,2)(17,4)
\drawline[AHnb=0](15,1)(18,4)
\drawline[AHnb=0](16,1)(18,3)
\drawline[AHnb=0](17,1)(18,2)

\put(14.8,-.4){\scriptsize {\tiny $w_2$}}
\put(16,-.4){\scriptsize {\tiny $w_1$}}
\put(17,-.4){\scriptsize {\tiny $w_3$}}
\put(18,-.4){\scriptsize {\tiny $w_2$}}

\put(14.4,.7){\scriptsize {\tiny $w_5$}}
\put(16.1,.7){\scriptsize {\tiny $w_4$}}
\put(17.1,.7){\scriptsize {\tiny $w_6$}}
\put(18.1,.7){\scriptsize {\tiny $w_5$}}

\put(14.4,1.7){\scriptsize {\tiny $u_3$}}
\put(16.1,1.7){\scriptsize {\tiny $u_2$}}
\put(17.1,1.7){\scriptsize {\tiny $u_1$}}
\put(18.1,1.7){\scriptsize {\tiny $u_3$}}

\put(14.4,2.7){\scriptsize {\tiny $u_4$}}
\put(16.1,2.7){\scriptsize {\tiny $u_6$}}
\put(17.1,2.7){\scriptsize {\tiny $u_5$}}
\put(18.1,2.7){\scriptsize {\tiny $u_4$}}

\put(14.4,3.7){\scriptsize {\tiny $w_3$}}
\put(16.1,3.7){\scriptsize {\tiny $w_2$}}
\put(17.1,3.7){\scriptsize {\tiny $w_1$}}
\put(18.1,3.7){\scriptsize {\tiny $w_3$}}

\put(14.4,-1.2){\scriptsize {\tiny Figure 5.2.6 : $M(3,4,2)$}}

\end{picture}

\vspace{1.75cm}
\begin{center}
	
	\noindent \textbf{Table \ref{table:2}} : DSEMs of type $T_2$ on the torus for $|V(M)| \leq 24$
	
	\renewcommand{\arraystretch}{1.1}
	\begin{tabular}{|p{1.2cm}|p{4cm}|p{3.8cm}|p{3cm}|}
		
		\hline
		$|V(M)|$ & Isomorphic classes & Length of cycles & No. of maps\\
		\hline
		12 &  $M(3,4,0)$ & $( 3, \{ 4, 4 \}, 4 )$ & 2\\
		\cline{2-3}
		& $M(3,4,1)$, $M(3,4,2)$ & $( 3, \{ 12, 12 \}, 5 )$ & \\	
		\cline{2-3}
		\hline

		16 &  $M(4,4,0)$, $M(4,4,1)$ & $( 4, \{ 4, 16 \}, 4 )$& 2\\
		\cline{2-3}
		& $M(4,4,2)$, $M(4,4,3)$ & $( 4, \{ 8, 16 \}, 6 )$ & \\
		\cline{2-3}
		
		\hline
		20 & $M(5,4,0)$, $M(5,4,2)$ & $( 5, \{ 4, 20 \}, 4 )$ &3\\
		\cline{2-3}
		& $M(5,4,3)$, $M(5,4,4)$  &$( 5, \{ 20, 20 \}, 7 )$ &\\
		\cline{2-3}
		& $M(5,4,1)$ & $( 5, \{ 20, 20 \}, 5 )$  &\\
		\cline{2-3}
		\hline
		24 & $M(6,4,0)$, $M(6,4,3)$ & $( 6, \{ 4, 8 \}, 4 )$ & 5\\
		\cline{2-3}
		&$M(6,4,1)$, $M(6,4,2)$ &$( 6, \{ 12, 24 \}, 5 )$&\\
		\cline{2-3}
		& $M(6,4,4)$, $M(6,4,5)$ &$( 6, \{ 12, 24 \}, 8)$&\\
		\cline{2-3}
		& $M(3,8,0)$ &$( 3, \{ 8, 8 \}, 8 )$ & \\
		\cline{2-3}
		& $M(3,8,1)$, $M(3,8,2)$ & $( 3, \{ 24, 24 \}, 9 )$ &  \\
		\hline
		
	\end{tabular}
	\label{table:2}
\end{center}

\subsection{DSEMs of type $T_3 = [f_1^{({f_1}^4.{f_2}^3)}:f_2^{({f_1}^3.{f_2}.{f_1}^3.{f_2} )}]$, where $(f_1, f_2)=((3^3.4^2),(4^4))$} \label{s4.3}
Let $M_3$ be a DSEM of type $T_3$ with the vertex set $V(M_3)$. Then, for the existence of $M_3$, we get $|V_{(3^3.4^2)}| = 2|V_{(4^4)}|$.  Now we define following three types of paths in $M_3$ as follows.



\begin{deff} \normalfont
	
	A path $P_{1} = P( \ldots, y_{i-1},y_{i},y_{i+1}, \ldots)$ in $M_3$ is of type $X_{1}$ if either of the conditions follows: 
	$(i)$ every vertex of $P_1$ has the face-sequence $(4^4)$ or
	$(ii)$ each vertex of $P_1$ has the face-sequence $(3^3.4^2)$ such that all the triangles (squares) that incident on its inner vertices lie on the one side  of the path $P_1$. For example, see the paths indicated by black color in Figure 5.3.1.
	
\end{deff}

\begin{deff} \normalfont A path $P_{2} = P(\ldots, z_{i-1},z_{i},z_{i+1}, \ldots)$ in $M_3$, such that $z_{i-1}, z_{i}, z_{i+1}$ are inner vertices of $P_{2}$ or an extended path of $P_{2}$, is of type $X_{2}$ (shown by green colored  paths in Figure 5.3.1), if either of the following three conditions follows for each vertex of the path.

	\begin{enumerate}

		\item  If ${\rm lk}(z_{i})=C_7(\boldsymbol{m},z_{i-1},\boldsymbol{n},o,z_{i+1},p,q)$ and  ${\rm lk}(z_{i-1})= C_8(\boldsymbol{q},m,\boldsymbol{r},z_{i-2},\boldsymbol{s},n,\boldsymbol{o},z_i)$,  then ${\rm lk}(z_{i+1}) \linebreak =C_7(\boldsymbol{t},z_{i+2},\boldsymbol{u},p,z_{i},o,v)$.
		
		\item If ${\rm lk}(z_{i})=C_7(\boldsymbol{m},z_{i+1},\boldsymbol{n},o,z_{i-1},p,q)$ and  ${\rm lk}(z_{i-1})=C_7(r,\boldsymbol{s},z_{i-2},\boldsymbol{t},p,z_{i},o)$, then ${\rm lk}(z_{i+1})=C_8(\boldsymbol{o},z_{i},\boldsymbol{q},m,\boldsymbol{u},z_{i+2},\boldsymbol{v},n)$. 
		
		\item  If ${\rm lk}(z_{i})= C_8(\boldsymbol{q},r,\boldsymbol{m},z_{i-1},\boldsymbol{n},o,\boldsymbol{p},z_{i+1})$ and ${\rm lk}(z_{i-1})=C_7(\boldsymbol{o},z_{i},\boldsymbol{r},m,z_{i-2},t,n)$, then ${\rm lk}(z_{i+1})=C_7(\boldsymbol{r},z_{i},\boldsymbol{o},p,z_{i+2},s,q)$.
		
	\end{enumerate}
	
\end{deff}

\begin{deff} \normalfont A path $P_{3} = P(\ldots, w_{i-1},w_{i},w_{i+1}, \ldots)$ in $M_3$, such that $w_{i-1}, w_{i}, w_{i+1}$ are inner vertices of $P_{3}$ or an extended path of $P_{3}$, is of type $X_{3}$ (shown by red colored  paths in Figure 5.3.1), if either of the following three conditions follows for each vertex of the path.

	\begin{enumerate}

		\item  If ${\rm lk}(w_{i})=C_7(\boldsymbol{m},w_{i-1},\boldsymbol{n},o,p,w_{i+1},q)$ and  ${\rm lk}(w_{i-1})= C_8(\boldsymbol{q},m,\boldsymbol{r},w_{i-2},\boldsymbol{s},n,\boldsymbol{o},w_i)$,  then ${\rm lk}(w_{i+1}) =C_7(\boldsymbol{t},w_{i+2},\boldsymbol{u},v,q,w_{i},p)$.
		
		\item If ${\rm lk}(w_{i})=C_7(\boldsymbol{m},w_{i+1},\boldsymbol{n},o,p,w_{i-1},q)$ and  ${\rm lk}(w_{i-1})=C_7(\boldsymbol{r},w_{i-2},\boldsymbol{s},t,q,w_{i},p)$, then ${\rm lk}(w_{i+1})=C_8(\boldsymbol{o},w_{i},\boldsymbol{q},m,\boldsymbol{u},w_{i+2},\boldsymbol{v},n)$. 
		
		\item  If ${\rm lk}(w_{i})= C_8(\boldsymbol{q},r,\boldsymbol{m},w_{i-1},\boldsymbol{n},o,\boldsymbol{p},w_{i+1})$ and ${\rm lk}(w_{i-1})=C_7(\boldsymbol{o},w_{i},\boldsymbol{r},m,t,w_{i-2},n)$, then ${\rm lk}(w_{i+1})=C_7(\boldsymbol{r},w_{i},\boldsymbol{o},p,u,w_{i+2},q)$.
		
	\end{enumerate}
	
\end{deff}

\vspace{-.5cm}

\begin{picture}(0,0)(-15,25)
\setlength{\unitlength}{6mm}

\drawpolygon(0,0)(6,0)(6,2)(0,2)


\drawline[AHnb=0](1,0)(1,2)
\drawline[AHnb=0](2,0)(2,2)
\drawline[AHnb=0](3,0)(3,2)
\drawline[AHnb=0](4,0)(4,2)
\drawline[AHnb=0](5,0)(5,2)

\drawline[AHnb=0](.5,3)(.5,3.5)
\drawline[AHnb=0](1.5,3)(1.5,3.5)

\drawline[AHnb=0](2.5,3)(2.5,3.5)
\drawline[AHnb=0](3.5,3)(3.5,3.5)

\drawline[AHnb=0](4.5,3)(4.5,3.5)
\drawline[AHnb=0](5.5,3)(5.5,3.5)


\drawline[AHnb=0](6,0)(6.5,0)
\drawline[AHnb=0](0,0)(-.5,0)
\drawline[AHnb=0](6,3)(6.5,3)
\drawline[AHnb=0](0,3)(-.5,3)

\drawline[AHnb=0](-.5,1)(6.5,1)
\drawline[AHnb=0](-.5,2)(6.5,2)

\drawline[AHnb=0](-.5,3)(6.5,3)


\drawline[AHnb=0](0,2)(.5,3)
\drawline[AHnb=0](1,2)(.5,3)

\drawline[AHnb=0](1,2)(1.5,3)
\drawline[AHnb=0](2,2)(1.5,3)

\drawline[AHnb=0](2,2)(2.5,3)
\drawline[AHnb=0](3,2)(2.5,3)

\drawline[AHnb=0](3,2)(3.5,3)
\drawline[AHnb=0](4,2)(3.5,3)

\drawline[AHnb=0](4,2)(4.5,3)
\drawline[AHnb=0](5,2)(4.5,3)

\drawline[AHnb=0](5,2)(5.5,3)
\drawline[AHnb=0](6,2)(5.5,3)


\drawline[AHnb=0](0,0)(.25,-.5)
\drawline[AHnb=0](1,0)(1.25,-.5)
\drawline[AHnb=0](2,0)(2.25,-.5)
\drawline[AHnb=0](3,0)(3.25,-.5)
\drawline[AHnb=0](4,0)(4.25,-.5)
\drawline[AHnb=0](5,0)(5.25,-.5)
\drawline[AHnb=0](6,0)(6.25,-.5)

\drawline[AHnb=0](0,0)(-.25,-.5)
\drawline[AHnb=0](1,0)(.75,-.5)
\drawline[AHnb=0](2,0)(1.75,-.5)
\drawline[AHnb=0](3,0)(2.75,-.5)
\drawline[AHnb=0](4,0)(3.75,-.5)
\drawline[AHnb=0](5,0)(4.75,-.5)
\drawline[AHnb=0](6,0)(5.75,-.5)

\drawline[AHnb=0](0,2)(-.25,2.5)
\drawline[AHnb=0](6,2)(6.25,2.5)


\drawpolygon[fillcolor=green](-.25,-.5)(0,0)(0,2)(.5,3)(.5,3.5)(.4,3.5)(.4,3)(-.1,2)(-.1,0)(-.35,-.5)

\drawpolygon[fillcolor=green](.75,-.5)(1,0)(1,2)(1.5,3)(1.5,3.5)(1.4,3.5)(1.4,3)(.9,2)(.9,0)(.65,-.5)

\drawpolygon[fillcolor=green](1.75,-.5)(2,0)(2,2)(2.5,3)(2.5,3.5)(2.4,3.5)(2.4,3)(1.9,2)(1.9,0)(1.65,-.5)

\drawpolygon[fillcolor=green](2.75,-.5)(3,0)(3,2)(3.5,3)(3.5,3.5)(3.4,3.5)(3.4,3)(2.9,2)(2.9,0)(2.65,-.5)

\drawpolygon[fillcolor=green](3.75,-.5)(4,0)(4,2)(4.5,3)(4.5,3.5)(4.4,3.5)(4.4,3)(3.9,2)(3.9,0)(3.65,-.5)

\drawpolygon[fillcolor=green](4.75,-.5)(5,0)(5,2)(5.5,3)(5.5,3.5)(5.4,3.5)(5.4,3)(4.9,2)(4.9,0)(4.65,-.5)

\drawpolygon[fillcolor=green](5.75,-.5)(6,0)(6,2)(6.25,2.5)(6.15,2.5)(5.9,2)(5.9,0)(5.65,-.5)

\drawpolygon[fillcolor=red](.25,-.5)(0,0)(0,2)(-.25,2.5)(-.15,2.5)(.1,2)(.1,0)(.35,-.5)

\drawpolygon[fillcolor=red](1.25,-.5)(1,0)(1,2)(.5,3)(.5,3.5)(.6,3.5)(.6,3)(1.1,2)(1.1,0)(1.35,-.5)

\drawpolygon[fillcolor=red](2.25,-.5)(2,0)(2,2)(1.5,3)(1.5,3.5)(1.6,3.5)(1.6,3)(2.1,2)(2.1,0)(2.35,-.5)

\drawpolygon[fillcolor=red](3.25,-.5)(3,0)(3,2)(2.5,3)(2.5,3.5)(2.6,3.5)(2.6,3)(3.1,2)(3.1,0)(3.35,-.5)

\drawpolygon[fillcolor=red](4.25,-.5)(4,0)(4,2)(3.5,3)(3.5,3.5)(3.6,3.5)(3.6,3)(4.1,2)(4.1,0)(4.35,-.5)

\drawpolygon[fillcolor=red](5.25,-.5)(5,0)(5,2)(4.5,3)(4.5,3.5)(4.6,3.5)(4.6,3)(5.1,2)(5.1,0)(5.35,-.5)

\drawpolygon[fillcolor=red](6.25,-.5)(6,0)(6,2)(5.5,3)(5.5,3.5)(5.6,3.5)(5.6,3)(6.1,2)(6.1,0)(6.35,-.5)

\drawpolygon[fillcolor=black](-.5,0)(6.5,0)(6.5,.1)(-.5,.1)

\drawpolygon[fillcolor=black](-.5,1)(6.5,1)(6.5,1.1)(-.5,1.1)

\drawpolygon[fillcolor=black](-.5,2)(6.5,2)(6.5,2.1)(-.5,2.1)

\drawpolygon[fillcolor=black](-.5,3)(6.5,3)(6.5,3.1)(-.5,3.1)

\put(2.4,3.75) {\scriptsize {\tiny $X_3$}}
\put(3.4,3.75) {\scriptsize {\tiny $X_2$}}
\put(7.1,-.1) {\scriptsize {\tiny $X_1$}}

\put(-.5,-1.5){\scriptsize {\tiny Figure 5.3.1 : Paths of type $X_1, X_2, X_3$}} 
\end{picture}

\begin{picture}(0,0)(-72,30)
\setlength{\unitlength}{6.75mm}

\drawpolygon(0,1)(3,1)(3,4)(0,4)
\drawpolygon(5,1)(8,1)(8,4)(5,4)
\drawpolygon(10,1)(12,1)(12,4)(10,4)


\drawline[AHnb=0](1,1)(1,4)
\drawline[AHnb=0](2,1)(2,4)
\drawline[AHnb=0](3,1)(3,4)
\drawline[AHnb=0](6,1)(6,4)
\drawline[AHnb=0](7,1)(7,4)
\drawline[AHnb=0](11,1)(11,4)


\drawline[AHnb=0](0,1)(3,1)
\drawline[AHnb=0](0,2)(3,2)
\drawline[AHnb=0](0,3)(3,3)

\drawline[AHnb=0](5,1)(8,1)
\drawline[AHnb=0](5,2)(8,2)
\drawline[AHnb=0](5,3)(8,3)

\drawline[AHnb=0](10,1)(12,1)
\drawline[AHnb=0](10,2)(12,2)
\drawline[AHnb=0](10,3)(12,3)


\drawline[AHnb=0](0,3)(1,4)
\drawline[AHnb=0](1,3)(2,4)
\drawline[AHnb=0](2,3)(3,4)

\drawline[AHnb=0](5,3)(6,4)
\drawline[AHnb=0](6,3)(7,4)
\drawline[AHnb=0](7,3)(8,4)

\drawline[AHnb=0](10,3)(11,4)
\drawline[AHnb=0](11,3)(12,4)

\put(-.2,.8){\scriptsize {\tiny $v_1$}}
\put(.7,.75){\scriptsize {\tiny $v_2$}}
\put(1.7,.75){\scriptsize {\tiny $v_3$}}
\put(2.7,.75){\scriptsize {\tiny $v_4$}}
\put(5,.75){\scriptsize {\tiny $v_{k}$}}
\put(5.8,.75){\scriptsize {\tiny $v_{k+1}$}}
\put(6.8,.75){\scriptsize {\tiny $v_{k+2}$}}
\put(8,.75){\scriptsize {\tiny $v_{k+3}$}}

\put(10,.75){\scriptsize {\tiny $v_{i-2}$}}
\put(11,.75){\scriptsize {\tiny $v_{i}$}}
\put(12,.75){\scriptsize {\tiny $v_{1}$}}

\put(-.5,1.75){\scriptsize {\tiny $w_1$}}
\put(.5,1.7){\scriptsize {\tiny $w_2$}}
\put(1.5,1.7){\scriptsize {\tiny $w_3$}}
\put(2.5,1.7){\scriptsize {\tiny $w_4$}}
\put(5.1,1.7){\scriptsize {\tiny $w_{k}$}}
\put(6.1,1.7){\scriptsize {\tiny $w_{k+1}$}}
\put(7.1,1.7){\scriptsize {\tiny $w_{k+2}$}}
\put(8.1,1.7){\scriptsize {\tiny $w_{k+3}$}}
\put(10.1,1.7){\scriptsize {\tiny $w_{i-2}$}}
\put(11.1,1.7){\scriptsize {\tiny $w_{i}$}}
\put(12.1,1.7){\scriptsize {\tiny $w_{1}$}}

\put(-.5,2.7){\scriptsize {\tiny $x_1$}}
\put(.5,2.7){\scriptsize {\tiny $x_2$}}
\put(1.5,2.7){\scriptsize {\tiny $x_3$}}
\put(2.5,2.7){\scriptsize {\tiny $x_4$}}

\put(5.1,2.7){\scriptsize {\tiny $x_{k}$}}
\put(6.1,2.7){\scriptsize {\tiny $x_{k+1}$}}
\put(7.1,2.7){\scriptsize {\tiny $x_{k+2}$}}
\put(8.1,2.7){\scriptsize {\tiny $x_{k+3}$}}

\put(10.1,2.75){\scriptsize {\tiny $x_{i-1}$}}
\put(11.2,2.7){\scriptsize {\tiny $x_{i}$}}
\put(12.1,2.7){\scriptsize {\tiny $x_1$}}

\put(-.5,4.2){\scriptsize {\tiny $v_{k+1}$}}
\put(.5,4.2){\scriptsize {\tiny $v_{k+2}$}}
\put(1.5,4.2){\scriptsize {\tiny $v_{k+3}$}}
\put(2.5,4.2){\scriptsize {\tiny $v_{k+4}$}}
\put(4.8,4.2){\scriptsize {\tiny $v_{i}$}}
\put(5.7,4.2){\scriptsize {\tiny $v_{1}$}}
\put(6.7,4.2){\scriptsize {\tiny $v_{2}$}}
\put(7.8,4.2){\scriptsize {\tiny $v_{3}$}}
\put(9.7,4.2){\scriptsize {\tiny $v_{k-1}$}}
\put(10.8,4.2){\scriptsize {\tiny $v_{k}$}}
\put(11.5,4.2){\scriptsize {\tiny $v_{k+1}$}}

\put(3.8,1){\scriptsize $\ldots$}
\put(3.8,1){\scriptsize $\ldots$}
\put(3.8,2){\scriptsize $\ldots$}
\put(3.8,2){\scriptsize $\ldots$}
\put(3.8,3){\scriptsize $\ldots$}
\put(3.8,3){\scriptsize $\ldots$}
\put(3.8,3.8){\scriptsize $\ldots$}
\put(3.8,3.8){\scriptsize $\ldots$}

\put(9.1,1){\scriptsize $\ldots$}
\put(9.1,1){\scriptsize $\ldots$}
\put(9.1,2){\scriptsize $\ldots$}
\put(9.1,2){\scriptsize $\ldots$}
\put(9.1,3){\scriptsize $\ldots$}
\put(9.1,3){\scriptsize $\ldots$}
\put(9.1,3.8){\scriptsize $\ldots$}
\put(9.1,3.8){\scriptsize $\ldots$}

\put(4.5,0){\scriptsize {\tiny {\bf Figure 5.3.2:} $M(i,3,k)$ }} 

\end{picture}

\vspace{3.25cm}

Following Remark \ref{r1}, and Lemma \ref{l4.1.3}, for a maximal path $P$ of the type $X_\alpha$, $\alpha \in \{1,2,3\}$ there is an edge $e$ in $M_3$ such that  $P\cup e$ is a non-contractible cycle of respective type. The cycles of types $X_{2}$ and $X_{3}$ define same type of cycles as they are mirror image of each other. Following Equation (1) of Section \ref{s4.1}, there is a cycle of type $X_4$ in $M_3$. As in Section \ref{s4.1}, we define an $M(i,j,k)$ representation for $M_3$ for some $i, j, k$. See for example $M(i,3,k)$ in Figure 5.3.2.

\begin{lem}\label{l4.3.1} The DSEM $M_3$ of type $T_3$ admits an $M(i,j,k)$-representation iff the following holds: $(i)$ $i \geq 3$ and $j=3m$, where $m \in \mathbb{N} $, $(ii)$ $ij \geq 9 $, $(iii)$ $ 0 \leq k \leq i-1 $.
\end{lem}
\noindent{\bf Proof.} Let $M_3$ be a DSEM of  type $T_3$ with $n$ vertices. Note that the map exists if $|V(3^3.4^2)|$ = $2|V(4^4)|$. Now proceeding, as in Lemma \ref{l4.1.7}, we get all possible values of $i$, $j$ and $k$ of $M(i,j,k)$. Thus the proof. \hfill $\Box$

\smallskip

For $ t \in \{1,2\}$, let $M_{t}$ be DSEMs of type $T_3$ on $n_t$ number of vertices with $n_1 = n_2$. Let $M_{t}(i_{t}, j_{t}, k_{t})$ be a representation of $M_t$ and $Q^t_{\alpha}$ be cycles of type $X_{\alpha}, \ \alpha \in \{1,2,3,4\}$. If $l^t_{\alpha}$ = length of the cycle of type $X_{\alpha}$ in $M_{t}(i_t,j_t,k_t)$ then we say that $M_{t}(i_t,j_t,k_t)$ has cycle-type $(l^t_1, l^t_2, l^t_3, l^t_4)$ if $l^t_2 \leq l^t_3$ or $(l^t_1, l^t_3, l^t_2, l^t_4)$  if $l^t_3 < l^t_2$. Now, we show the following lemma.

\begin{lem}\label{l4.3.2} The DSEMs $M_{1} \cong  M_{2}$ iff they have same cycle-type.
\end{lem}
\noindent{\bf Proof.} Suppose $M_{1}$ and  $M_{2}$ be two DSEMs of the type $[3^3.4^2:4^4]_1$ with same cardinality such that they have same cycle-type. This means $l^1_{1} = l^2_{1}$, $\{l^1_{2},l^1_{3}\}=\{l^2_{2},l^2_{3}\}$ and $l^1_{4} = l^2_{4} $.

\smallskip

\noindent{Claim.} $M_{1}(i_1,j_1,k_1) \cong M_{2}(i_2,j_2,k_2)$. 

\smallskip

By the definition, $M_{t}(i_{t},j_{t},k_{t})$ has $j_t$ number of $X_{1}$ type disjoint horizontal cycles of length $i_t$, say,  
$Q_{0}=C_{i_1}(w_{0,0},w_{0,1},\ldots,w_{0,i_{1}-1}), Q_{1}=C_{i_1}(w_{1,0},w_{1,1},\ldots,w_{1,i_{1}-1}),\ldots,Q_{j_{1}-1}=C_{i_1}(w_{{j_1}-1,0},w_{{j_1}-1,1}, \linebreak \ldots,w_{j_{1}-1,i_{1}-1})$ in $M_{1}(i_1,j_1,k_1)$  and 
$Q_{0}'=C_{i_2}(z_{0,0}$, $z_{0,1},\ldots,z_{0,i_{2}-1})$, $Q_{1}'= C_{i_2}(z_{1,0},z_{1,1},\ldots, \linebreak z_{1,i_{2}-1}),\ldots, Q_{j_{2}-1}'=C_{i_2}(z_{j_{2}-1,0},z_{j_{2}-1,1},\ldots,z_{j_{2}-1,i_{2}-1})$ in $M_{2}(i_2,j_2,k_2)$. Then, 

\smallskip

\noindent{\bf Case 1:} If $(i_{1}, j_{1}, k_{1})=(i_{2}, j_{2}, k_{2})$ then $i_{1}=i_{2}, j_{1}=j_{2}, k_{1}=k_{2}$. Define an isomorphism  $f :V(M_{1}(i_1,j_1$, $k_1)) \to V(M_{2}(i_2,j_2,k_2))$ such that $f(w_{u,v})=z_{u,v}$ for $ 0 \leq u \leq j_{1}-1 $ and $ 0 \leq v \leq i_{1}-1 $. So by $f$, $M_{1}(i_1,j_1,k_1) \cong M_{2}(i_2,j_2,k_2)$.

\smallskip

\noindent{\bf Case 2:} If $ i_{1} \neq i_{2}$, then it contradicts the fact that $l^1_{1} = l^2_{1}$. Thus $i_{1}=i_{2} $.

\smallskip

\noindent{\bf Case 3:} If $ j_{1} \neq j_{2} $, it implies that $n_{1}=i_{1}j_{1} \neq n_{2}=i_{2}j_{2}$ as $ i_{1}=i_{2}$. A contradiction. So, $j_{1}=j_{2}$.

\smallskip

\noindent{\bf Case 4:} If $k_{1} \neq k_{2} $ then by $l^1_{4} = l^2_{4} $ we see length$(Q^1_4)$ = length$(Q^2_4)$. This means min$\{k_1+j_1,j_1+(i_1-k_1-j_{1}/3)\}$ =  min$\{k_2+j_2,j_2+(i_2-k_2-j_{2}/3)\}$. Since $i_{1}=i_{2}, j_{1}=j_{2}$ and $k_{1} \neq k_{2}$, $k_1+j_1 \neq k_2+j_2$ and $j_1+(i_1-k_1-j_{1}/3) \neq j_2+(i_2-k_2-j_{2}/3)$. This implies $k_1+j_1=i_2-k_2-j_{2}/3+j_{2}=i_1-k_2+2j_{1}/3$ as $i_{1}=i_{2}$ and $j_{1}=j_{2}$. That is, $k_{2}=i_{1}-k_{1}-j_{1}/3$.
Now identify $M_{2}(i_2,j_2,k_2)$ along the vertical boundary $P(z_{0,0},z_{1,0},\ldots , z_{{j_{2}-1},0},z_{0,k_{2}})$ and then cut along the path $P(z_{0,0},z_{1,0},z_{2,0},z_{3,1}, \ldots , z_{{j_{2}-1,j_{2}/3-1}},z_{0,k_{2}+j_{2}/3})$ of type $X_{2}$ through vertex $z_{0,0}$. This gives another representation of $M_2$, say $R$, with a map $f_1: V(M_{2}(i_2,j_2,k_2)) \to V(R)$ such that $f_1({z_{u,v}})=z_{u,(i_{2}-v+ \lfloor t/3 \rfloor)(mod\,i_{2})}$ for $0 \leq u \leq j_{2}-1$ and $0 \leq v \leq i_{2}-1 $. In $R$ the lower and upper horizontal cycles are $Q'=C_{i_2}(z_{0,0},z_{0,i_{2}-1},z_{0,i_{2}-2}, \ldots , z_{0,1})$ and $Q''=C_{i_2}(z_{0,k_{2}+j_{2}/3}, z_{0,k_{2}+j_{2}/3-1}, \ldots, z_{0,k_{2}+j_{2}/3+1})$ respectively. The path $P(z_{0,0},z_{0,i_{2}-1},z_{0,i_{2}-2}, \ldots , z_{0,k_{2}+j_{2}/3})$ in $Q'$ has length $i_{2}-k_{2}-j_{2}/3$. Note that $R$ has $j_{2}$ number of horizontal cycles of length $i_{2}$. So, $R=M_2(i_{2},j_{2},i_{2}-k_{2}-j_{2}/3)$. Note that $i_{2}-k_{2}-j_{2}/3=i_{2}-(i_{1}-k_{1}-j_{1}/3) -j_{2}/3=k_{1}$ as $i_{2}=i_{1},j_{2}=j_{1}$. Thus, $(i_{2},j_{2},i_{2}-k_{2}-j_{2}/3) = (i_{1}, j_{1}, k_{1})$. Therefore, $M_{1}(i_1,j_1,k_1) \cong M_{2}(i_2,j_2,k_2)$. So the claim follows. Hence, $M_{1} \cong M_{2}$.

\smallskip
The converse part follows similarly as the converse part of Lemma \ref{l4.1.8}. Thus the proof. \hfill $\Box$

\smallskip

Now computing the DSEMs for the first four admissible values of $|V(M_3)|$, we get Table \ref{table:3}. We illustrate this computation below for $|V(M_3)| = 9$. 

\begin{eg} \normalfont Let $M_3$ be a DSEM of type $T_3$ with $9$ vertices on the torus. By Lemma \ref{l4.3.1}, $M$ has three $M(i,j,k)$ representation, namely, $M(3, 3, 0), M(3, 3, 1)$ and $M(3, 3, 2)$, see Figures 5.3.3, 5.3.4, and 5.3.5 respectively. In $M(3,3,0)$, $Q^1_1 = C_3(x_{1}, x_{2}, x_{3})$ is a $X_{1}$ type cycle, $Q^1_2 =C_9(x_{1},x_{4},x_{7},x_{2},x_{5},x_{8},x_{3},x_{6},x_{9})$ and $Q^1_3 = C_3(x_{1},x_{4},x_{7})$ are $X_{2}$ type cycles and $Q^1_4 = C_3(x_{1},x_{4}$, $x_{7})$ is a $X_{4}$ type cycle. In $M(3,3,1)$, $Q^2_1 = C_3(y_{1}, y_{2}, y_{3})$ is a $X_{1}$ type cycle, $Q^2_2 = C_{9}(y_{1},y_{4},y_{7},y_{3}$, $y_{6},y_{9},y_{2},y_{5},y_{8})$ and $Q^2_3 = C_{9}(y_{1},y_{4},y_{7},y_{2},y_{5},y_{8},y_{3},y_{6},y_{9})$ are $X_{2}$ type cycles and $Q^2_4 = C_4(y_{2}$, $y_{5},y_{8},y_{1})$ is a $X_{4}$ type cycle. In $M(3,3,2)$, $Q^3_1 = C_{3}(z_{1},z_{2},z_{3})$ is a $X_{1}$ type cycle, $Q^3_2 = C_{3}(z_{1},z_{4},z_{7})$ and $Q^3_3 = C_{9}(z_{1},z_{4},z_{7},z_{3},z_{6},z_{9},z_{2},z_{5},z_{8})$  are $X_{2}$ type cycle and $Q^3_4 = C_3(z_{3},z_{6},z_{9})$ is a $X_{4}$ type cycle. 	
	
	Observe that type $X_{1}$ cycles have the same length and type $X_{2}$ cycles have at most two different lengths. Since length$(Q^2_4)$ $\neq$ length$(Q^r_{4})$ for $r \in \{1,3\}$, $M(3, 3, 1) \ncong  M(3, 3, 0), M(3, 3, 2)$. Observe that, length $(Q^1_1)$ = length$(Q^3_1)$, $\{$length$(Q^1_2)$, length$(Q^1_3)\} = \{$ length$(Q^3_2)$, length$(Q^3_3)\}$ and length$(Q^1_4)$ =length$(Q^3_4)$. Now identifying the vertical boundary of $M(3,3,0)$ and cutting along the path $P(x_{1},x_{4},x_{7},x_{2})$ leads to Figure 5.3.6, i.e., $M(3,3,2)$. By the isomorphism map define in Lemma \ref{l4.3.2}, $M(3,3,0) \cong M(3,3,2)$. Therefore, there are two DSEMs of type $T_3$ with $9$ vertices on the torus upto isomorphism.
	
\end{eg}

\begin{picture}(0,0)(-2,25.5)
\setlength{\unitlength}{7mm}

\drawpolygon(0,1)(3,1)(3,4)(0,4)


\drawline[AHnb=0](1,1)(1,4)
\drawline[AHnb=0](2,1)(2,4)


\drawline[AHnb=0](0,1)(3,1)
\drawline[AHnb=0](0,2)(3,2)
\drawline[AHnb=0](0,3)(3,3)


\drawline[AHnb=0](0,3)(1,4)
\drawline[AHnb=0](1,3)(2,4)
\drawline[AHnb=0](2,3)(3,4)

\put(-.3,.8){\scriptsize {\tiny $x_1$}}
\put(1,.75){\scriptsize {\tiny $x_2$}}
\put(2,.75){\scriptsize {\tiny $x_3$}}
\put(3,.75){\scriptsize {\tiny $x_1$}}

\put(-.5,1.75){\scriptsize {\tiny $x_4$}}
\put(1.1,1.7){\scriptsize {\tiny $x_5$}}
\put(2.1,1.7){\scriptsize {\tiny $x_6$}}
\put(3.1,1.7){\scriptsize {\tiny $x_4$}}

\put(-.5,2.7){\scriptsize {\tiny $x_7$}}
\put(1.1,2.7){\scriptsize {\tiny $x_8$}}
\put(2.1,2.7){\scriptsize {\tiny $x_9$}}
\put(3.1,2.7){\scriptsize {\tiny $x_7$}}

\put(-.3,4.12){\scriptsize {\tiny $x_1$}}
\put(.9,4.12){\scriptsize {\tiny $x_2$}}
\put(1.9,4.12){\scriptsize {\tiny $x_3$}}
\put(3,4.12){\scriptsize {\tiny $x_1$}}

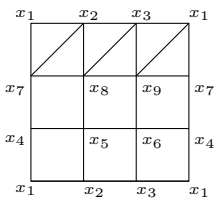
\put(-.5,0){\scriptsize {\tiny {\bf Figure 5.3.3:} $M(3,3,0)$ }} 

\end{picture}

\begin{picture}(0,0)(-8,20.5)
\setlength{\unitlength}{7mm}

\drawpolygon(5,1)(8,1)(8,4)(5,4)


\drawline[AHnb=0](6,1)(6,4)
\drawline[AHnb=0](7,1)(7,4)


\drawline[AHnb=0](5,1)(8,1)
\drawline[AHnb=0](5,2)(8,2)
\drawline[AHnb=0](5,3)(8,3)


\drawline[AHnb=0](5,3)(6,4)
\drawline[AHnb=0](6,3)(7,4)
\drawline[AHnb=0](7,3)(8,4)

\put(4.6,.75){\scriptsize {\tiny $y_1$}}
\put(5.9,.75){\scriptsize {\tiny $y_2$}}
\put(6.9,.75){\scriptsize {\tiny $y_3$}}
\put(8,.75){\scriptsize {\tiny $y_1$}}

\put(4.5,1.7){\scriptsize {\tiny $y_4$}}
\put(6.1,1.7){\scriptsize {\tiny $y_5$}}
\put(7.1,1.7){\scriptsize {\tiny $y_6$}}
\put(8.1,1.7){\scriptsize {\tiny $y_4$}}

\put(4.5,2.7){\scriptsize {\tiny $y_7$}}
\put(6.1,2.7){\scriptsize {\tiny $y_8$}}
\put(7.1,2.7){\scriptsize {\tiny $y_9$}}
\put(8.1,2.7){\scriptsize {\tiny $y_7$}}

\put(4.6,4.12){\scriptsize {\tiny $y_2$}}
\put(5.8,4.12){\scriptsize {\tiny $y_3$}}
\put(6.8,4.12){\scriptsize {\tiny $y_1$}}
\put(8,4.12){\scriptsize {\tiny $y_2$}}

\put(4.6,-.01){\scriptsize {\tiny {\bf Figure 5.3.4:} $M(3,3,1)$ }} 

\end{picture}

\begin{picture}(0,0)(-18,16)
\setlength{\unitlength}{7mm}

\drawpolygon(10,1)(13,1)(13,4)(10,4)


\drawline[AHnb=0](11,1)(11,4)
\drawline[AHnb=0](12,1)(12,4)

\drawline[AHnb=0](10,1)(13,1)
\drawline[AHnb=0](10,2)(13,2)
\drawline[AHnb=0](10,3)(13,3)


\drawline[AHnb=0](10,3)(11,4)
\drawline[AHnb=0](11,3)(12,4)
\drawline[AHnb=0](12,3)(13,4)

\put(9.6,.75){\scriptsize {\tiny $z_1$}}
\put(11,.75){\scriptsize {\tiny $z_2$}}
\put(12,.75){\scriptsize {\tiny $z_3$}}
\put(13,.75){\scriptsize {\tiny $z_1$}}

\put(9.5,1.7){\scriptsize {\tiny $z_4$}}
\put(11.1,1.7){\scriptsize {\tiny $z_5$}}
\put(12.1,1.7){\scriptsize {\tiny $z_6$}}
\put(13.1,1.7){\scriptsize {\tiny $z_4$}}

\put(9.5,2.7){\scriptsize {\tiny $z_7$}}
\put(11.1,2.7){\scriptsize {\tiny $z_8$}}
\put(12.1,2.7){\scriptsize {\tiny $z_9$}}
\put(13.1,2.7){\scriptsize {\tiny $z_1$}}

\put(9.6,4.1){\scriptsize {\tiny $z_3$}}
\put(10.9,4.1){\scriptsize {\tiny $z_1$}}
\put(11.9,4.1){\scriptsize {\tiny $z_2$}}
\put(13.1,4.1){\scriptsize {\tiny $z_3$}}

\put(9.6,0){\scriptsize {\tiny {\bf Figure 5.3.5:} $M(3,3,2)$ }} 

\end{picture}

\begin{picture}(0,0)(-58,11)
\setlength{\unitlength}{7mm}

\drawpolygon(10,1)(13,1)(13,4)(10,4)


\drawline[AHnb=0](11,1)(11,4)
\drawline[AHnb=0](12,1)(12,4)

\drawline[AHnb=0](10,1)(13,1)
\drawline[AHnb=0](10,2)(13,2)
\drawline[AHnb=0](10,3)(13,3)


\drawline[AHnb=0](10,3)(11,4)
\drawline[AHnb=0](11,3)(12,4)
\drawline[AHnb=0](12,3)(13,4)

\put(9.8,.75){\scriptsize {\tiny $x_1$}}
\put(11,.75){\scriptsize {\tiny $x_3$}}
\put(12,.75){\scriptsize {\tiny $x_2$}}
\put(13,.75){\scriptsize {\tiny $x_1$}}

\put(9.5,1.7){\scriptsize {\tiny $x_4$}}
\put(11.1,1.7){\scriptsize {\tiny $x_6$}}
\put(12.1,1.7){\scriptsize {\tiny $x_5$}}
\put(13.1,1.7){\scriptsize {\tiny $x_4$}}

\put(9.5,2.7){\scriptsize {\tiny $x_7$}}
\put(11.1,2.7){\scriptsize {\tiny $x_9$}}
\put(12.1,2.7){\scriptsize {\tiny $x_8$}}
\put(13.1,2.7){\scriptsize {\tiny $x_7$}}

\put(9.6,4.1){\scriptsize {\tiny $x_2$}}
\put(10.9,4.1){\scriptsize {\tiny $x_1$}}
\put(11.9,4.1){\scriptsize {\tiny $x_3$}}
\put(13.1,4.1){\scriptsize {\tiny $x_2$}}

\put(9.6,0){\scriptsize {\tiny {\bf Figure 5.3.6:} $M(3,3,2)$ }} 

\end{picture}



\vspace{.85cm}

\begin{center}
	
	\noindent \textbf{Table \ref{table:3}}: DSEMs of type $T_3$ on the torus for $|V(M)| \leq 18$ 
	
	\renewcommand{\arraystretch}{1.1}
	\begin{tabular}{|p{1.2cm}|p{4cm}|p{3.8cm}|p{3cm}|}
		
		\hline
		$|V(M)|$ & Isomorphic classes & Length of cycles & No of maps\\
		
		\hline
		9 &  $M(3,3,0)$, $M(3,3,2)$ & $( 3, \{ 3, 9 \}, 3 )$ & 2\\
		
		\cline{2-3}
		
		& $M(3,3,1)$ & $( 3, \{ 9, 9 \}, 4 )$ & \\
		\cline{2-3}

		\hline
		12 &  $M(4,3,0)$, $M(4,3,3)$ & $( 4, \{ 3, 12 \}, 3 )$& 2\\
		\cline{2-3}
		& $M(4,3,1)$, $M(4,3,2)$ & $( 4, \{ 6,12 \}, 4 )$ & \\
		
		\hline

		\hline
		15 & $M(5,3,0)$, $M(5,3,4)$ & $( 5, \{ 3, 15 \}, 3 )$ &3\\
		\cline{2-3}
		& $M(5,3,1)$, $M(5,3,3)$ & $(5, \{ 15, 15 \}, 4)$ &\\
		\cline{2-3}
		& $M(5,3,2)$ & $(5, \{ 15, 15 \}, 5)$  &\\
		\cline{2-3}
		\hline
			\end{tabular}
	
\end{center}

\begin{center}
\begin{tabular}{ |p{1.2cm}|p{4cm}|p{3.8cm}|p{3cm}| }

		\hline
		18 & $M(6,3,0)$, $M(6,3,5)$& $( 6, \{ 3, 18 \}, 3 )$ & 5\\
		\cline{2-3}
		&$M(6,3,1)$, $M(6,3,4)$ &$( 6, \{ 9, 18 \}, 4 )$&\\
		\cline{2-3}
		& $M(6,3,2)$, $M(6,3,3)$ &$( 6, \{ 6, 9 \}, 5)$&\\
		\cline{2-3}
		& $M(3,6,0)$, $M(3,6,1)$ &$(3, \{ 6,18 \}, 6 )$&\\
		\cline{2-3}
		& $M(3,6,2)$ & $( 3, \{ 18, 18 \}, 8 )$ &  \\
		\cline{2-3}
		\hline

	\end{tabular}
	\label{table:3}
\end{center}

\subsection{DSEMs of type $T_4= [f_1^{({f_1}^4.{f_2}^3)}:f_2^{({f_1}^3.{f_2}^5)}]$, where $(f_1, f_2)=((3^3.4^2),(4^4))$} \label{s4.4}

Let $M_4$ be a DSEM of type $T_4$ with the vertex set $V(M_4)$. For the existence of $M_4$, we have $|V_{(4^4)}| = |V_{(3^3.4^2)}| = 2k$ for some $k \in \mathbb{N}$.  Now we define following three types of paths in $M_4$ as follows.

\begin{picture}(60,0)(-17,23)
\setlength{\unitlength}{5mm}

\drawpolygon(0,0)(6,0)(6,3)(0,3)


\drawline[AHnb=0](1,0)(1,3)
\drawline[AHnb=0](2,0)(2,3)
\drawline[AHnb=0](3,0)(3,3)
\drawline[AHnb=0](4,0)(4,3)
\drawline[AHnb=0](5,0)(5,3)

\drawline[AHnb=0](.5,4)(.5,4.5)
\drawline[AHnb=0](1.5,4)(1.5,4.5)

\drawline[AHnb=0](2.5,4)(2.5,4.5)
\drawline[AHnb=0](3.5,4)(3.5,4.5)

\drawline[AHnb=0](4.5,4)(4.5,4.5)
\drawline[AHnb=0](5.5,4)(5.5,4.5)


\drawline[AHnb=0](6,0)(6.5,0)
\drawline[AHnb=0](0,0)(-.5,0)
\drawline[AHnb=0](6,3)(6.5,3)
\drawline[AHnb=0](0,3)(-.5,3)

\drawline[AHnb=0](-.5,1)(6.5,1)
\drawline[AHnb=0](-.5,2)(6.5,2)

\drawline[AHnb=0](-.5,4)(6.5,4)


\drawline[AHnb=0](0,3)(.5,4)
\drawline[AHnb=0](1,3)(.5,4)

\drawline[AHnb=0](1,3)(1.5,4)
\drawline[AHnb=0](2,3)(1.5,4)

\drawline[AHnb=0](2,3)(2.5,4)
\drawline[AHnb=0](3,3)(2.5,4)

\drawline[AHnb=0](3,3)(3.5,4)
\drawline[AHnb=0](4,3)(3.5,4)

\drawline[AHnb=0](4,3)(4.5,4)
\drawline[AHnb=0](5,3)(4.5,4)

\drawline[AHnb=0](5,3)(5.5,4)
\drawline[AHnb=0](6,3)(5.5,4)


\drawline[AHnb=0](0,0)(.25,-.5)
\drawline[AHnb=0](1,0)(1.25,-.5)
\drawline[AHnb=0](2,0)(2.25,-.5)
\drawline[AHnb=0](3,0)(3.25,-.5)
\drawline[AHnb=0](4,0)(4.25,-.5)
\drawline[AHnb=0](5,0)(5.25,-.5)
\drawline[AHnb=0](6,0)(6.25,-.5)

\drawline[AHnb=0](0,0)(-.25,-.5)
\drawline[AHnb=0](1,0)(.75,-.5)
\drawline[AHnb=0](2,0)(1.75,-.5)
\drawline[AHnb=0](3,0)(2.75,-.5)
\drawline[AHnb=0](4,0)(3.75,-.5)
\drawline[AHnb=0](5,0)(4.75,-.5)
\drawline[AHnb=0](6,0)(5.75,-.5)

\drawline[AHnb=0](0,3)(-.25,3.5)
\drawline[AHnb=0](6,3)(6.25,3.5)


\drawpolygon[fillcolor=black](-.5,0)(6.5,0)(6.5,.1)(-.5,.1)

\drawpolygon[fillcolor=black](-.5,1)(6.5,1)(6.5,1.1)(-.5,1.1)

\drawpolygon[fillcolor=black](-.5,2)(6.5,2)(6.5,2.1)(-.5,2.1)

\drawpolygon[fillcolor=black](-.5,3)(6.5,3)(6.5,3.1)(-.5,3.1)

\drawpolygon[fillcolor=black](-.5,4)(6.5,4)(6.5,4.1)(-.5,4.1)

\drawpolygon[fillcolor=green](-.25,-.5)(0,0)(0,3)(.5,4)(.5,4.5)(.4,4.5)(.4,4)(-.1,3)(-.1,0)(-.35,-.5)

\drawpolygon[fillcolor=green](.75,-.5)(1,0)(1,3)(1.5,4)(1.5,4.5)(1.4,4.5)(1.4,4)(.9,3)(.9,0)(.65,-.5)

\drawpolygon[fillcolor=green](1.75,-.5)(2,0)(2,3)(2.5,4)(2.5,4.5)(2.4,4.5)(2.4,4)(1.9,3)(1.9,0)(1.65,-.5)

\drawpolygon[fillcolor=green](2.75,-.5)(3,0)(3,3)(3.5,4)(3.5,4.5)(3.4,4.5)(3.4,4)(2.9,3)(2.9,0)(2.65,-.5)

\drawpolygon[fillcolor=green](3.75,-.5)(4,0)(4,3)(4.5,4)(4.5,4.5)(4.4,4.5)(4.4,4)(3.9,3)(3.9,0)(3.65,-.5)

\drawpolygon[fillcolor=green](4.75,-.5)(5,0)(5,3)(5.5,4)(5.5,4.5)(5.4,4.5)(5.4,4)(4.9,3)(4.9,0)(4.65,-.5)

\drawpolygon[fillcolor=green](5.75,-.5)(6,0)(6,3)(6.25,3.5)(6.15,3.5)(5.9,3)(5.9,0)(5.65,-.5)

\drawpolygon[fillcolor=red](.25,-.5)(0,0)(0,3)(-.25,3.5)(-.15,3.5)(.1,3)(.1,0)(.35,-.5)

\drawpolygon[fillcolor=red](1.25,-.5)(1,0)(1,3)(.5,4)(.5,4.5)(.6,4.5)(.6,4)(1.1,3)(1.1,0)(1.35,-.5)

\drawpolygon[fillcolor=red](2.25,-.5)(2,0)(2,3)(1.5,4)(1.5,4.5)(1.6,4.5)(1.6,4)(2.1,3)(2.1,0)(2.35,-.5)

\drawpolygon[fillcolor=red](3.25,-.5)(3,0)(3,3)(2.5,4)(2.5,4.5)(2.6,4.5)(2.6,4)(3.1,3)(3.1,0)(3.35,-.5)

\drawpolygon[fillcolor=red](4.25,-.5)(4,0)(4,3)(3.5,4)(3.5,4.5)(3.6,4.5)(3.6,4)(4.1,3)(4.1,0)(4.35,-.5)

\drawpolygon[fillcolor=red](5.25,-.5)(5,0)(5,3)(4.5,4)(4.5,4.5)(4.6,4.5)(4.6,4)(5.1,3)(5.1,0)(5.35,-.5)

\drawpolygon[fillcolor=red](6.25,-.5)(6,0)(6,3)(5.5,4)(5.5,4.5)(5.6,4.5)(5.6,4)(6.1,3)(6.1,0)(6.35,-.5)

\put(2.4,4.75) {\scriptsize {\tiny $Y_3$}}
\put(3.4,4.75) {\scriptsize {\tiny $Y_2$}}
\put(7.1,-.1) {\scriptsize {\tiny $Y_1$}}

\put(-.95,-1.8){\scriptsize {\tiny Figure 5.4.1 : Paths of type $Y_1, Y_2, Y_3$ }} 

\end{picture}

\begin{picture}(0,0)(-66,20)
\setlength{\unitlength}{6.5mm}

\drawpolygon(0,0)(3,0)(3,4)(0,4)
\drawpolygon(5,0)(8,0)(8,4)(5,4)
\drawpolygon(10,0)(12,0)(12,4)(10,4)


\drawline[AHnb=0](1,0)(1,4)
\drawline[AHnb=0](2,0)(2,4)
\drawline[AHnb=0](3,0)(3,4)
\drawline[AHnb=0](6,0)(6,4)
\drawline[AHnb=0](7,0)(7,4)
\drawline[AHnb=0](11,0)(11,4)


\drawline[AHnb=0](0,1)(3,1)
\drawline[AHnb=0](0,2)(3,2)
\drawline[AHnb=0](0,3)(3,3)

\drawline[AHnb=0](5,1)(8,1)
\drawline[AHnb=0](5,2)(8,2)
\drawline[AHnb=0](5,3)(8,3)

\drawline[AHnb=0](10,1)(12,1)
\drawline[AHnb=0](10,2)(12,2)
\drawline[AHnb=0](10,3)(12,3)


\drawline[AHnb=0](0,3)(1,4)
\drawline[AHnb=0](1,3)(2,4)
\drawline[AHnb=0](2,3)(3,4)

\drawline[AHnb=0](5,3)(6,4)
\drawline[AHnb=0](6,3)(7,4)
\drawline[AHnb=0](7,3)(8,4)

\drawline[AHnb=0](10,3)(11,4)
\drawline[AHnb=0](11,3)(12,4)

\put(-.2,-.4){\scriptsize {\tiny $u_1$}}
\put(.9,-.4){\scriptsize {\tiny $u_2$}}

\put(1.9,-.4){\scriptsize {\tiny $u_3$}}
\put(2.8,-.4){\scriptsize {\tiny $u_4$}}
\put(4.8,-.4){\scriptsize {\tiny $u_{k}$}}
\put(5.6,-.4){\scriptsize {\tiny $u_{k+1}$}}
\put(6.6,-.4){\scriptsize {\tiny $u_{k+2}$}}
\put(7.6,-.4){\scriptsize {\tiny $u_{k+3}$}}
\put(9.8,-.4){\scriptsize {\tiny $u_{i-1}$}}
\put(10.9,-.4){\scriptsize {\tiny $u_{i}$}}
\put(11.9,-.4){\scriptsize {\tiny $u_1$}}

\put(-.4,.8){\scriptsize {\tiny $v_1$}}
\put(.6,.75){\scriptsize {\tiny $v_2$}}
\put(1.6,.75){\scriptsize {\tiny $v_3$}}
\put(2.6,.75){\scriptsize {\tiny $v_4$}}
\put(5.1,.75){\scriptsize {\tiny $v_{k}$}}
\put(6.1,.75){\scriptsize {\tiny $v_{k+1}$}}
\put(7.1,.75){\scriptsize {\tiny $v_{k+2}$}}
\put(8.2,.75){\scriptsize {\tiny $v_{k+3}$}}

\put(10.1,.75){\scriptsize {\tiny $v_{i-2}$}}
\put(11.1,.75){\scriptsize {\tiny $v_{i}$}}
\put(12.1,.75){\scriptsize {\tiny $v_{1}$}}

\put(-.44,1.75){\scriptsize {\tiny $w_1$}}
\put(.5,1.7){\scriptsize {\tiny $w_2$}}
\put(1.5,1.7){\scriptsize {\tiny $w_3$}}
\put(2.5,1.7){\scriptsize {\tiny $w_4$}}
\put(5.1,1.7){\scriptsize {\tiny $w_{k}$}}
\put(6.05,1.7){\scriptsize {\tiny $w_{k+1}$}}
\put(7.05,1.7){\scriptsize {\tiny $w_{k+2}$}}
\put(8.1,1.7){\scriptsize {\tiny $w_{k+3}$}}
\put(10.1,1.7){\scriptsize {\tiny $w_{i-2}$}}
\put(11.1,1.7){\scriptsize {\tiny $w_{i}$}}
\put(12.1,1.7){\scriptsize {\tiny $w_{1}$}}

\put(-.4,2.7){\scriptsize {\tiny $x_1$}}
\put(.6,2.7){\scriptsize {\tiny $x_2$}}
\put(1.6,2.7){\scriptsize {\tiny $x_3$}}
\put(2.6,2.7){\scriptsize {\tiny $x_4$}}
\put(5.1,2.7){\scriptsize {\tiny $x_{k}$}}
\put(6.1,2.7){\scriptsize {\tiny $x_{k+1}$}}
\put(7.05,2.7){\scriptsize {\tiny $x_{k+2}$}}
\put(8.1,2.7){\scriptsize {\tiny $x_{k+3}$}}
\put(10.1,2.75){\scriptsize {\tiny $x_{i-1}$}}
\put(11.2,2.7){\scriptsize {\tiny $x_{i}$}}
\put(12.1,2.7){\scriptsize {\tiny $x_1$}}

\put(-.5,4.2){\scriptsize {\tiny $u_{k+1}$}}
\put(.5,4.2){\scriptsize {\tiny $u_{k+2}$}}
\put(1.5,4.2){\scriptsize {\tiny $u_{k+3}$}}
\put(2.5,4.2){\scriptsize {\tiny $u_{k+4}$}}
\put(4.8,4.2){\scriptsize {\tiny $u_{i}$}}
\put(5.7,4.2){\scriptsize {\tiny $u_{1}$}}
\put(6.7,4.2){\scriptsize {\tiny $u_{2}$}}
\put(7.8,4.2){\scriptsize {\tiny $u_{3}$}}
\put(9.8,4.2){\scriptsize {\tiny $u_{k-1}$}}
\put(10.8,4.2){\scriptsize {\tiny $u_{k}$}}
\put(11.5,4.2){\scriptsize {\tiny $u_{k+1}$}}

\put(3.8,0){\scriptsize $\ldots$}
\put(3.8,0){\scriptsize $\ldots$}
\put(3.8,1){\scriptsize $\ldots$}
\put(3.8,1){\scriptsize $\ldots$}
\put(3.8,2){\scriptsize $\ldots$}
\put(3.8,2){\scriptsize $\ldots$}
\put(3.8,3){\scriptsize $\ldots$}
\put(3.8,3){\scriptsize $\ldots$}
\put(3.8,3.8){\scriptsize $\ldots$}
\put(3.8,3.8){\scriptsize $\ldots$}

\put(9.1,0){\scriptsize $\ldots$}
\put(9.1,0){\scriptsize $\ldots$}
\put(9.1,1){\scriptsize $\ldots$}
\put(9.1,1){\scriptsize $\ldots$}
\put(9.1,2){\scriptsize $\ldots$}
\put(9.1,2){\scriptsize $\ldots$}
\put(9.1,3){\scriptsize $\ldots$}
\put(9.1,3){\scriptsize $\ldots$}
\put(9.1,3.8){\scriptsize $\ldots$}
\put(9.1,3.8){\scriptsize $\ldots$}

\put(3,-1.2){\scriptsize {\tiny {\bf Figure 5.4.2:} $M(i,4,k)$ representation}} 

\end{picture}

\vspace{3cm}




\begin{deff} \normalfont
	
	A path $P_{1} = P( \ldots, y_{i-1},y_{i},y_{i+1}, \ldots)$ in $M_4$ is of type $Y_{1}$ if either of the conditions follows: 
	$(i)$ every vertex of $P_1$ has the face-sequence $(4^4)$ or
	$(ii)$ each vertex of $P_1$ has the face-sequence $(3^3.4^2)$ such that all the triangles (square) that incident on its inner vertices lie on the one side  of the path $P_1$. For example, see the paths indicated by black color in Figure 5.4.1.
	
\end{deff}

\begin{deff} \normalfont A path $P_{2} = P(\ldots, z_{i-1},z_{i},z_{i+1}, \ldots)$ in $M_4$, such that $z_{i-1}, z_{i}, z_{i+1}$ are inner vertices of $P_{2}$ or an extended path of $P_{2}$, is of type $Y_{2}$ (shown by green colored  paths in Figure 5.4.1), if either of the following four conditions follows for each vertex of the path.

	\begin{enumerate} 
		\item If ${\rm lk}(z_{i})=C_7(\boldsymbol{m},z_{i-1},\boldsymbol{n},o,z_{i+1},p,q)$ and ${\rm lk}(z_{i-1})=C_8(\boldsymbol{q},m,\boldsymbol{r},z_{i-2},\boldsymbol{s},n,\boldsymbol{o},z_i)$, then ${\rm lk}(z_{i+1}) \linebreak=C_7(\boldsymbol{t},z_{i+2}, \boldsymbol{u},p,z_{i},o,v)$, ${\rm lk}(z_{i+2})=C_8(\boldsymbol{w},z_{i+3}, \boldsymbol{x},u,\boldsymbol{p},z_{i+1},\boldsymbol{v},t)$. 
		
		\item If ${\rm lk}(z_{i})=C_7(\boldsymbol{m},z_{i+1},\boldsymbol{n},o,z_{i-1},p,q)$ and ${\rm lk}(z_{i-1})=C_7(\boldsymbol{s},z_{i-2},\boldsymbol{t},p,z_{i},o,r)$, then ${\rm lk}(z_{i+1})=C_8(\boldsymbol{o},z_{i}, \boldsymbol{q},m, \boldsymbol{u},z_{i+2},\boldsymbol{v},n)$, ${\rm lk}(z_{i+2})=C_8(\boldsymbol{n},z_{i+1}, \boldsymbol{m},u,\boldsymbol{w},z_{i+3},\boldsymbol{x},v)$.
		
		\item If  ${\rm lk}(z_{i})=C_8(\boldsymbol{m},z_{i-1}, \boldsymbol{n},o,\boldsymbol{p},z_{i+1},\boldsymbol{q},r)$ and ${\rm lk}(z_{i-1})=C_7(\boldsymbol{o},z_{i},\boldsymbol{r},m,z_{i-2},s,n)$, then ${\rm lk}(z_{i+1}) \linebreak=C_8(\boldsymbol{r},z_{i}, \boldsymbol{o},p,\boldsymbol{t},z_{i+2},\boldsymbol{u},q)$, ${\rm lk}(z_{i+2})=C_7(\boldsymbol{q},z_{i+1}, \boldsymbol{p},t,z_{i+3},v,u)$.
		
		\item If ${\rm lk}(z_{i})=C_8(\boldsymbol{m},z_{i-1}, \boldsymbol{n},o,\boldsymbol{p},z_{i+1},\boldsymbol{q},r)$ and ${\rm lk}(z_{i-1})=C_8(\boldsymbol{o},z_{i}, \boldsymbol{r},m,\boldsymbol{s},z_{i-2},\boldsymbol{t},n)$, then \linebreak 
		${\rm lk}(z_{i+1}) =C_7(\boldsymbol{r},z_{i},\boldsymbol{o},p,z_{i+2},v,q)$, ${\rm lk}(z_{i+2})=C_7(\boldsymbol{w}, z_{i+3},\boldsymbol{x},v,z_{i+1},p,u)$.
		
	\end{enumerate}
\end{deff}

\begin{deff} \normalfont A path $P_{3} = P(\ldots, w_{i-1},w_{i},w_{i+1}, \ldots)$ in $M_4$, such that $w_{i-1}, w_{i}, w_{i+1}$ are inner vertices of $P_{3}$ or an extended path of $P_{3}$, is of type $Y_{3}$ (shown by red colored  paths in Figure 5.4.1), if either of the following four conditions follows for each vertex of the path.
	
	\begin{enumerate} 
		
		\item If ${\rm lk}(z_{i})=C_7(\boldsymbol{m},z_{i-1},\boldsymbol{n},o,p,z_{i+1},q)$ and ${\rm lk}(z_{i-1})=C_8(\boldsymbol{q},m,\boldsymbol{r},z_{i-2},\boldsymbol{s},n,\boldsymbol{o},z_i)$, then ${\rm lk}(z_{i+1}) \linebreak=C_7(\boldsymbol{t},z_{i+2}, \boldsymbol{u},v,q,z_{i},p)$, ${\rm lk}(z_{i+2})=C_8(\boldsymbol{x},z_{i+3}, \boldsymbol{y},u,\boldsymbol{v},w_{i+1},\boldsymbol{p},t)$. 
		
		\item If ${\rm lk}(z_{i})=C_7(\boldsymbol{m},z_{i+1},\boldsymbol{n},o,p,z_{i-1},q)$ and ${\rm lk}(z_{i-1})=C_7(\boldsymbol{r},z_{i-2},\boldsymbol{s},t,q,z_{i},p)$, then ${\rm lk}(z_{i+1})=C_8(\boldsymbol{o},z_{i}, \boldsymbol{q},m, \boldsymbol{u},z_{i+2},\boldsymbol{v},n)$, ${\rm lk}(z_{i+2})=C_8(\boldsymbol{n},z_{i+1}, \boldsymbol{m},u,\boldsymbol{x},z_{i+3},\boldsymbol{y},v)$.

		\item If  ${\rm lk}(z_{i})=C_8(\boldsymbol{m},z_{i-1}, \boldsymbol{n},o,\boldsymbol{p},z_{i+1},\boldsymbol{q},r)$ and ${\rm lk}(z_{i-1})=C_7(\boldsymbol{o},z_{i},\boldsymbol{r},m,s,z_{i-2},n)$, then ${\rm lk}(z_{i+1}) \linebreak=C_8(\boldsymbol{r},z_{i}, \boldsymbol{o},p,\boldsymbol{t},z_{i+2},\boldsymbol{u},q)$, ${\rm lk}(z_{i+2})=C_7(\boldsymbol{q},z_{i+1}, \boldsymbol{p},t,v,z_{i+3},u)$.
		
		\item If ${\rm lk}(z_{i})=C_8(\boldsymbol{m},z_{i-1}, \boldsymbol{n},o,\boldsymbol{p},z_{i+1},\boldsymbol{q},r)$ and ${\rm lk}(z_{i-1})=C_8(\boldsymbol{o},z_{i}, \boldsymbol{r},m,\boldsymbol{s},z_{i-2},\boldsymbol{t},n)$, then \linebreak 
		${\rm lk}(z_{i+1}) =C_7(\boldsymbol{r},z_{i},\boldsymbol{o},p,u,z_{i+2},q)$, ${\rm lk}(z_{i+2})=C_7(\boldsymbol{x}, z_{i+3},\boldsymbol{y},v,q,z_{i+1},u)$.
	\end{enumerate}
\end{deff}

Consider a maximal path $P(v_1, v_2, \ldots, v_i)$ of the type $Y_\alpha$, $\alpha \in \{1,2,3\}$. Following Remark \ref{r1} and Lemma \ref{l4.1.3}, there is an edge $v_i$-$v_1$ in $M_4$ such that  $P(v_1, v_2, \ldots, v_i) \cup \{v_i$-$v_1\}$ is a non-contractible cycle $Q=C_i(v_1, v_2, \ldots, v_i)$. The cycles of type $Y_{2}$ and $Y_{3}$ define same type cycles as they are mirror image of each other. By Equation (1) of Section \ref{s4.1}, define a cycle of type $Y_4$ in $M_4$. Let $u \in V(M_4)$ and $Q_{\alpha}$ be cycles of type $Y_{\alpha}$ through $u$. As in Section \ref{s4.1}, we define an $M(i,j,k)$ representation of $M_4$ for some $i, j, k$. For this, we first cut $M_4$ along the cycle $Q_{1}$ and then cut it along the cycle $Q_{3}$. Without loss of generality, assume that the faces incident on the base horizontal cycle are quadrangular, for examples see $M(i,4,k)$ in Figure 5.4.2. Then we prove the following lemma.

\begin{lem}\label{l4.4.1} The DSEM $M$ of type $T_4$ admits an $M(i,j,k)$ representation iff the following holds: $(i)$ $ i \geq 3$ and $j=4m$, where $m \in \mathbb{N} $, $(ii)$ $ij \geq 12 $, $(iii)$ $0 \leq k \leq i-1$.
\end{lem}
\noindent{\bf Proof.} Let $M$ be a DSEM of  type $T_4$ with $n$ vertices. Then by using the fact that $|V(3^3.4^2)|$ = $|V(4^4)|$, and proceeding as in Lemma \ref{l4.2.1}, we get all the admissible value of $i,j$ and $k$ of $M(i,j,k)$. Thus the proof. \hfill $\Box$

\smallskip

For $ t \in \{1,2\}$, let $M_{t}$ be DSEMs of type $T_4$ on $n_t$ number of vertices with representation $M_{t}(i_{t}, j_{t}, k_{t})$. Suppose $n_1=n_2$. Let $Q^t_{\alpha}$ be the cycles of type $Y_{\alpha}$ and  $l^t_{\alpha}$ = length of the cycle of type $Y_{\alpha}$, for $\alpha \in \{1,2,3,4\}$, in $M_{t}(i_t,j_t,k_t)$. We say that $M_{t}(i_t,j_t,k_t)$ has cycle-type $(l^t_1, l^t_2, l^t_3, l^t_4)$ if $l^t_2 \leq l^t_3$ or $(l^t_1, l^t_3, l^t_2, l_{t,4})$  if $l^t_3 < l^t_2$. 
Similarly as in Section \ref{s4.1}, we have

\begin{lem}\label{l4.4.2} The DSEMs $M_{1} \cong  M_{2}$ iff they have same cycle-type.
\end{lem}



\smallskip


Now computing the DSEMs for the first four admissible values of $|V(M_4)|$, we get Table \ref{table:4}. We illustrate this computation for $|V(M_4)| = 12$ as follows. 

\begin{eg} \normalfont Let $M_4$ be a DSEM of type $T_4$ with $12$ vertices on the torus. By Lemma \ref{l4.4.1}, $M_4$ has three $M(i,j,k)$ representation, namely, $M(3, 4, 0), M(3, 4, 1)$ and $M(3, 4, 2)$, see Figures 5.4.3, 5.4.4, and 5.4.5 respectively. In $M(3, 4, 0)$, $Q^1_1 = C_3(x_{1}, x_{2}, x_{3})$ is a $Y_{1}$ type cycle, $Q^1_2 =C_{12}(x_{1},x_{4},x_{7},x_{10},x_2,x_5,x_8,x_{11},x_3,x_6,x_9,x_{12})$ and $Q^1_3 = C_4(x_{1},x_{4},x_{7},x_{10})$ are $Y_{2}$ type cycles and $Q^1_4 = C_4(x_{1},x_{4},x_{7},x_{10})$ is a $Y_{4}$ type cycle. In $M( 3,4,1)$, $Q^2_1 = C_3(y_{1}, y_{2}, y_{3})$ is a $Y_{1}$ type cycle, $Q^2_2 = C_{12}(y_{1},y_{4},y_{7},y_{10},y_{3},y_{6},y_{9},y_{12},y_{2},y_{5},y_{8},y_{11})$ and $Q^2_3 = C_{12}(y_{1},y_{4},y_{7},y_{10},y_{2},y_{5},y_{8},y_{11},y_{3}$, $y_{6},y_{9},y_{12})$ are $Y_{2}$ type cycles and $Q^2_4 = C_5(y_{2},y_{5},y_{8},y_{11}, y_{1})$ is a $Y_{4}$ type cycle. In $M( 3,4,2)$, $Q^3_1 = C_{3}(z_{1},z_{2},z_{3})$ is a $Y_{1}$ type cycle, $Q^3_2 = C_{4}(z_{1},z_{4},z_{7},z_{10})$ and $Q^3_3 = C_{12}(z_{1},z_{4},z_{7},z_{10},z_{3},z_{6},z_{9},z_{12}$, $z_{2},z_{5},z_{8},z_{11})$ are $Y_{2}$ type cycles and $Q^3_4 = C_4(z_{3},z_{6},z_{9},z_{12})$ is a $Y_{4}$ type cycle. 	
	
	In $M(i,j,k)$, observe that type $Y_{1}$ cycles have the same length and type $Y_{2}$ cycles have at most two different lengths. Since length$(Q^2_4)$ $\neq$ length$(Q^r_{4})$ for $r \in \{1,3\}$, $M(3, 4, 1) \ncong  M(3, 4, 0), \linebreak M(3, 4, 2)$. Observe that, length $(Q^1_1)$ = length$(Q^3_1)$, \{ length$(Q^1_2)$, length$(Q^1_3)$\} = \{  length $(Q^3_2)$, length $(Q^3_3)$\} and length$(Q^1_4)$ =length$(Q^3_4)$. Now identifying $M(3,4,0)$ along the vertical boundary and cutting along
	the path $P(x_{1},x_{4},x_{17},x_{10},x_{2})$ leads to Figure 5.4.6, i.e., $M(3,4,2)$. By the isomorphism map define in Lemma \ref{l4.4.2}, $M(3,4,0) \cong M(3,4,2)$. Therefore, there are two DSEMs, up to isomorphism, of type $T_4$ with $12$ vertices on the torus.
	
\end{eg}

\begin{picture}(0,0)(-8,20)
\setlength{\unitlength}{5.8mm}

\drawpolygon(0,0)(3,0)(3,4)(0,4)


\drawline[AHnb=0](1,0)(1,4)
\drawline[AHnb=0](2,0)(2,4)


\drawline[AHnb=0](0,1)(3,1)
\drawline[AHnb=0](0,2)(3,2)
\drawline[AHnb=0](0,3)(3,3)


\drawline[AHnb=0](0,3)(1,4)
\drawline[AHnb=0](1,3)(2,4)
\drawline[AHnb=0](2,3)(3,4)

\put(-.2,-.4){\scriptsize {\tiny $x_1$}}
\put(.9,-.4){\scriptsize {\tiny $x_2$}}
\put(1.9,-.4){\scriptsize {\tiny $x_3$}}
\put(2.8,-.4){\scriptsize {\tiny $x_1$}}

\put(-.5,.8){\scriptsize {\tiny $x_4$}}
\put(1.1,.7){\scriptsize {\tiny $x_5$}}
\put(2.1,.7){\scriptsize {\tiny $x_6$}}
\put(3.1,.7){\scriptsize {\tiny $x_4$}}

\put(-.5,1.7){\scriptsize {\tiny $x_7$}}
\put(1.1,1.7){\scriptsize {\tiny $x_8$}}
\put(2.1,1.7){\scriptsize {\tiny $x_9$}}
\put(3.1,1.7){\scriptsize {\tiny $x_7$}}

\put(-.7,2.7){\scriptsize {\tiny $x_{10}$}}
\put(1.1,2.7){\scriptsize {\tiny $x_{11}$}}
\put(2.1,2.7){\scriptsize {\tiny $x_{12}$}}
\put(3.1,2.7){\scriptsize {\tiny $x_{10}$}}

\put(-.4,4.1){\scriptsize {\tiny $x_1$}}
\put(.9,4.1){\scriptsize {\tiny $x_2$}}
\put(1.9,4.1){\scriptsize {\tiny $x_3$}}
\put(3,4.1){\scriptsize {\tiny $x_1$}}

\put(-.8,-1.2){\scriptsize {\tiny {\bf Figure 5.4.3:} $M(3,4,0)$}} 

\end{picture}

\begin{picture}(0,0)(-16,15)
\setlength{\unitlength}{5.8mm}

\drawpolygon(5,0)(8,0)(8,4)(5,4)


\drawline[AHnb=0](6,0)(6,4)
\drawline[AHnb=0](7,0)(7,4)


\drawline[AHnb=0](5,1)(8,1)
\drawline[AHnb=0](5,2)(8,2)
\drawline[AHnb=0](5,3)(8,3)


\drawline[AHnb=0](5,3)(6,4)
\drawline[AHnb=0](6,3)(7,4)
\drawline[AHnb=0](7,3)(8,4)

\put(4.8,-.4){\scriptsize {\tiny $y_1$}}
\put(6,-.4){\scriptsize {\tiny $y_2$}}
\put(7,-.4){\scriptsize {\tiny $y_3$}}
\put(8,-.4){\scriptsize {\tiny $y_1$}}

\put(4.5,.7){\scriptsize {\tiny $y_4$}}
\put(6.1,.7){\scriptsize {\tiny $y_5$}}
\put(7.1,.7){\scriptsize {\tiny $y_6$}}
\put(8.2,.7){\scriptsize {\tiny $y_4$}}

\put(4.5,1.7){\scriptsize {\tiny $y_7$}}
\put(6.05,1.7){\scriptsize {\tiny $y_8$}}
\put(7.05,1.7){\scriptsize {\tiny $y_9$}}
\put(8.1,1.7){\scriptsize {\tiny $y_7$}}

\put(4.3,2.7){\scriptsize {\tiny $y_{10}$}}
\put(6.1,2.7){\scriptsize {\tiny $y_{11}$}}
\put(7.1,2.7){\scriptsize {\tiny $y_{12}$}}
\put(8.1,2.7){\scriptsize {\tiny $y_{10}$}}

\put(4.5,4.12){\scriptsize {\tiny $y_2$}}
\put(5.8,4.12){\scriptsize {\tiny $y_3$}}
\put(6.8,4.12){\scriptsize {\tiny $y_1$}}
\put(8,4.12){\scriptsize {\tiny $y_2$}}
\put(4.5,-1.2){\scriptsize {\tiny {\bf Figure 5.4.4:} $M(3,4,1)$ }} 

\end{picture}

\begin{picture}(0,0)(-56,10)
\setlength{\unitlength}{5.8mm}

\drawpolygon(5,0)(8,0)(8,4)(5,4)


\drawline[AHnb=0](6,0)(6,4)
\drawline[AHnb=0](7,0)(7,4)


\drawline[AHnb=0](5,1)(8,1)
\drawline[AHnb=0](5,2)(8,2)
\drawline[AHnb=0](5,3)(8,3)

\drawline[AHnb=0](5,3)(6,4)
\drawline[AHnb=0](6,3)(7,4)
\drawline[AHnb=0](7,3)(8,4)

\put(4.8,-.4){\scriptsize {\tiny $z_1$}}
\put(6,-.4){\scriptsize {\tiny $z_2$}}
\put(7,-.4){\scriptsize {\tiny $z_3$}}
\put(8,-.4){\scriptsize {\tiny $z_1$}}

\put(4.5,.7){\scriptsize {\tiny $z_4$}}
\put(6.1,.7){\scriptsize {\tiny $z_5$}}
\put(7.1,.7){\scriptsize {\tiny $z_6$}}
\put(8.1,.7){\scriptsize {\tiny $z_4$}}

\put(4.5,1.7){\scriptsize {\tiny $z_7$}}
\put(6.05,1.7){\scriptsize {\tiny $z_8$}}
\put(7.05,1.7){\scriptsize {\tiny $z_9$}}
\put(8.1,1.7){\scriptsize {\tiny $z_7$}}

\put(4.3,2.7){\scriptsize {\tiny $z_{10}$}}
\put(6.1,2.7){\scriptsize {\tiny $z_{11}$}}
\put(7.05,2.7){\scriptsize {\tiny $z_{12}$}}
\put(8.1,2.7){\scriptsize {\tiny $z_{10}$}}

\put(4.5,4.15){\scriptsize {\tiny $z_3$}}
\put(5.8,4.15){\scriptsize {\tiny $z_1$}}
\put(6.9,4.15){\scriptsize {\tiny $z_2$}}
\put(8.1,4.15){\scriptsize {\tiny $z_3$}}

\put(4.5,-1.2){\scriptsize {\tiny {\bf Figure 5.4.5:} $M(3,4,2)$ }} 

\end{picture}

\begin{picture}(0,0)(-94,5)
\setlength{\unitlength}{5.8mm}

\drawpolygon(5,0)(8,0)(8,4)(5,4)


\drawline[AHnb=0](6,0)(6,4)
\drawline[AHnb=0](7,0)(7,4)


\drawline[AHnb=0](5,1)(8,1)
\drawline[AHnb=0](5,2)(8,2)
\drawline[AHnb=0](5,3)(8,3)


\drawline[AHnb=0](5,3)(6,4)
\drawline[AHnb=0](6,3)(7,4)
\drawline[AHnb=0](7,3)(8,4)

\put(4.8,-.4){\scriptsize {\tiny $x_1$}}
\put(6,-.4){\scriptsize {\tiny $x_3$}}
\put(7,-.4){\scriptsize {\tiny $x_2$}}
\put(8,-.4){\scriptsize {\tiny $x_1$}}

\put(4.5,.7){\scriptsize {\tiny $x_4$}}
\put(6.1,.7){\scriptsize {\tiny $x_6$}}
\put(7.1,.7){\scriptsize {\tiny $x_5$}}
\put(8.1,.7){\scriptsize {\tiny $x_4$}}

\put(4.5,1.7){\scriptsize {\tiny $x_7$}}
\put(6.05,1.7){\scriptsize {\tiny $x_9$}}
\put(7.05,1.7){\scriptsize {\tiny $x_8$}}
\put(8.1,1.7){\scriptsize {\tiny $x_7$}}

\put(4.3,2.7){\scriptsize {\tiny $x_{10}$}}
\put(6.1,2.7){\scriptsize {\tiny $x_{12}$}}
\put(7.05,2.7){\scriptsize {\tiny $x_{11}$}}
\put(8.1,2.7){\scriptsize {\tiny $x_{10}$}}

\put(4.5,4.15){\scriptsize {\tiny $x_2$}}
\put(5.8,4.15){\scriptsize {\tiny $x_1$}}
\put(6.9,4.15){\scriptsize {\tiny $x_3$}}
\put(8.1,4.15){\scriptsize {\tiny $x_2$}}

\put(4.5,-1.2){\scriptsize {\tiny {\bf Figure 5.4.6:} $M(3,4,2)$ }} 

\end{picture}



\vspace{1.5cm}

\begin{center}
	
	\noindent \textbf{Table \ref{table:4}} : DSEMs of type $T_4$ on the torus for $|V(M)| \leq 24$
	
	\renewcommand{\arraystretch}{1.1}
	\begin{tabular}{ |p{1.2cm}|p{4cm}|p{3.8cm}|p{3cm}|}
		
		\hline
		$|V(M)|$ & Isomorphic classes & Length of cycles & No of maps\\
		\hline
		12 &  $M(3,4,0)$, $M(3,4,2)$& $( 3, \{ 4, 12 \}, 4 )$ & 2\\
		
		\cline{2-3}
		& $M(3,4,1)$ & $( 3, \{ 12, 12 \}, 5 )$ & \\	
		
		\cline{2-3}
		
		\hline	
		
		16 &  $M(4,4,0)$, $M(4,3,3)$ & $( 4, \{ 4, 16 \}, 4 )$& 2\\
		
		\cline{2-3}
		
		& $M(4,4,1)$, $M(4,4,2)$ & $( 4, \{ 8, 16 \}, 5 )$ & \\
		
		\cline{2-3}
		
		\hline
		20 & $M(5,4,0)$, $M(5,4,4)$ & $( 5, \{ 4, 20 \}, 4 )$ &3\\
		
		\cline{2-3}

		& $M(5,4,1)$, $M(5,4,3)$  &$( 5, \{ 20, 20 \}, 5 )$ &\\
		
		\cline{2-3}
		
		& $M(5,4,2)$ & $( 5, \{ 20, 20 \}, 6 )$  &\\
		
		\cline{2-3}

		\hline	
		24 & $M(6,4,0)$, $M(6,4,5)$ & $( 6, \{ 4, 24 \}, 4 )$ & 5\\
		
		\cline{2-3}
		
		&$M(6,4,1)$, $M(6,4,4)$ &$( 6, \{ 12, 24 \}, 5 )$&\\
		
		\cline{2-3}
		
		& $M(6,4,2)$, $M(6,4,3)$ &$( 6, \{ 8, 12 \}, 6)$&\\
		
		\cline{2-3}
		
		& $M(3,8,0)$, $M(3,8,1)$ &$( 3, \{ 8, 24 \}, 8 )$ & \\
		
		\cline{2-3}
		
		& $M(3,8,2)$ & $( 3, \{ 24, 24 \}, 10 )$ &  \\
		\hline
		
	\end{tabular}
	\label{table:4}
\end{center}

\bigskip
\noindent{\bf Proof of Theorem \ref{t2} and Corollary \ref{co2}:} Follows from the results and tables given in Sections \ref{s4.1}-\ref{s4.4}. \hfill$\Box$

\vspace{.4cm}

\noindent{\bf Remark:} As, we know that out of infinitely many types of semi-equivelar maps on the plane, only eleven types (corresponding to the Archimedean tilings of the plane) give the same type  semi-equivelar maps on the torus. In case of doubly semi-equivelar maps, one can construct several types (other than the 22 types discussed here) on the torus by stacking suitably the faces of semi-equivelar maps on the torus. Thus a natural questions occurs here: are there exist finitely many types of doubly semi-equivelar maps on the torus?

\section{Acknowledgement} \label{s4.8} The first author is thankful to the Ministry of Human Resource Development, New Delhi (India), for financial support. The second author expresses his thanks to IIIT Allahabad for providing the facility and resources to carry this research work.

\newpage
\end{document}